\documentclass[10pt]{elsarticle}
\usepackage{lineno,hyperref}
\usepackage{booktabs}
\usepackage{multirow}
\usepackage{amsthm}
\usepackage{hyperref}
\usepackage{crossreftools}
\usepackage{physics}
\usepackage{graphicx}
\usepackage{subfigure}
\usepackage{float}
\usepackage{bm}
\usepackage{lipsum}
\usepackage{listings} 
\usepackage{color}
\usepackage{amsfonts,enumerate,amsmath,amssymb}
\def\qed{\strut\hfill $\Box$}
\newtheorem{thm}{Theorem}[section]

\newtheorem{prop}[thm]{Proposition}
\newtheorem{lem}[thm]{Lemma}

\newtheorem{rem}[thm]{Remark}
\newtheorem{defn}[thm]{Definition}

\newcommand{\thmref}[1]{Theorem~{\rm \ref{#1}}}
\newcommand{\lemref}[1]{Lemma~{\rm \ref{#1}}}
\newcommand{\propref}[1]{Proposition~{\rm \ref{#1}}}

\newcommand{\remref}[1]{Remark~{\rm \ref{#1}}}

\def\para#1{\vskip .4\baselineskip\noindent{\bf #1}}
\newcommand{\vertiii}[1]{{\left\vert\kern-0.25ex\left\vert\kern-0.25ex\left\vert #1 
\right\vert\kern-0.25ex\right\vert\kern-0.25ex\right\vert}}
\def\para#1{\vskip .4\baselineskip\noindent{\bf #1}}
\numberwithin{equation}{section}
\begin{document}
\begin{frontmatter} 
\title{Large and moderate deviation for rough slow-fast system under level 3 geometric rough path}
\author[mymainaddress,mysecondaddress]{Xiaoyu Yang}
\ead{yangxiaoyu@yahoo.com}

\author[mymainaddress,mythirdaddress]{Yong Xu\corref{mycorrespondingauthor}}\cortext[mycorrespondingauthor]{Corresponding author}\ead{hsux3@nwpu.edu.cn}

\address[mymainaddress]{School of Mathematics and Statistics, Northwestern Polytechnical University, Xi'an, 710072, China}
\address[mysecondaddress]{Faculty of Mathematics, Kyushu University, Fukuoka, 8190395, Japan}
\address[mythirdaddress]{MOE Key Laboratory of Complexity Science in Aerospace, Northwestern Polytechnical University, Xi'an, 710072, China}

\begin{abstract}
This work is to give the large deviation principle (LDP) and moderate deviation principle (MDP) for a slow-fast system driven by the mixed fractional Brownian motion (FBM) $(1/4,1/3)$ via the level-3 geometric rough path (RP).  Firstly, we provide a different way to lift the translation of mixed FBM in the Cameron-Martin direction to the geometric RP. The LDP turns to the weak convergence of the controlled system by averaging the fast one. Lacking the invariant measure for the controlled fast one, a ``replaced process” whose limiting measure could be decoupled from the control term and owes the exponential ergodicity is constructed. Here, the stability under level-3 RPs is established under more elaborate estimates. Besides, different from LDP we give the MDP, that requires more intricate bounded estimates of the deviation component. The MDP for level-2 RP system is recovered as a special case.

\vskip 0.08in 
\noindent{\bf AMS Math Classification.}
{60F10, 60G15, 60H10}.
\vskip 0.08in 
\end{abstract} 
\end{frontmatter}

\tableofcontents

\section{Introduction}\label{sec-1}
Next, we consider the slow-fast system with mixed FBM in the time interval $[0,T]$:
\begin{eqnarray}\label{1-4}
\left
\{
\begin{array}{ll}
dX^{\varepsilon, \delta}_t = f(X^{\varepsilon, \delta}_t, Y^{\varepsilon, \delta}_t)dt + \sqrt \varepsilon \sigma( X^{\varepsilon, \delta}_t)db^H_{t},\\
dY^{\varepsilon, \delta}_t =\frac{1}{\delta} F( X^{\varepsilon, \delta}_t, Y^{\varepsilon, \delta}_t)dt + \frac{1}{{\sqrt \delta}}G( X^{\varepsilon, \delta}_t, Y^{\varepsilon, \delta}_t)dw_{t},\\
(X^{\varepsilon, \delta}_0, Y^{\varepsilon, \delta}_0)=(x, y)\in \mathbb{R}^{{m}}\times \mathbb{R}^{{n}}.
\end{array}
\right.
\end{eqnarray}
Here, $(b^H,w)$ is the mixed FBM with Hurst parameter $H\in (\frac{1}{4},\frac{1}{3})$. $\varepsilon, \delta\in(0,1]$ are small parameters which are assumed that $\delta=o(\varepsilon)$. The coefficients $f:\mathbb{R}^{{m}}\times \mathbb{R}^{{n}}\to\mathbb{R}^{{m}}$, $F:\mathbb{R}^{{m}}\times \mathbb{R}^{{n}}\to\mathbb{R}^{{n}}$, {$\sigma: \mathbb{R}^{{m}}\to \mathbb{R}^{{m}\times{d}}$ and $G: \mathbb{R}^{{m}}\times \mathbb{R}^{{n}}\to \mathbb{R}^{{n}\times{e}} $} are nonlinear functions, which will be assumed ``smooth" enough in the Section 5. Set $\Omega=\mathcal{C}_0([0,T]: \mathbb{R}^{d+e})$. The complete probability space $(\Omega, \mathcal{F}, \mathbb{P})$ supports $b^H$ and $w$ exists independently, where $\mathbb{P}$ is the unique probability measure on $\Omega$ and $\mathcal{F}=\mathcal{B}(\mathcal{C}_0([0,T]: \mathbb{R}^{d+e}))$ is the $\mathbb{P}$-completion of the Borel $\sigma$-field. Then, we consider the canonical filtration given by $\{\mathcal{F}_t^H: t \in[0,T]\}$, where $\mathcal{F}_t^H=\sigma\{(b_s^H,w_s): 0 \leq s \leq t\} \vee \mathcal{N}$ and $\mathcal{N}$ is the set of the $\mathbb{P}$-negligible events. This kind of slow-fast systems with different time scales are widely applied in diverse areas, such as feedback control of mobile robot and so on \cite{2008Pavliotis,Kuehn2015}. 

The FBM is self-similar and possesses long-range dependence \cite{2008Biagini,2008Mishura}. However, due to the property whereby the Hurst parameter $H$ describes the raggedness of the resultant motion with higher values leading to smoother motion, the definition of the solution mapping to the stochastic system driven by FBM depends on the choice of Hurst parameter $H$. When $H\in(1/2,1)$, the solution to the stochastic differential equations (SDE) driven by mixed FBM can be defined using the generalized Riemann–Stieltjes integral \cite{Guerra}. When $H < 1/2$ and the driver path is an additive FBM, the issue can be addressed by using the generalized Riemann–Stieltjes integral. However, in the case of multiplicative noise, it will be more complicated and the generalized Riemann–Stieltjes integral will be invalid. The rough path (RP) theory, a type of pathwise integral theory that does not require filtration, Markov or martingale theory, is efficient in solving this ``rough" case (see references \cite{1998Lyons, 2002Lyons,2010Friz,2020Friz}. Due to the result that the FBM and mixed FBM could be lifted to geometric rough path (GRP) when the Hurst parameter $H\in (1/4,1/2)$ \cite{2010Friz}, the solution to the RDE could be well-defined in the rough integral. 

According to the averaging principle, under suitable conditions, the effective dynamics of a slow system can be described by an averaged system as the time-scale parameter tends to zero \cite{2021Pei}, that is
\begin{eqnarray}\label{4-1}
d\bar X_t = \bar f(\bar X_t)dt,
\end{eqnarray}
with $X_0=x\in \mathbb{R}^{{m}}$. 
where $\bar{f}(x)=\int_{\mathbb{R}^{n}}f(x, { y})\mu^{x}(d{ y})$ for $x\in\mathbb{R}^m$. Here, $\mu^x$ is the invariant probability measure of the fast variable and the ``frozen" slow one, as described in references. When $H\in (1/3,1/2)$, the averaging principle for slow-fast SDEs or SPDEs are given in \cite{2023Pei, 2022Inahama,Pei2025}. When the sample paths of FBM are rougher, that is $H\in (1/4,1/3)$, the averaging principle for the slow-fast system is also constructed \cite{2025Inahama}
The averaging principle can be viewed as a form of the law of large numbers. 

Then, our aims are to study the fluctuation of $X^{\varepsilon, \delta}$ around the limit $X$. Firstly, the deviation component between $X^{\varepsilon, \delta}$ and $\bar X$ is defined as below, 
\begin{eqnarray}\label{02}
Z^{\varepsilon,\delta}_t = \frac{X^{\varepsilon, \delta}_t -\bar X_t }{\sqrt{\varepsilon}h(\varepsilon)},\quad Z^{\varepsilon}_0=0.
\end{eqnarray}
Here, $h:(0,1]\to (0,\infty)$ is continuous function. The case $h(\varepsilon)\equiv1$ corresponds to the central limit theorem (CLT), while the case $h(\varepsilon)=\frac{1}{\sqrt{\varepsilon}}$ corresponds to the large deviation principle (LDP), the case $\lim_{\varepsilon \to 0}\sqrt{\varepsilon}h(\varepsilon)=0$ is called the moderate deviation principle (MDP). The MDP are always considered as a bridge between the CLT and the LDP. Specially, it is equivalent between the MDP for $\{X^{\varepsilon,\delta}\}_{\varepsilon,\delta\in(0,1]}$ and the LDP for $\{Z^{\varepsilon,\delta}\}_{\varepsilon,\delta\in(0,1]}$ but with more complicated estimates for the deviation component.

The family
$\{X^{\varepsilon,\delta}\}_{\varepsilon\in(0,1]}\in \mathcal{C}^{\alpha}\left([0,T], \mathbb{R}^m\right)$ is said to satisfy an LDP on $ \mathcal{C}^{\alpha}\left([0,T], \mathbb{R}^m\right)$ with speed $a(\varepsilon)=\varepsilon$ and a good 
rate function $I: \mathcal{C}^{\alpha}\left([0,T], \mathbb{R}^m\right)\rightarrow [0, \infty]$ ($1-H <\alpha <1/2$) if the following two conditions hold:
{\begin{itemize}
\item For each closed subset $F$ of $C^{\alpha}\left([0,T], \mathbb{R}^m\right)$,
$$\limsup _{\varepsilon \rightarrow 0} a(\varepsilon) \log \mathbb{P}\big(X^{{\varepsilon,\delta}} \in F\big) \leqslant-\inf _{X \in F} I(X).$$
\item For each open subset $G$ of $C^{\alpha}\left([0,T], \mathbb{R}^m\right)$,
$$\liminf _{\varepsilon \rightarrow 0} a(\varepsilon) \log \mathbb{P}\big(X^{{\varepsilon,\delta}} \in G\big) \geqslant-\inf _{X \in G} I(X).$$
\end{itemize}}
The above LDP is the known Freidlin-Wentzell LDP \cite{1984Random}, which is based on the so-called Freidlin-Wentzell exponential estimates \cite{1994Peszat}. When $a(\varepsilon)=1/\varepsilon^2$, it is the MDP. Since it is equivalent between the MDP for the slow component and LDP for the deviation component, LDP will be introduced first. Then, the variational representation and weak convergence method are proposed, the advantage of which is that it avoids the complicated exponential estimates \cite{1998Dupuis,2011Variational,2011Dupuis,BP_book}. After that, the other method, based on nonlinear semigroups and viscosity solution theory, has also been proposed \cite{FFK, FK}.

Our focus is on the asymptotic behavior of the slow component, $X^{\varepsilon,\delta}$ around its averaged system, i.e. the LDP and MDP. If $f(x,y)=f(x)$ in the above slow-fast system \eqref{1-4}, then it will degenerate into a single-times cale system. When $H<1/2$, the LDP for this type of system has been well established \cite{2015Inahama,2007Friz, 2010Friz,2006Millet}. Of course, when the driver path is a geometric RP which is lifted by BM, the LDP also was proved holds \cite{2002Ledoux,2013Inahama}. The LDP and MDP also could be obtained by applying the Lyons' continuity theorem and the property of the product measure between two measure families \cite{2022Yang,Emanuela,2024Inahama}. However, the above method is invalid in the multi-time scale system due to the perturbation from the fast varying process. Thanks to the variational representation Wiener space, the weak convergence method is possible to slove the LDP issue for the slow-fast system. 
When $H=1/2$, the aforementioned slow-fast system degenerates into one driven by standard Brownian motion (BM). The LDP has been established by applying the weak convergence method \cite{2012Dupuis}. When the Hurst parameter $H>1/2$, thanks to the variational representation for the random functional in abstract Wiener space \cite{2009Zhang}, the weak convergence method for the LDP of stochastic system driven by FBM was obstructed \cite{2020Budhiraja}, which will be possible to slove the LDP issue for the slow-fast system. 

Subsequently, the LDP for slow-fast system with BM were established by combining weak convergence and Khasminskii averaging technique \cite{2020Large}, or the viable pair method \cite{2017Morse,2019Hu}. Based on this, the LDP for other types of slow-fast (partial) systems was also established \cite{Liu2020, 2023Hong}. 
The situation will be more complicated when the Hurst parameter $H\neq 1/2$. When the Hurst parameter is in the range $H\in(1/2,1)$, the LDP for the aforementioned slow-fast system driven by mixed FBM was established via variational representation for FBM and the weak convergence method \cite{2023Inahama}. 
Furthermore, it was demonstrated that the LDP can be obtained by combining the weak convergence and viable pair method \cite{2021Bourguin,2025Gailus}, which can also be extended to the MDP \cite{2023Gasteratos}.

When the Hurst parameter satisfies $H\in(1/3,1/2)$, the LDP has been established via the weak convergence method combined with the Khasminskii averaging technique \cite{2025Yang}. To date, however, no LDP or MDP results are available for the above slow-fast system driven by mixed FBM with $H\in(1/4,1/3)$. Unlike the previous level-2 driver, the new driving path arises from an GRP derived from the mixed FBM with $H\in(1/4,1/3)$. The first part is dedicated to the LDP to the slow-fast system with a level-3 random GRP, which is lifted from the mixed FBM with Hurst parameter $H\in(1/4,1/3)$. To this end, we again employ the variational representation and the weak convergence method. Hence, we first to show that the Cameron--Martin translation of the mixed FBM can still be lifted to a GRP via the anisotropic RP. The anisotropic RP was developed in \cite{2016Gyurko}. In fact, this approach drops the Gaussianity assumption on the first element of the mixed process, which is different from \cite[Lemma 15.58]{2010Friz}. However, the presence of third-level RPs significantly complicates the weak convergence analysis, necessitating more delicate controlled RP estimates. Ultimately, we prove that the controlled slow component converges weakly to the solution of the skeleton equation. Since the controlled fast variable admits no invariant measure, we construct an auxiliary fast process that possesses an invariant measure and enjoys exponential ergodicity, which serves as a surrogate for the controlled fast variable. Exploiting the continuity of the solution map, we then establish the desired weak convergence result. In this setting, the solution to the skeleton equation remains well-defined in the Young sense within a variational framework. Finally, we present the MDP, which can be obtained under certain additional assumptions. Unlike the LDP, the MDP requires a more delicate and involved analysis, together with certain boundedness estimates for the deviation process.

We will now outline this paper. In Section 2, we introduce some notation and preliminaries, as well as providing some estimates on the RDE driven by level 3 GRP. Section 3 provides some necessary results in GRP. Section 4 is devoted to the introduction for the slow-fast system under mixed FBM. Section 5 gives the large and moderate deviation results to slow-fast systems with level 3 random GRP.
Throughout this paper, $c$, $C$, $c_1$, $C_1$, $\cdots$ denote certain positive constants that may vary from line to line. The dependence of the constant on parameters will be explicit if necessary. {$\mathbb{N}=\{1,2, \ldots\}$} and time horizon $T >0$. Let $\lfloor a \rfloor$ denote the integer part for any positive number $a$. $\hookrightarrow$ stands for the continuous embedding between different spaces.

\section{Rough differential equations driven by level 3 rough paths}\label{sec-2}
\textbf{Notations}
Firstly, we introduce the notations that will be used throughout the paper.
Let $[a,b]\subset[0,T]$ and $\Delta_{[a,b]}:=\{(s,t)\in \mathbb{R}^2|a\leq s\leq t\leq b\}$. We write $\Delta_{T}$ simply when $[a,b]=[0,T]$. Denote $\nabla$ be the standard gradient on a Euclidean space. We assume $\mathcal{V}$, $\mathcal{V}_i, i\in \mathbb{N}$ and $\mathcal{W}$ are both Euclidean spaces. We denote the truncated tensor algebra of degree $k$ over $\mathcal{V}$ by $T^{k}(\mathcal{V})=\mathbb{R} \oplus \mathcal{V} \oplus \mathcal{V}^{\otimes 2}\oplus \cdots \mathcal{V}^{\otimes k}$ and $\otimes$ is the tensor product. All notation and results in this section are deterministic; no probabilistic structure is involved.
\begin{itemize}
\item \textbf{(Continuous space)} The set of continuous functions is denoted by $\varphi:[a,b]\to \mathcal{V}$. With the norm $\|\varphi\|_{\infty}=\sup _{t \in[a, b]}|\varphi_t|<\infty$, it is a Banach space. The set of continuous functions starts from $0$ is denoted by $\mathcal{C}_0([a,b],\mathcal{V})$.
\item \textbf{(H\"older continuous space)}
Set $\eta \in(0,1]$. For a path $\varphi:[a,b]\mapsto\mathcal{V}$, define the $\eta-$H\"older semi-norm by
$$
\|\varphi\|_{\eta-\operatorname{hld},[a,b]}:=\sup _{a \leq s<t \leq b} \frac{|\varphi_t-\varphi_s|}{(t-s)^\eta}<\infty.
$$
The $\varphi$ is called $\eta$-H\"older continuous. 
The set of all $\eta$-H\"older continuous path is denoted by $\mathcal{C}^{\eta-\operatorname{hld}}([a,b],\mathcal{V})$.
Define its Banach norm by $ |\varphi_a|_{\mathcal{V}}+
\|\varphi\|_{\eta-\operatorname{hld},[a,b]}$. The set of $\eta$-H\"older continuous functions starts from $0$ is denoted by $\mathcal{C}_0^{\eta-\operatorname{hld}}([a,b],\mathcal{V})$.
\item For a continuous map $\psi:\Delta_{[a,b]}\to \mathcal{V}$, we set 
$$ \|\psi\|_{\eta-\operatorname{hld},[a,b]}:=\sup _{a \leq s<t \leq b} \frac{|\psi_t-\psi_s|}{(t-s)^\eta}.$$
the set of above $\psi$ with $\|\psi\|_{\eta-\operatorname{hld},[a,b]}<\infty$ is denoted by $\mathcal{C}_2^{\eta-\operatorname{hld}}([a,b],\mathcal{V})$. It is a Banach space with the norm $\|\psi\|_{\eta-\operatorname{hld},[a,b]}$.
\item \textbf{(Variation space)} Set $1\leq p<\infty $. Denote $\mathcal{C}^{p-\operatorname{var}}\left([a,b],\mathcal{V}\right)=\left\{\varphi \in \mathcal{C}\left([a,b],\mathcal{V}\right):\|\varphi\|_{p \text {-var }}<\infty\right\}$ with the usual $p$-variation semi-norm $\|\varphi\|_{p \text {-var }}$. The set of $p$-variation functions starts from $0$ is denoted by $\mathcal{C}_0^{p \text{-var}}([a,b],\mathcal{V})$.

\item \textbf{(Besov space)} For continuous function $\phi:[a,b] \rightarrow \mathcal{V}$ and $\delta \in(0,1) $ and $p \in(1, \infty)$, the Besov space is denoted by $W^{\delta, p}\left([a,b],\mathcal{V}\right)$ equipped with the norm as follows:
\begin{equation}\label{2-2}
\|\phi\|_{W^{\delta, p}}=\|\phi\|_{L^p}+\left(\iint_{[a,b]^2} \frac{\left|\phi_t-\phi_s\right|^p}{|t-s|^{1+\delta p}} d s d t\right)^{1 / p} <\infty.
\end{equation}
When $\eta'=\delta-1/p>0$, the continuous imbedding $W^{\delta, p}([a,b],\mathcal{V}) \hookrightarrow \mathcal{C}^{\eta'-\operatorname{hld}}([a,b],\mathcal{V})$ holds and more could refer to \cite[Theorem 2]{2006Friz} for further detail.
\item \textbf{{($C^k$-function and $C^k_b$-function)}}
Let $U\subset \mathcal{V}$ be an open set. For $k\in \mathbb{N}$, denote $C^k(U,\mathcal{W})$ by the set of $C^k$-functions from $U$ to $\mathcal{W}$. The set of $C^k$-bounded functions whose derivatives up to order $k-$ are also bounded is denoted by $C_b^k(U,\mathcal{W})$. The space $C_b^k(U,\mathcal{W})$ is a Banach space with the norm $\|\varphi\|_{C_b^k}:=\sum_{i=0}^{k}\|\nabla^i\varphi\|_\infty<\infty$.
\item $L(\mathcal{W},\mathcal{V})$ denotes the set of bounded linear maps from $\mathcal{W}$ to $\mathcal{V}$. We set the $k$-bounded linear maps from $\mathcal{V}_1\times \cdots\times \mathcal{V}_k$ to $\mathcal{W}$ by $L^{(k)}(\mathcal{V}_1,\cdots, \mathcal{V}_k; \mathcal{W})\cong L( \mathcal{V}^{\otimes k}, \mathcal{W})$.
\item \textbf{{(Signature of path)}} For any continuous path $x\in\mathcal{C}_0^{\eta-\operatorname{hld}}([a,b],\mathcal{V})$, one can construct its $m$-layer roughness path, denoted $S_m:=(S^1, S^2, \cdots, S^m)$, where $S^k$ is computed as follows for any $1\le k\le m$.
\begin{eqnarray}\label{exm3}
S_{s, t}^k=\int_{s<t_1<\cdots<t_k<t} {~d} x_{t_1} \otimes \cdots \otimes {~d} x_{t_k},
\end{eqnarray}
then
\begin{eqnarray}\label{exm4}
S_m(x)_{s, t}= S_m(x)_{s, u} \otimes S_m(x)_{u,t}.
\end{eqnarray}
\end{itemize}

\subsection{Rough Path and Geometric Rough Path(GRP)}
In this subsection, let $\alpha\in (1/4,1/3]$ and $k:=\lfloor 1 / \alpha\rfloor$. 
Now, we recall the definition of rough path (RP) and geometric RP (GRP).
\begin{defn}(\textbf{Rough Path})\label{RPS}
A continuous map $\mathbf{X}=\big(1, X^{1}, X^{2},\cdots, X^{k}\big): \Delta \rightarrow T^{k}(\mathcal{V}) $ is said to be a $\mathcal{V}$-valued RP of roughness $k$ if it satisfies the following conditions: \textbf{(A)}: For any $s \leq u \leq t$, $\mathbf{X}_{s, t}=\mathbf{X}_{s, u} \otimes \mathbf{X}_{u, t}$. \textbf{(B)}: $\|X^{i}\|_{i\alpha \mathrm{-hld}}<\infty$ for $i\in \{1,\cdots,k\}$.
\end{defn}
The $0$-th element $1$ is obviously omitted and the RP is denoted by $\mathbf{X}=\left(X^{1}, \cdots X^{k}\right)$. Below, we set $\vertiii{\mathbf{X}}_{\alpha-\operatorname{hld}}:=\sum_{i=1}^{k}\|X^{i}\|_{i\alpha \mathrm{-hld}}$. The set of all $\mathcal{V}-$valued RPs with $1/4<\alpha<1/3$ is denoted by $\Omega_{\alpha}(\mathcal{V})$. Equipped with the $\alpha$-H\"older distance, the above RP space is a complete space. It is directly to verify that $\Omega_\alpha(\mathcal{V}) \hookrightarrow \Omega_\beta(\mathcal{V})$ for $1/4<\beta \leq \alpha \leq 1/3$. 
The distance between two different RPs $\mathbf{X}\in \Omega_\alpha(\mathcal{V})$ and $\tilde{\mathbf{X}}\in \Omega_\alpha(\mathcal{V})$ is defined by
$\rho_\alpha(\mathbf{X}, \tilde{\mathbf{X}}):=\sum_{i=1}^{k}\|X^i-\tilde{X}^i\|_{i\alpha \mathrm{-hld}}$. 

Now we recall the definition of geometric rough path.
\begin{defn}(\textbf{Geometric Rough Path})\label{GRP}
A geometric rough path (GRP) is the closure of the $S_k(\mathcal{C}_0^1(\mathcal{V}))$ under the $\alpha$-H\"older topology. 
\end{defn}
The set of all GRP is denoted by $G\Omega_\alpha(\mathcal{V})$, which is a complete and separable space.
A RP $\mathbf{X}=\big(1, X^{1}, X^{2},\cdots, X^{k}\big): \Delta \rightarrow T^{k}(\mathcal{V})\in G\Omega_\alpha(\mathcal{V})$ is called geometric rough path (GRP) if it satisfies the algebraic property which is also called Shuffle relation.
Then we explain the Shuffle relation with $\mathcal{V}=\mathbf{R}^m$ and $k=3$, that is, for all $p,q,r=\{1,\cdots, m\}$ and $\Delta_{T}$,
\begin{eqnarray*}\label{Shuffle}
X_{s,t}^{1,p}X_{s,t}^{1,p}&=&X_{s,t}^{2,pq}+X_{s,t}^{2,qp}\cr
X_{s,t}^{1,p}X_{s,t}^{2,qr}&=&X_{s,t}^{3,pqr}+X_{s,t}^{3,qpr}+X_{s,t}^{3,qrp}.
\end{eqnarray*}

\subsection{Controlled RP with reference level 3 RP}
In this subsection, we mainly introduce the controlled RP driven by reference RP with $k=3$ (The situation of $k = 2$ has been already presented comprehensively in \cite{2020Friz}).

Firstly, we recall the definition of controlled RP with respect to the reference RP. It says that $(Y, Y^{\dagger}, Y^{\dagger\dagger}, Y^{\sharp}, Y^{\sharp\sharp})$
is a $\mathcal{W}$-valued controlled RP with respect to $\mathbf{X}=\left(X^{1}, X^{2}, X^{3}\right)\in \Omega_{\alpha}(\mathcal{V})$ if it satisfies the following conditions:
$$
Y_t-Y_s=Y_s^{\dagger} X_{s, t}^1+Y_t^{\dagger\dagger} X_{s, t}^2+Y^{\sharp}_{s, t}$$ 
$$Y^{\dagger}_t-Y^{\dagger}_s=Y_s^{\dagger\dagger} X_{s, t}^1+Y^{\sharp\sharp}_{s, t}$$
for all $(s, t) \in \triangle_{[a, b]} $ and
\begin{eqnarray*}
(Y, Y^{\dagger}, Y^{\dagger\dagger}, Y^{\sharp}, Y^{\sharp\sharp})&\in& \mathcal{C}^{\alpha-\operatorname{hld}}([a, b], \mathcal{W}) \times \mathcal{C}^{\alpha-\operatorname{hld}}([a, b], L(\mathcal{V}, \mathcal{W}))\times \mathcal{C}^{\alpha-\operatorname{hld}}([a, b], L(\mathcal{V}, L(\mathcal{V}, \mathcal{W}))) \cr &&\times \mathcal{C}_2^{3\alpha-\operatorname{hld}}([a, b], \mathcal{W})\times \mathcal{C}_2^{2\alpha-\operatorname{hld}}([a, b], L(\mathcal{V}, \mathcal{W}))
\end{eqnarray*}
Let $\mathcal{Q}_{\mathbf{X}}^\alpha([a, b], \mathcal{W})$ stand for the set of all above controlled RPs and always denote the controlled RP $(Y, Y^{\dagger}, Y^{\dagger\dagger})$ for simplicity. Denote the semi-norm of controlled RP $(Y, Y^{\dagger}, Y^{\dagger\dagger}, Y^{\sharp}, Y^{\sharp\sharp})\in \mathcal{Q}_{\mathbf{X}}^\alpha([a, b], \mathcal{W})$ by 
$$\|Y, Y^{\dagger}, Y^{\dagger\dagger}\|_{\mathcal{Q}_{\mathbf{X}}^\alpha,[a, b]}=\|Y^{\dagger\dagger}\|_{{\alpha-\operatorname{hld}},[a, b]}
+\|Y^{\sharp}\|_{{2\alpha-\operatorname{hld}},[a, b]}
+\|Y^{\sharp\sharp}\|_{{3\alpha-\operatorname{hld}},[a, b]}.$$ 
The controlled RP space $\mathcal{Q}_{\mathbf{X}}^\alpha([a, b], \mathcal{W})$ is a Banach space equipped with the norm $|Y_a|_{\mathcal{W}}+|Y_a^{\dagger}|_{L(\mathcal{V}, \mathcal{W})}+|Y_a^{\dagger\dagger}|_{L(\mathcal{V}, L(\mathcal{V}, \mathcal{W}))}+\|(Y, Y^{\dagger}, Y^{\dagger\dagger})\|_{\mathcal{Q}_{\mathbf{X}}^\alpha,[a, b]}$. 
Due to some computation in \cite[Section 3.1]{2025Inahama}, it deduces that 
\begin{eqnarray}\label{crp2}
\|Y^{\dagger}\|_{\alpha,[a,b]}\le C(1+\|X^{1}\|_{\alpha})(|Y_{a}^{\dagger\dagger}|+\|(Y,Y^{\dagger},Y^{\dagger\dagger})\|_{\mathcal{Q}_{\mathbf{X}}^{\alpha},[a,b]})
\end{eqnarray}
and 
\begin{eqnarray}\label{crp3}
\|Y\|_{\alpha,[a,b]}\le C^{\prime}(1+\|X^{1}\|_{\alpha}^{2}+\|X^{2}\|_{2\alpha})(|Y_{a}^{\dagger}|+|Y_{a}^{\dagger\dagger}|+\|(Y,Y^{\dagger},Y^{\dagger\dagger})\|_{\mathcal{Q}_{X}^{\alpha},[a,b]})
\end{eqnarray}
for some constants $C, C'$ only depend on $b-a, \alpha$.

The distance between two different controlled RPs $(Y, Y^{\dagger},Y^{\dagger\dagger})\in \mathcal{Q}_{\mathbf{X}}^\alpha([a, b], \mathcal{W})$ and $(\tilde{Y}, \tilde{Y}^{\dagger},\tilde{Y}^{\dagger\dagger})\in \mathcal{Q}_{\tilde{\mathbf{X}}}^\alpha([a, b], \mathcal{W})$ is denoted by
$$d_{\mathbf{X}, \tilde{\mathbf{X}}, 3\alpha}\big(Y, Y^{\dagger}, Y^{\dagger\dagger} ; \tilde{Y}, \tilde{Y}^{\dagger}, \tilde{Y}^{\dagger\dagger}\big) \stackrel{\text { def }}{=}\|Y^{\dagger\dagger}-\tilde{Y}^{\dagger\dagger}\|_{{\alpha-\operatorname{hld}},[a, b]}
+\|Y^{\sharp}-\tilde{Y}^{\sharp}\|_{{2\alpha-\operatorname{hld}},[a, b]}
+\|Y^{\sharp\sharp}-\tilde{Y}^{\sharp\sharp}\|_{{3\alpha-\operatorname{hld}},[a, b]}.$$
Then, we state that the function of controlled RP is also a controlled RP for a given RP \cite[Example 3.1]{2025Inahama}.
\begin{rem}\label{function-CRP}
Let $1/4<\alpha<1/3$ and $[a,b]\subset[0,T]$. For a given RP $\mathbf{X}=\left(X^{1}, X^{2}, X^{3}\right)\in \Omega_{\alpha}(\mathcal{V})$ and suppose $(Y, Y^{\dagger},Y^{\dagger\dagger}) \in \mathcal{Q}_{\mathbf{X}}^\alpha([a, b],\mathcal{W})$, for a function $\Phi:\mathcal{W}\mapsto \mathcal{W'}$, we have $(\Phi(Y), \Phi(Y)^{\dagger}, \Phi(Y)^{\dagger\dagger} ) \in \mathcal{Q}_{\mathbf{X}}^\alpha([a, b], \mathcal{W'})$, here 
\begin{eqnarray}\label{crp}
\Phi(Y)_t&:=&\Phi(Y_t)\cr
\Phi(Y)^{\dagger}_t&:=&\nabla \Phi(Y_t)Y^{\dagger}_t\cr
\Phi(Y)^{\dagger\dagger}_t&:=&\nabla \Phi(Y_t)Y^{\dagger\dagger}_t+\nabla^2 \Phi(Y_t)\langle Y^{\dagger}_t\cdot,Y^{\dagger}_t\star\rangle,
\end{eqnarray}
with 
\begin{eqnarray}\label{crp-remain1}
\Phi(Y)^{\sharp}_{s,t}&:=&\nabla g(Y_s)\langle Y_{s,t}^{\sharp}\rangle+\nabla^2 \Phi(Y_s)\langle Y_{s}^{\dagger}X^1_{s,t},Y_{s}^{\dagger\dagger}X^2_{s,t}\rangle+\frac{1}{2}\nabla^2 \Phi(Y_s) \langle Y_{s}^{\dagger\dagger}X^2_{s,t},Y_{s}^{\dagger\dagger}X^2_{s,t}\rangle\cr
&&+\nabla^2 \Phi(Y_s)\langle Y_{s}^{\dagger}X^1_{s,t}+Y_{s}^{\dagger\dagger}X^2_{s,t},Y_{s,t}^{\sharp}\rangle+\frac{1}{2}\nabla^2 \Phi(Y_s)\langle Y_{s,t}^{\sharp},Y_{s,t}^{\sharp}\rangle\cr
&&+\frac{1}{2}\int_{0}^{1}d\theta(1-\theta)^2\nabla^3 \Phi(Y_s+\theta Y_{s,t})\langle Y_{s,t},Y_{s,t},Y_{s,t}\rangle,
\end{eqnarray}
and
\begin{eqnarray}\label{crp-remain2}
\Phi(Y)^{\sharp\sharp}_{s,t}&:=&\nabla \Phi(Y_s)\langle Y_{s,t}^{\sharp\sharp}\rangle+\nabla^2 \Phi(Y_s) \langle Y_{s}^{\dagger\dagger}X^2_{s,t}+Y_{s}^{\sharp},Y_{s}^{\dagger}\rangle\cr
&&+\nabla^2 \Phi(Y_s)\langle Y_{s}^{\dagger}X^1_{s,t}+Y_{s}^{\dagger\dagger}X^2_{s,t}+Y_{s,t}^{\sharp},Y_{s}^{\dagger\dagger}X^1_{s,t}+Y_{s}^{\sharp\sharp}\rangle\cr
&&+\frac{1}{2}\int_{0}^{1}d\theta(1-\theta)^2\nabla^3 \Phi(Y_s+\theta Y_{s,t})\langle Y_{s,t},Y_{s,t},Y_{t}^{\dagger}\rangle.
\end{eqnarray}
(Note that the second term in $\Phi(Y)^{\dagger\dagger}$ is symmetric in $\cdot$ and $\star$.)
\end{rem} 
In the following, we state that the integration of controlled RP against a reference RP is again a controlled RP, whose precise proof refers to \cite[Proposition 3.2]{2025Inahama}.
\begin{rem}\label{integral-CRP}
Let $1/4<\alpha<1/3$ and $[a,b]\subset[0,T]$. For a given RP $\mathbf{X}=\left(X^{1}, X^{2}, X^{3}\right)\in \Omega_{\alpha}(\mathcal{V})$ and controlled RP $(Y, Y^{\dagger},Y^{\dagger\dagger}) \in \mathcal{Q}_{\mathbf{X}}^\alpha([a, b], L(\mathcal{V}, \mathcal{W}))$, we have $\left(\int_a^{\cdot} Y_u d {\mathbf{X}}_u, Y,Y^{\dagger}\right) \in \mathcal{Q}_{\mathbf{X}}^\alpha([a, b], \mathcal{W})$, moreover, 
\begin{eqnarray}\label{integral-crp}
&&\big|\int_s^t Y_u d {\mathbf{X}}_u-(Y_{s}X^1_{s,t}+Y_{s}^{\dagger}X^2_{s,t}+Y_{s}^{\dagger\dagger}X^3_{s,t})\big|\cr
&\le& 2^{4\alpha}\zeta(4\alpha) (t-s)^{4\alpha}\big(\|Y^{\sharp}\|_{{2\alpha-\operatorname{hld}},[a, b]}\|X^{1}\|_{\alpha \mathrm{-hld}} 
+\|Y^{\sharp\sharp}\|_{{3\alpha-\operatorname{hld}},[a, b]} \|X^{2}\|_{2\alpha \mathrm{-hld}}
\big.\cr
&&\big.+\|Y^{\dagger\dagger}\|_{{\alpha-\operatorname{hld}},[a, b]}\|X^{3}\|_{3\alpha \mathrm{-hld}}\big) 
\end{eqnarray}
where $\zeta$ is the Riemann zeta function.
\end{rem}

\subsection{Stability estimate for RDE driven by level 3 RP}
In this section we give stability statements for RDE driven by level 3 RP. Even the RDE and controlled RP theory for the level 2 situation is well-known \cite{2020Friz}, but our RDEs is driven by level 3 RPs where elaborate somewhat tedious computations are required and there are few works stating it \cite{2022Boedihardjo,2025Inahama}. We will therefore provide a detailed explanation below.

We begin by proving the stability of the composition, adapting the argument from \cite[Theorem 7.6]{2020Friz} to our framework.
\begin{prop}\label{function-stability}
Let $1/4<\alpha<1/3$ and give two reference RP $\mathbf{X}=(X^{1}, X^{2}, X^{3}), \tilde{\mathbf{X}}=(\tilde{X}^{1}, \tilde{X}^{2},\tilde{X}^{3})\in \Omega_{\alpha}(\mathcal{V})$ and consider suppose $(Y, Y^{\dagger},Y^{\dagger\dagger}) \in \mathcal{Q}_{\mathbf{X}}^\alpha([0, T], \mathcal{W})$ and $(\tilde{Y}, \tilde{Y}^{\dagger},\tilde{Y}^{\dagger\dagger}) \in \mathcal{Q}_{\tilde{\mathbf{X}}}^\alpha([0, T],\mathcal{W})$ in a bounded set, that is 
in a bounded set, that is
$$
|Y_a|_{\mathcal{W}}+|Y_0^{\dagger}|_{L(\mathcal{V}, \mathcal{W})}+|Y_0^{\dagger\dagger}|_{L(\mathcal{V}, L(\mathcal{V}, \mathcal{W}))}+\|(Y, Y^{\dagger}, Y^{\dagger\dagger})\|_{\mathcal{Q}_{\mathbf{X}}^\alpha,[a, b]}\le M,\quad \vertiii{\mathbf{X}}_{\alpha-\operatorname{hld}}\le M
$$
with an identical bound for $\tilde{\mathbf{X}}$ and $(\tilde{Y}, \tilde{Y}^{\dagger},\tilde{Y}^{\dagger\dagger})$.

Then, for a function $\Phi:\mathcal{W}\mapsto \mathcal{W'}$ of $C^4$, define $(Z, Z^{\dagger},Z^{\dagger\dagger}):=(\Phi(Y),\nabla \Phi(Y)Y^{\dagger}, \nabla \Phi(Y)Y^{\dagger\dagger}+\nabla^2 \Phi(Y)\langle Y^{\dagger}\cdot,Y^{\dagger}\star\rangle)\in \mathcal{Q}_{\mathbf{X}}^\alpha([0, T], \mathcal{W'})$. Similarly, we could define $(\tilde{Z},\tilde{Z}^{\dagger},\tilde{Z}^{\dagger\dagger})\in \mathcal{Q}_{\tilde{\mathbf{X}}}^\alpha([0, T], \mathcal{W'})$. Then, one has the local Lipschitz estimates,
\begin{eqnarray}\label{stability-com}
\|Z, Z^{\dagger}, Z^{\dagger\dagger} ; \tilde{Z}, \tilde{Z}^{\dagger}, \tilde{Z}^{\dagger\dagger}\|_{\mathbf{X}, \tilde{\mathbf{X}}, 3 \alpha} &\le& C_M\big(\rho_\alpha(\mathbf{X}, \tilde{\mathbf{X}})+|Y_0-\tilde{Y}_0|+|Y_0^{\dagger}-\tilde{Y}_0^{\dagger}|+|Y_0^{\dagger\dagger}-\tilde{Y}_0^{\dagger\dagger}|\big. \cr
&&\big.+\|Y, Y^{\dagger}, Y^{\dagger\dagger} ; \tilde{Y}, \tilde{Y}^{\dagger}, \tilde{Y}^{\dagger\dagger}\|_{\mathbf{X}, \tilde{\mathbf{X}}, 3 \alpha}\big),
\end{eqnarray}
moreover, 
\begin{eqnarray}\label{stability-com1}
\|Z-\tilde{Z}\|_{\alpha \mathrm{-hld}} &\le& C_M\big(\rho_\alpha(\mathbf{X}, \tilde{\mathbf{X}})+|Y_0-\tilde{Y}_0|+|Y_0^{\dagger}-\tilde{Y}_0^{\dagger}|+|Y_0^{\dagger\dagger}-\tilde{Y}_0^{\dagger\dagger}|\big. \cr
&&\big.+\|Y, Y^{\dagger}, Y^{\dagger\dagger} ; \tilde{Y}, \tilde{Y}^{\dagger}, \tilde{Y}^{\dagger\dagger}\|_{\mathbf{X}, \tilde{\mathbf{X}}, 3 \alpha}\big),
\end{eqnarray}
for some suitable constant $C_M=C(M,\alpha, \Phi)>0$.
\end{prop}
\para{Proof}. Due to that $\Phi$ is of $C^4$ and some direct computation, we have
\begin{eqnarray*}\label{sc-1}
|Z_0^{\dagger}-\tilde{Z}_0^{\dagger}|&=&|\nabla\Phi(Y_0)Y_0^{\dagger}-\nabla\Phi(\tilde{Y}_0)\tilde{Y}_0^{\dagger}|\leq C_M(|Y_0-\tilde{Y}_0|+|Y_0^{\dagger}-\tilde{Y}_0^{\dagger}|)\cr
|Z_0^{\dagger\dagger}-\tilde{Z}_0^{\dagger\dagger}|&=&|\Phi(Y_0)Y^{\dagger\dagger}_0+\nabla^2 \Phi(Y_0)\langle Y^{\dagger}_0\cdot,Y^{\dagger}_0\star\rangle-\Phi(\tilde{Y}_0)\tilde{Y}^{\dagger\dagger}_0-\nabla^2 \Phi(Y_0)\langle \tilde{Y}^{\dagger}_0\cdot,\tilde{Y}^{\dagger}_0\star\rangle|\cr
&\le& C_M(|Y_0-\tilde{Y}_0|+|Y_0^{\dagger}-\tilde{Y}_0^{\dagger}|+|Y_0^{\dagger\dagger}-\tilde{Y}_0^{\dagger\dagger}|).
\end{eqnarray*}
Then, it suffices to construct the estimate \eqref{stability-com} and the second one \eqref{stability-com1} would follow it. So it needs to estimate 
$$\|Z^{\dagger\dagger}-\tilde Z^{\dagger\dagger}\|_{{\alpha-\operatorname{hld}},[a, b]}
+\|Z^{\sharp}-\tilde Z^{\sharp}\|_{{2\alpha-\operatorname{hld}},[a, b]}
+\|Y^{\sharp\sharp}-\tilde Z^{\sharp\sharp}\|_{{3\alpha-\operatorname{hld}},[a, b]}.$$
For function $\Phi$, we set
$$\Phi(x)-\Phi(y)=g_\Phi(x,y)(x-y), \quad g_\Phi(x,y):=\int_0^1\nabla\Phi(\theta x+(1-\theta)y)d\theta.$$
Firstly, by direct computation, it has
\begin{eqnarray*}\label{sc2}
Z^{\dagger\dagger}_{s,t}&=&\nabla \Phi(Y)_{s,t}Y^{\dagger\dagger}_{t}+\nabla\Phi(Y)_{s}Y^{\dagger\dagger}_{s,t}+g_{\nabla^2\Phi}(Y)_{s,t}\langle Y^{\dagger}_{t}\cdot,Y^{\dagger}_{t}\star\rangle\cr
&&+\nabla^2 \Phi(Y)_{s}\langle Y^{\dagger}_{s,t}\cdot,Y^{\dagger}_{t}\star\rangle+\nabla^2 \Phi(Y)_{s}\langle Y^{\dagger}_{s}\cdot,Y^{\dagger}_{s,t}\star\rangle.
\end{eqnarray*}
Then, 
\begin{eqnarray*}\label{sc2}
Z^{\dagger\dagger}_{s,t}-\tilde{Z}^{\dagger\dagger}_{s,t}
&=& [\nabla \Phi(Y)_{s,t}-\nabla \Phi(\tilde Y)_{s,t}]Y^{\dagger\dagger}_{t}+\nabla \Phi(\tilde Y)_{s,t}[Y^{\dagger\dagger}_{t}-\tilde Y^{\dagger\dagger}_{t}]\cr
&&+ [\nabla \Phi(Y)_{s}-\nabla \Phi(\tilde Y)_{s}]Y^{\dagger\dagger}_{s,t}+\nabla \Phi(\tilde Y)_{s}[Y^{\dagger\dagger}_{s,t}-\tilde Y^{\dagger\dagger}_{s,t}]\cr
&&+[g_{\nabla^2\Phi}(Y)_{s,t}-g_{\nabla^2\Phi}(\tilde Y)_{s,t}]\langle Y^{\dagger}_{t}\cdot,Y^{\dagger}_{t}\star\rangle+g_{\nabla^2\Phi}(\tilde Y)_{s,t}(\langle Y^{\dagger}_{t}\cdot,Y^{\dagger}_{t}\star\rangle-\langle \tilde Y^{\dagger}_{t}\cdot,\tilde Y^{\dagger}_{t}\star\rangle)\cr
&&+ [\nabla^2 \Phi(Y)_{s}-\nabla^2 \Phi(\tilde Y)_{s}]\langle Y^{\dagger}_{s,t}\cdot,Y^{\dagger}_{t}\star\rangle+\nabla^2 \Phi(\tilde Y)_{s}[\langle Y^{\dagger}_{s,t}\cdot,Y^{\dagger}_{t}\star\rangle-\langle \tilde Y^{\dagger}_{s,t}\cdot,\tilde Y^{\dagger}_{t}\star\rangle]\cr
&&+ [\nabla^2 \Phi(Y)_{s}-\nabla^2 \Phi(\tilde Y)_{s}]\langle Y^{\dagger}_{s}\cdot,Y^{\dagger}_{s,t}\star\rangle+\nabla^2 \Phi(\tilde Y)_{s}[\langle Y^{\dagger}_{s}\cdot,Y^{\dagger}_{s,t}\star\rangle-\langle \tilde Y^{\dagger}_{s}\cdot, \tilde Y^{\dagger}_{s,t}\star\rangle]\cr
&=:&\mathcal{I}_1+\mathcal{I}_2+\mathcal{I}_3+\mathcal{I}_4+\mathcal{I}_5.
\end{eqnarray*}
For the first term, with setting $\Delta_t=Y_t-\tilde Y_t$, we obtain that 
\begin{eqnarray}\label{sc3}
&&|\nabla\Phi(Y)_{s,t}-\nabla\Phi({\tilde Y})_{s,t}|\cr
&=&|g_{\nabla \Phi}(Y_t, \tilde Y_t)(\Delta_t-\Delta_s)+(g_{\nabla \Phi}(Y_t, \tilde Y_t)-g_{\nabla \Phi}(Y_s, \tilde Y_s))\Delta_s|\cr
&\le& \|\nabla^2\Phi\|_\infty|Y_{s,t}-\tilde Y_{s,t}|+C\|\nabla^3\Phi\|_\infty(Y_{s,t},\tilde Y_{s,t})|_{\mathcal{W}\times \mathcal{W}}|Y_s-\tilde Y_s|\cr
&\le&\|\nabla^2\Phi\|_\infty|Y_{s,t}-\tilde Y_{s,t}|+C\|\nabla^3\Phi\|_\infty(|Y_{s,t}|+|\tilde Y_{s,t}|)\|Y-\tilde Y\|_{\infty;[0,T]}\cr
&\le& C|t-s|^\alpha(\|\nabla^2\Phi\|_\infty \|Y-\tilde Y\|_{\alpha \mathrm{-hld}}+\|\nabla^3\Phi\|_\infty(\|Y\|_{\alpha \mathrm{-hld}}+\|\tilde Y\|_{\alpha \mathrm{-hld}})\|Y-\tilde Y\|_{\infty;[0,T]}).
\end{eqnarray} 
Meanwhile, $|\nabla \Phi(\tilde Y)_{s,t}[Y^{\dagger\dagger}_{t}-\tilde Y^{\dagger\dagger}_{t}]|\le C|t-s|^\alpha \|\nabla^3\Phi\|_\infty \|Y-\tilde Y\|_{\infty;[0,T]}.$
Combining with the estimate \eqref{crp3}, it is not too difficult to see that 
\begin{eqnarray}\label{sc4}
\|\mathcal{I}_1\|_{\alpha \mathrm{-hld}} &\le& C_M\big(\rho_\alpha(\mathbf{X}, \tilde{\mathbf{X}})+|Y_0-\tilde{Y}_0|+|Y_0^{\dagger}-\tilde{Y}_0^{\dagger}|+|Y_0^{\dagger\dagger}-\tilde{Y}_0^{\dagger\dagger}|\big. \cr
&&\big.+\|Y, Y^{\dagger}, Y^{\dagger\dagger} ; \tilde{Y}, \tilde{Y}^{\dagger}, \tilde{Y}^{\dagger\dagger}\|_{\mathbf{X}, \tilde{\mathbf{X}}, 3 \alpha}\big).
\end{eqnarray} 
The same bound for $\mathcal{I}_2$ is also easy to verify. Then, for the term $g_{\nabla^2\Phi}(Y)_{s,t}-g_{\nabla^2\Phi}(\tilde Y)_{s,t}$, the result similar to \eqref{sc3} also holds but $\|\nabla^3 \Phi\|_\infty$is insteaded by $\|\nabla^4 \Phi\|_\infty$. Therefore, with \eqref{crp2} and \eqref{crp3}, the result \eqref{sc4} also holds for terms $\mathcal{I}_i$ with $i=3,4,5$.

Then, we turns to estimate $\|Z^{\sharp}-\tilde Z^{\sharp}\|_{{2\alpha-\operatorname{hld}},[0, T]}$. Then, with aid of the estimate \eqref{crp-remain1}, we have
\begin{eqnarray}\label{sc5}
\Phi(Y)^{\sharp}_{s,t}-\Phi(\tilde Y)^{\sharp}_{s,t}&=&\nabla \Phi(Y_s)\langle Y_{s,t}^{\sharp}\rangle-\nabla \Phi(\tilde Y_s)\langle \tilde Y_{s,t}^{\sharp}\rangle\cr
&&+\nabla^2 \Phi(Y_s)\langle Y_{s}^{\dagger}X^1_{s,t},Y_{s}^{\dagger\dagger}X^2_{s,t}\rangle-\nabla^2 \Phi(\tilde Y_s)\langle \tilde Y_{s}^{\dagger}\tilde X^1_{s,t},\tilde Y_{s}^{\dagger\dagger}\tilde X^2_{s,t}\rangle\cr
&&+\frac{1}{2}\nabla^2 \Phi(Y_s) \langle Y_{s}^{\dagger\dagger}X^2_{s,t},Y_{s}^{\dagger\dagger}X^2_{s,t}\rangle-\frac{1}{2}\nabla^2 \Phi(\tilde Y_s) \langle \tilde Y_{s}^{\dagger\dagger}\tilde X^2_{s,t},\tilde Y_{s}^{\dagger\dagger}\tilde X^2_{s,t}\rangle\cr
&&+\nabla^2 \Phi(Y_s)\langle Y_{s}^{\dagger}X^1_{s,t}+Y_{s}^{\dagger\dagger}X^2_{s,t},Y_{s,t}^{\sharp}\rangle-\nabla^2 \Phi(\tilde Y_s)\langle \tilde Y_{s}^{\dagger}\tilde X^1_{s,t}+\tilde Y_{s}^{\dagger\dagger}\tilde X^2_{s,t},\tilde Y_{s,t}^{\sharp}\rangle\cr
&&+\frac{1}{2}\nabla^2 \Phi(Y_s)\langle Y_{s,t}^{\sharp},Y_{s,t}^{\sharp}\rangle-\frac{1}{2}\nabla^2 \Phi(\tilde Y_s)\langle \tilde Y_{s,t}^{\sharp},\tilde Y_{s,t}^{\sharp}\rangle\cr
&&+\frac{1}{2}\int_{0}^{1}d\theta(1-\theta)^2\nabla^3 \Phi(Y_s+\theta Y_{s,t})\langle Y_{s,t},Y_{s,t},Y_{s,t}\rangle\cr
&&-\frac{1}{2}\int_{0}^{1}d\theta(1-\theta)^2\nabla^3 \Phi(\tilde Y_s+\theta \tilde Y_{s,t})\langle \tilde Y_{s,t},\tilde Y_{s,t},\tilde Y_{s,t}\rangle.
\end{eqnarray}
By taking similar estimates to \eqref{sc3}--\eqref{sc4}, it is not too difficult to deduce that
\begin{eqnarray}\label{sc6}
\|Z^{\sharp}-\tilde Z^{\sharp}\|_{{2\alpha-\operatorname{hld}},[0, T]} &\le& C_M\big(\rho_\alpha(\mathbf{X}, \tilde{\mathbf{X}})+|Y_0-\tilde{Y}_0|+|Y_0^{\dagger}-\tilde{Y}_0^{\dagger}|+|Y_0^{\dagger\dagger}-\tilde{Y}_0^{\dagger\dagger}|\big. \cr
&&\big.+\|Y, Y^{\dagger}, Y^{\dagger\dagger} ; \tilde{Y}, \tilde{Y}^{\dagger}, \tilde{Y}^{\dagger\dagger}\|_{\mathbf{X}, \tilde{\mathbf{X}}, 3 \alpha}\big).
\end{eqnarray} 
With \eqref{crp-remain2} and same approach in \eqref{sc3}, one has the above bound also holds for $\|Z^{\sharp\sharp}-\tilde Z^{\sharp\sharp}\|_{{3\alpha-\operatorname{hld}},[0, T]}$. Hence, the result \eqref{stability-com} holds. \qed

Then, we give the stability statement of rough integration.
\begin{prop}\label{integral-stability}
Let $1/4<\alpha<1/3$ and give two reference RPs $\mathbf{X}=(X^{1}, X^{2}, X^{3}), \tilde{\mathbf{X}}=(\tilde{X}^{1}, \tilde{X}^{2},\tilde{X}^{3})\in \Omega_{\alpha}(\mathcal{V})$ and consider $(Y, Y^{\dagger},Y^{\dagger\dagger}) \in \mathcal{Q}_{\mathbf{X}}^\alpha([0, T], \mathcal{L}(\mathcal{V},\mathcal{W}))$ and $(\tilde{Y}, \tilde{Y}^{\dagger},\tilde{Y}^{\dagger\dagger}) \in \mathcal{Q}_{\tilde{\mathbf{X}}}^\alpha([0, T],\mathcal{L}(\mathcal{V},\mathcal{W}))$ in a bounded set. Define $$(Z,Z^{\dagger}, Z^{\dagger\dagger}):=\big(\int_0^\cdot Y d\mathbf{X}, Y, Y^{\dagger}\big) \in \mathcal{Q}_{\mathbf{X}}^\alpha([0, T], \mathcal{W})$$
we also define $(\tilde{Z},\tilde{Z}^{\dagger}, \tilde{Z}^{\dagger})$ in a similar way. Then, the local Lipschitz estimate holds as below,
\begin{eqnarray}\label{stability-inte}
\|Z, Z^{\dagger}, Z^{\dagger\dagger} ; \tilde{Z}, \tilde{Z}^{\dagger}, \tilde{Z}^{\dagger\dagger}\|_{\mathbf{X}, \tilde{\mathbf{X}}, 3 \alpha} &\le& C_M\big(\rho_\alpha(\mathbf{X}, \tilde{\mathbf{X}})+|Y_0-\tilde{Y}_0|+|Y_0^{\dagger}-\tilde{Y}_0^{\dagger}|+|Y_0^{\dagger\dagger}-\tilde{Y}_0^{\dagger\dagger}|\big. \cr
&&\big.+T^\alpha\|Y, Y^{\dagger}, Y^{\dagger\dagger} ; \tilde{Y}, \tilde{Y}^{\dagger}, \tilde{Y}^{\dagger\dagger}\|_{\mathbf{X}, \tilde{\mathbf{X}}, 3 \alpha}\big),
\end{eqnarray}
moreover,
\begin{eqnarray}\label{stability-inte1}
\|Z-\tilde{Z}\|_{\alpha \mathrm{-hld}} &\le& C_M\big(\rho_\alpha(\mathbf{X}, \tilde{\mathbf{X}})+|Y_0-\tilde{Y}_0|+|Y_0^{\dagger}-\tilde{Y}_0^{\dagger}|+|Y_0^{\dagger\dagger}-\tilde{Y}_0^{\dagger\dagger}|\big. \cr
&&\big.+T^\alpha\|Y, Y^{\dagger}, Y^{\dagger\dagger} ; \tilde{Y}, \tilde{Y}^{\dagger}, \tilde{Y}^{\dagger\dagger}\|_{\mathbf{X}, \tilde{\mathbf{X}}, 3 \alpha}\big),
\end{eqnarray}
for some suitable constant $C_M=C(M,\alpha)>0$.
\end{prop}
\para{Proof}. Firstly, note that
\begin{eqnarray}\label{si2}
&&|Y^{\dagger}_{s,t}-\tilde{Y}^{\dagger}_{s,t}|\cr
&=&|(Y_{0,s}^{\dagger\dagger}-Y_{0}^{\dagger\dagger})X^1_{s,t}+(\tilde{Y}_{0,s}^{\dagger\dagger}+\tilde{Y}_{0}^{\dagger\dagger})\tilde{X}^1_{s,t}+Y^{\sharp\sharp}_{s, t}-\tilde Y^{\sharp\sharp}_{s, t}|\cr
&\le& C|t-s|^\alpha(|Y^{\dagger\dagger}_0-\tilde{Y}^{\dagger\dagger}_0|+\|\mathbf{X}-\tilde{\mathbf{X}}\|_{\alpha \mathrm{-hld}}+\|Y^{\dagger\dagger}_{0,\cdot}-\tilde{Y}^{\dagger\dagger}_{0,\cdot}\|_\infty+\|Y^{\sharp\sharp}-\tilde Y^{\sharp\sharp}\|_{\alpha \mathrm{-hld}})\cr
&\le& C|t-s|^{\alpha}\big(|Y^{\dagger\dagger}_0-\tilde{Y}^{\dagger\dagger}_0|+\|\mathbf{X}-\tilde{\mathbf{X}}\|_{\alpha \mathrm{-hld}}+T^{\alpha}(\|Y^{\dagger\dagger}-\tilde{Y}^{\dagger\dagger}\|_{\alpha \mathrm{-hld}}+\|Y^{\sharp\sharp}-\tilde Y^{\sharp\sharp}\|_{3\alpha \mathrm{-hld}})\big).
\end{eqnarray}
Furthermore, we could obtain
\begin{eqnarray}\label{si4}
&&|Y_{s,t}-\tilde{Y}_{s,t}|\cr
&=&|(Y_{0,s}^{\dagger}-Y_{0}^{\dagger})X^1_{s,t}+(\tilde{Y}_{0,s}^{\dagger}+\tilde{Y}_{0}^{\dagger})\tilde{X}^1_{s,t}+(Y_{0,s}^{\dagger\dagger}-Y_{0}^{\dagger\dagger})X^2_{s,t}+(\tilde{Y}_{0,s}^{\dagger\dagger}+\tilde{Y}_{0}^{\dagger\dagger})\tilde{X}^2_{s,t}+Y^{\sharp}_{s, t}-\tilde Y^{\sharp}_{s, t}|\cr
&\le& C|t-s|^\alpha(|Y^{\dagger}_0-\tilde{Y}^{\dagger}_0|+|Y^{\dagger\dagger}_0-\tilde{Y}^{\dagger\dagger}_0|+\|\mathbf{X}-\tilde{\mathbf{X}}\|_{\alpha \mathrm{-hld}}+\|Y^{\dagger}_{0,\cdot}-\tilde{Y}^{\dagger}_{0,\cdot}\|_\infty+\|Y^{\dagger\dagger}_{0,\cdot}-\tilde{Y}^{\dagger\dagger}_{0,\cdot}\|_\infty+\|Y^{\sharp}-\tilde Y^{\sharp}\|_{\alpha \mathrm{-hld}})\cr
&\le& C|t-s|^{\alpha}\big(|Y^{\dagger}_{0}-\tilde{Y}^{\dagger}_{0}|+|Y^{\dagger\dagger}_0-\tilde{Y}^{\dagger\dagger}_0|+\|\mathbf{X}-\tilde{\mathbf{X}}\|_{\alpha \mathrm{-hld}}\big.\cr
&&\big.+T^{\alpha}(\|Y^{\dagger}-\tilde{Y}^{\dagger}\|_{\alpha \mathrm{-hld}}+\|Y^{\dagger\dagger}-\tilde{Y}^{\dagger\dagger}\|_{\alpha \mathrm{-hld}}+\|Y^{\sharp}-\tilde Y^{\sharp}\|_{3\alpha \mathrm{-hld}})\big).
\end{eqnarray}
By combining with the estimate \eqref{si2}, it deduces 
\begin{eqnarray*}\label{si5}
\|Y-\tilde{Y}\|_{\alpha \mathrm{-hld}} \le C(\|\mathbf{X}-\tilde{\mathbf{X}}\|_{\alpha \mathrm{-hld}}+|Y_{0}^{\dagger}-\tilde{Y}_{0}^{\dagger}|+|Y_{0}^{\dagger\dagger}-\tilde{Y}_{0}^{\dagger\dagger}|+T^{\alpha}\|Y,Y^{\dagger},Y^{\dagger\dagger};\tilde{Y},\tilde{Y}^{\dagger}, \tilde Y^{\dagger\dagger}\|_{\mathbf{X},\tilde{\mathbf{X}},3\alpha}).
\end{eqnarray*}
Then, it suffices to construct the estimate \eqref{stability-inte} and the second one \eqref{stability-inte1} would follow it. So we need to estimate 
$$\|Z^{\dagger\dagger}-\tilde Z^{\dagger\dagger}\|_{{\alpha-\operatorname{hld}},[0, T]}
+\|Z^{\sharp}-\tilde Z^{\sharp}\|_{{3\alpha-\operatorname{hld}},[0, T]}
+\|Z^{\sharp\sharp}-\tilde Z^{\sharp\sharp}\|_{{2\alpha-\operatorname{hld}},[0, T]}.$$
Firstly, we estimate the second term $\|Z^{\sharp}-\tilde Z^{\sharp}\|_{{3\alpha-\operatorname{hld}},[0, T]}$ in the above. Note that
\begin{eqnarray*}\label{si6}
Z^{\sharp}_{s,t}-\tilde Z^{\sharp}_{s,t} = Z_{s,t}-\tilde Z_{s,t}-(Y_s-\tilde Y_s)X^1_{s,t}-(Y^{\dagger}_s-\tilde Y^{\dagger}_s)X^2_{s,t}.
\end{eqnarray*}
Abbreviate $\mathcal{I} \Sigma:=Z_{s, t}$ and $\Sigma:=Y_s X^1_{s,t}+Y_s^{\dagger}X^2_{s,t}$. Moreover, $\mathcal{I} \tilde\Sigma$ and $\tilde\Sigma$ could be defined in a similar way with respect to $\tilde{Y}$. We set $\mathcal{Q}:=\Sigma-\tilde\Sigma$. After that, we obtain that
\begin{eqnarray}\label{si6}
|Z^{\sharp}_{s,t}-\tilde Z^{\sharp}_{s,t} |&=&\big|(\mathcal{I} \mathcal{Q})_{s, t}-\mathcal{Q}_{s, t}\big|+\big|Y^{\dagger\dagger}_s X^3_{s, t}-\tilde Y^{\dagger\dagger}_s X^3_{s, t}\big|\cr
&\le& C\|\delta \mathcal{Q}\|_{3 \alpha}|t-s|^{3 \alpha}+\big|Y^{\dagger\dagger}_s X^3_{s, t}-\tilde Y^{\dagger\dagger}_s X^3_{s, t}\big|.
\end{eqnarray}
Here, $\delta \mathcal{Q}_{s, u, t}:=\mathcal{I} \mathcal{Q})_{s, t}-\mathcal{Q}_{s, t}$, so we have $\delta \mathcal{Q}_{s, u, t}=\tilde{Y}^{\sharp}_{s,u} \tilde{X}^1_{u, t}-{Y}^{\sharp}_{s,u} {X}^1_{u, t}+\tilde Y^{\sharp\sharp}_{s,u}\tilde X^2_{u, t}-Y^{\sharp\sharp}_{s,u}X^2_{u, t}+ Y^{\dagger\dagger}_{s, u} \tilde{X}^3_{u, t}-Y^{\dagger\dagger}_{s, u} X^3_{u, t}$.
With the help of \remref{integral-CRP}, we obtain that $\|Z^{\sharp}-\tilde Z^{\sharp}\|_{{3\alpha-\operatorname{hld}},[0, T]}$ can be dominated by the right hand side of \eqref{stability-inte}. Next, we have
\begin{eqnarray*}\label{si7}
Z^{\sharp\sharp}_{s,t}=Z^{\dagger}_{s,t}-Z^{\dagger}_{s}X^1_{s,t}=Y_{s,t}-Y^{\dagger}_{s}X^1_{s,t}=Y^{\dagger\dagger}_{s}X^2_{s,t}+Y^{\sharp}_{s,t}.
\end{eqnarray*}
Along the same approach, it is not too difficult to verify that $\|Z^{\sharp\sharp}-\tilde Z^{\sharp\sharp}\|_{{2\alpha-\operatorname{hld}},[0, T]}$ can be dominated by the right hand side of \eqref{stability-inte}.
Finally, with aid of the fact that $Z^{\dagger\dagger}=Y^{\dagger}$, the result \eqref{stability-inte} holds. \qed

We now turn to the stability statement of the solution to the controlled RDE.
\begin{prop} \label{prop2.5}
Let $\xi\in \mathcal{W}$ and $1/4<\alpha<1/3$, give $\mathbf{X}=\left(X^{1}, X^{2}, X^{3}\right)\in \Omega_{\alpha}(\mathcal{V})$. Assume $(Y, \sigma(Y),\sigma(Y)^{\dagger})\in \mathcal{Q}_{\mathbf{X}}^\beta([0, T], \mathcal{W})$ with $1/4<\beta<\alpha$
be the (unique) solution to the following RDE
\begin{eqnarray}\label{2-39}
dY=f(Y_t)dt+\sigma(Y_t) d X_t, \quad Y_0=\xi \in \mathcal{W}.
\end{eqnarray}
Assume that $f:\mathcal{W} \to\mathcal{W}$ is globally bounded and Lipschitz continuous function and $\sigma: \mathcal{W}\to L(\mathcal{V},\mathcal{W})$ are assumed $C_b^4$. Similarly, let $(\tilde{Y}, \sigma(\tilde Y),\sigma(\tilde Y)^{\dagger}) \in \mathcal{Q}_{\tilde{\mathbf{X}}}^\beta([0, T], \mathcal{W})$ with initial value $(\tilde \xi, \sigma(\tilde\xi), \sigma(\tilde\xi)^{\dagger})$. Assume 
$$\vertiii{\mathbf{X}}_{\alpha-\operatorname{hld}},\vertiii{\tilde{\mathbf{X}}}_{\alpha-\operatorname{hld}}\leq M<\infty.$$
Then, we have the (local) Lipschitz estimates as following:
\begin{eqnarray}\label{2-40}
d_{\mathbf{X}, \tilde{\mathbf{X}}, 3\beta}(Y, \sigma(Y),\sigma(Y)^{\dagger} ; \tilde{Y}, \sigma(\tilde{Y}),\sigma(\tilde{Y})^{\dagger}) \leq C_{M,f,\sigma}\big(|\xi-\tilde{\xi}|+\rho_\alpha(\mathbf{X}, \tilde{\mathbf{X}})\big).
\end{eqnarray}
and
\begin{eqnarray}\label{2-41}
\|Y-\tilde{Y}\|_{\beta-\operatorname{hld}} \leq C_{M,f,\sigma}\big(|\xi-\tilde{\xi}|+\rho_\alpha(\mathbf{X}, \tilde{\mathbf{X}})\big).
\end{eqnarray}
Here, $C_{M,f,\sigma}=C(M,\alpha,\beta,L_f,\|\sigma\|_{C_b^4})>0$.
\end{prop}
\para{Proof}. According to the definition of controlled RP, we have
\begin{eqnarray}\label{2-52}
\|Y-\tilde{Y}\|_{\beta-\operatorname{hld}} \leq C(d_{\mathbf{X}, \tilde{\mathbf{X}},3 \beta}\big(Y, Y^{\dagger},\sigma(Y)^{\dagger} ; \tilde{Y}, \tilde{Y}^{\dagger},\sigma(\tilde{Y})^{\dagger}\big)+|\xi-\tilde{\xi}|+\rho_\alpha(\mathbf{X}, \tilde{\mathbf{X}})),
\end{eqnarray}
so it only needs to show that \eqref{2-40} holds, then \eqref{2-41} will follow it.
Let $0<\tau<T$ and we first show that \eqref{2-40} holds in the time interval $[0,\tau]$.
So, we set $\mathcal{M}_{[0, \tau]}^1,\mathcal{M}_{[0, \tau]}^2:\mathcal{Q}_{\mathbf{X}}^\beta([0, \tau], \mathcal{W})\mapsto \mathcal{Q}_{\mathbf{X}}^\beta([0, \tau], \mathcal{W})$ by
\begin{eqnarray}\label{2-42}
\mathcal{M}_{[0, \tau]}^1\left(Y, Y^{\dagger},Y^{\dagger\dagger}\right)&=&\left(\int_0^{\cdot} \sigma\left(Y_s\right) d \mathbf{X}_s, \sigma(Y),\sigma({Y})^{\dagger}\right),\cr
\mathcal{M}_{[0, \tau]}^2\left(Y, Y^{\dagger}, Y^{\dagger\dagger} \right)&=&\left(\int_0^{\cdot} f\left(Y_s\right) ds, 0,0\right)
\end{eqnarray}
and $(Z, Z^{\dagger},Z^{\dagger\dagger}):=\mathcal{M}_{[0, \tau]}^{\xi}:=(\xi, 0)+\mathcal{M}_{[0, \tau]}^1+\mathcal{M}_{[0, \tau]}^2$. Moreover, we stress the fact that the fixed point of $\mathcal{M}_{[0, \tau]}^{\xi}$ is the solution to the \eqref{2-39} on the time interval $[0,\tau]$ for $0<\tau\le T$. Thanks to the fixed point theorem, we arrive at
\begin{eqnarray}\label{2-43}
(Y, \sigma(Y),\sigma^{\dagger}(Y))=\left(Y, Y^{\dagger},Y^{\dagger\dagger}\right)=\left(Z, Z^{\dagger}, Z^{\dagger\dagger}\right)=(Z, \sigma(Y), \sigma^{\dagger}(Y)). 
\end{eqnarray}
By using the estimate in \propref{integral-stability} with minor adjustment by adding drift term, we see that 
\begin{eqnarray}\label{2-44}
&&d_{\mathbf{X}, \tilde{\mathbf{X}}, 3 \beta}\big(Y, Y^{\dagger}, Y^{\dagger\dagger} ; \tilde{Y}, \tilde{Y}^{\dagger}, \tilde{Y}^{\dagger\dagger}\big) \cr
& =&d_{\mathbf{X}, \tilde{\mathbf{X}}, 3 \beta}\big(Z, Z^{\dagger} , {Z}^{\dagger\dagger}; \tilde{Z}, \tilde{Z}^{\dagger}, \tilde{Z}^{\dagger\dagger}\big) \cr
&=&\|Z^{\dagger\dagger}-\tilde Z^{\dagger\dagger}\|_{{\alpha-\operatorname{hld}},[0, \tau]}
+\|Z^{\sharp}-\tilde Z^{\sharp}\|_{{3\alpha-\operatorname{hld}},[0, \tau]}
+\|Z^{\sharp\sharp}-\tilde Z^{\sharp\sharp}\|_{{2\alpha-\operatorname{hld}},[0, \tau]}\cr
& \lesssim& \rho_\alpha(\mathbf{X}, \tilde{\mathbf{X}})+|\xi-\tilde{\xi}|\cr
&&+\tau^\beta d_{\mathbf{X}, \tilde{\mathbf{X}}, 3 \beta}\big(\sigma(Y), {\sigma^\dagger(Y)}
,{\sigma^{\dagger\dagger}( Y)}; \sigma(\tilde Y), {\sigma^\dagger(\tilde Y)}, {\sigma^{\dagger\dagger}(\tilde Y)}\big)+L_f\tau^{\beta}\|Y-\tilde Y\|_{\beta-\operatorname{hld}}.
\end{eqnarray}
Next, with aid of the \propref{function-stability}, we observe that
\begin{eqnarray}\label{2-45}
&&d_{\mathbf{X}, \tilde{\mathbf{X}}, 3 \beta}\big(\sigma(Y), {\sigma^\dagger(Y)},{\sigma^{\dagger\dagger}(Y)}; \sigma(\tilde Y), {\sigma^\dagger(\tilde Y)},{\sigma^{\dagger\dagger}(\tilde Y)},{\sigma^{\dagger\dagger}(\tilde Y)}\big)\cr
&& \lesssim \rho_\alpha(\mathbf{X}, \tilde{\mathbf{X}})+|\xi-\tilde{\xi}|+ d_{\mathbf{X}, \tilde{\mathbf{X}}, 3 \beta}\big(Y, Y^{\dagger} ,Y^{\dagger\dagger}; \tilde{Y}, \tilde{Y}^{\dagger}, \tilde{Y}^{\dagger\dagger}\big).
\end{eqnarray}
Therefore, by using the above estimate \eqref{2-45} and \eqref{2-52}, we show that there exists a positive constant $C_{M,f,\sigma}:=C(M,\alpha,\beta,L_f,\|\sigma\|_{C_b^4})$ such that 
\begin{eqnarray}\label{2-46}
&&d_{\mathbf{X}, \tilde{\mathbf{X}}, 3 \beta}\big(Y, Y^{\dagger}, Y^{\dagger\dagger} ; \tilde{Y}, \tilde{Y}^{\dagger}, \tilde{Y}^{\dagger\dagger}\big) \cr
& \le& C_{M,f,\sigma}\big[\rho_\alpha(\mathbf{X}, \tilde{\mathbf{X}})+|\xi-\tilde{\xi}|+\tau^\beta d_{\mathbf{X}, \tilde{\mathbf{X}}, 3 \beta}\big(Y, Y^{\dagger}, {Y}^{\dagger\dagger} ; \tilde{Y}, \tilde{Y}^{\dagger}, \tilde{Y}^{\dagger\dagger}\big)+\tau^{\beta}\|Y-\tilde Y\|_{\beta-\operatorname{hld}} \big]\cr
& \le& C_{M,f,\sigma}\big[\rho_\alpha(\mathbf{X}, \tilde{\mathbf{X}})+|\xi-\tilde{\xi}|+\tau^\beta d_{\mathbf{X}, \tilde{\mathbf{X}}, 3 \beta}\big(Y, Y^{\dagger} , {Y}^{\dagger\dagger}; \tilde{Y}, \tilde{Y}^{\dagger}, \tilde{Y}^{\dagger\dagger}\big) \big]
\end{eqnarray}
holds.
We could take $\tau>0$ such that $C_{M,f,\sigma}\tau^\beta<1/2$, then
\begin{eqnarray}\label{2-47}
d_{\mathbf{X}, \tilde{\mathbf{X}}, 3 \beta}\big(Y, Y^{\dagger}, Y^{\dagger\dagger} ; \tilde{Y}, \tilde{Y}^{\dagger}, \tilde{Y}^{\dagger\dagger}\big) 
\le C_{M,f,\sigma}(\rho_\alpha(\mathbf{X}, \tilde{\mathbf{X}})+|\xi-\tilde{\xi}|). 
\end{eqnarray}
Then, by using \eqref{2-52}, we show that
\begin{eqnarray}\label{2-51}
\|Y-\tilde{Y}\|_{\beta-\operatorname{hld},[0,\tau]} &\leq&C(d_{\mathbf{X}, \tilde{\mathbf{X}}, 3 \beta}\big(Y, Y^{\dagger}, Y^{\dagger\dagger} ; \tilde{Y}, \tilde{Y}^{\dagger}, \tilde{Y}^{\dagger\dagger}\big) +|\xi-\tilde{\xi}|+\rho_\alpha(\mathbf{X}, \tilde{\mathbf{X}}))\cr
&\leq& C_{M,f,\sigma}\big(|\xi-\tilde{\xi}|+\rho_\alpha(\mathbf{X}, \tilde{\mathbf{X}})\big).
\end{eqnarray}
By taking iteration techniques, it arrives at the result that \eqref{2-40} and \eqref{2-41} hold at the full time interval $[0,T]$.
This proof is completed. \qed

Next, we consider the RDE where the drift term is related to a continuous path. Precisely, for an $\mathcal{V}$-valued continuous path $\psi:[0,T]\to \mathcal{S}$, consider the below RDE driven by RP $\mathbf{X}$ with initial value $\xi:\mathcal{W}$,
\begin{eqnarray}\label{psiequation}
Y_t=\xi+\int_0^t f(Y_s, \psi_s)ds+\int_0^t \sigma(Y_s)d\mathbf{X}_s, \quad Y_t^{\dagger}=\sigma(Y_t), \quad Y^{\dagger\dagger}=\nabla \sigma(Y_t)Y_t^{\dagger},
\end{eqnarray}
for all $t\in[0,T]$. The coefficient $f:\mathcal{W}\times \mathcal{S}\to\mathcal{W}$ is globally bounded and Lipschitz continuous function and $\sigma: \mathcal{W}\to L(\mathcal{V},\mathcal{W})$ is assumed $C_b^4$.

Then, we show properties of the solution map to the above RDE \eqref{psiequation}, whose proof comes from \cite[Proposition 3.3, Proposition 3.5]{2025Inahama} with some small extension, so details of the proof are omitted.
\begin{prop}\label{psi1}
Let $\xi\in \mathcal{W}$ and $1/4<\alpha<1/3$, give $\mathbf{X}=\left(X^{1}, X^{2}, X^{3}\right)\in \Omega_{\alpha}(\mathcal{V})$. There exists a unique global solution $(Y, \sigma(Y),\sigma(Y)^{\dagger})\in \mathcal{Q}_{\mathbf{X}}^\beta([0, T], \mathcal{W})$ with $1/4<\beta<\alpha$ to the RDE \eqref{psiequation}. Assume 
$\vertiii{\mathbf{X}}_{\alpha-\operatorname{hld}}\leq M<\infty$, then, there exists a $\iota>0$ independent of $\xi, \mathbf{X}, \psi, \sigma, f$ such that 
\begin{eqnarray}\label{psisolution}
\|Y\|_{\beta-\operatorname{hld}}\le c_{\alpha,\beta}\{(K+1){(M+1)}\}^{\iota},
\end{eqnarray}
where $c_{\alpha,\beta}>0$ is a constant depending on $\alpha, \beta$ and $K:=\|\sigma\|_{C_b^4} \vee \|f\|_\infty \vee L_f$.
\end{prop}

\begin{prop}\label{psi2}
Let $\xi\in \mathcal{W}$ and $1/4<\alpha<1/3$, give $\mathbf{X}=\left(X^{1}, X^{2}, X^{3}\right)\in \Omega_{\alpha}(\mathcal{V})$. Assume $(Y, \sigma(Y),\sigma(Y)^{\dagger})\in \mathcal{Q}_{\mathbf{X}}^\beta([0, T], \mathcal{W})$ with $1/4<\beta<\alpha$
be the (unique) solution to the above RDE \eqref{psiequation} and 
$\vertiii{\mathbf{X}}_{\alpha-\operatorname{hld}}\leq M<\infty$. Similarly, let $(\tilde{Y}, \sigma(\tilde Y),\sigma(\tilde Y)^{\dagger}) \in \mathcal{Q}_{{\mathbf{X}}}^\beta([0, T], \mathcal{W})$ is the solution map to the below RDE
\begin{eqnarray}\label{6-40}
\tilde Y_t=\xi+\int_0^t \tilde f(\tilde Y_s, \tilde \psi_s)ds+\int_0^t \sigma(\tilde Y_s)d\mathbf{X}_s, \quad \tilde Y_t^{\dagger}=\sigma(\tilde Y_t), \quad \tilde Y^{\dagger\dagger}=\nabla \sigma(\tilde Y_t)\tilde Y_t^{\dagger},
\end{eqnarray}
with initial value $( \xi, \sigma(\xi), \sigma(\xi)^{\dagger})$. $\tilde f:\mathcal{W}\times \mathcal{S}\to\mathcal{W}$ is also assumed globally bounded and Lipschitz continuous.
Then, for any bounded and Lipschitz map $g:\mathcal{W}\to\mathcal{W}$, set
\begin{eqnarray}\label{6-41}
Q_t:=(Y_t-\tilde Y_t)-\int_0^t (g(Y_s)-g(\tilde Y_s))ds-\int_0^t (\sigma(Y_s)-\sigma(\tilde Y_s))d\mathbf{X}_s, \quad t\in[0,T].
\end{eqnarray}
Then, $M\in \mathcal{C}^1(\mathcal{W})$ and the following estimate holds, 
\begin{eqnarray}\label{psiproperty}
\|Y-\tilde Y\|_{\beta-\operatorname{hld}} \le c \exp \big[c\left(K^{\prime}+1\right)^\iota\big(M+1\big)^\iota\big]\|M\|_{3\beta-\operatorname{hld}}
\end{eqnarray}
where $c>0$ is a constant depending on $\alpha, \beta$ and $K^{\prime}:=\|\sigma\|_{C_b^4} \vee \|f\|_\infty \vee L_f\vee \|\tilde f\|_\infty \vee L_{\tilde f}\vee \| g\|_\infty \vee L_{ g}$. And $\iota>0$ independent of $\xi, \mathbf{X}, \psi, \sigma, f, \tilde f, \tilde \psi, g, M$
\end{prop}

\section{Anisotropic geometric RP and GRP lifted by translation of mixed FBM in C-M space}
In this section, we show that the translation of mixed fractional Brownian motion (for $H\in (1/4,1/3)$) in the Cameron--Martin direction admits a lift to GRPs via the anisotropic geometric RP approach. Throughout this section, let $H\in (1/4,1/3)$ and $\alpha\in (0,H)$. We choose suitable $\lfloor 1 / H \rfloor<{p}<\lfloor 1 / H \rfloor+1$ and $(H+1/2)^{-1}<q<2$ such that $1/p+1/q>1$.
\subsection{Mixed fractional Brownian motion and C-M space}\label{sec-2-1}
This subsection features a brief overview of the mixed FBM of Hurst parameter $H$, and only focuses on the case of $H\in(1/4,1/3)$.
The continuous stochastic process $(b^H_t)_{t\in[0,T]}:=(b_t^{H,1},b_t^{H,2},\cdots,b_t^{H,d})_{t\in[0,T]}\in \mathbb{R}^{d}$ starting from 0 is called a FBM if it is a centered Gaussian process, satisfying that
$$\mathbb{E}\big[b_{t}^{H} b_{s}^{H}\big]=\frac{1}{2}\left[t^{2 H}+s^{2 H}-|t-s|^{2 H}\right]\times I_{d}, \quad(0\leq s\leq t \leq T),$$
where $I_{d}$ stands the identity matrix in $\mathbb{R}^{d\times d}$. 
Then, according to some direct computation, it deduces that 
$$\mathbb{E}\big[(b_{t}^{H}-b_{s}^{H})^{2}\big]=|t-s|^{2 H} \times I_{d},\quad (0\leq s\leq t \leq T).$$
From the Kolmogorov's continuity criterion, 
the trajectories of $b^H$ is of $\lfloor 1 / H \rfloor<{p}<\lfloor 1 / H \rfloor+1$-variation and $H'$-H\"older continuous ($H'\in(0,H)$) almost surely. 
The reproducing kernel Hilbert space of $b^H$ is denoted by $\mathcal{H}^{H,d}$.
Thanks to \cite[Proposition 3.4]{2013Inahama}, it admits that each element $g\in \mathcal{H}^{H,d}$ is of finite $(H+1/2)^{-1}<q<2$-variation and $H'$-H\"older continuous.
Next, we consider the $\mathbb{R}^{e}$-valued standard BM $(w_t)_{t\in[0,T]}:=(w_t^1,w_t^2,\cdots,w_t^{e})_{t\in[0,T]}$. We denote the mixed FBM by $(b_t^H, w_t)_{0\le t\le T}\in \mathbb{R}^{d+e}$. It is not too difficult to see that trajectories of mixed FBM $(b^H, w)$ is $\lfloor 1 / H \rfloor<{p}<\lfloor 1 / H \rfloor+1$-variation and $H'$-H\"older continuous ($H'\in(0,H)$) almost surely.

The reproducing kernel Hilbert space for $(w_t)_{t\in[0,T]}$, denoted by $\mathcal{H}^{\frac{1}{2},e}$, which is defined as below, 
$$\mathcal{H}^{\frac{1}{2},e}:=\big\{ k \in \mathcal{C}_0([0,T],\mathbb{R}^{e})\mid k _{t}=\int_{0}^{t} k _{s}^{\prime} d s \text { for } t\in[0,T] \text { with }\| k \|_{\mathcal{H}^{\frac{1}{2},e}}^{2}:=\int_{0}^{T}| k _{t}^{\prime}|_{\mathbb{R}^{e}}^{2} d t<\infty\big\}.$$
Denote $\mathcal{H}:={{\mathcal{H}^{H,d}}\oplus{\mathcal{H}^{\frac{1}{2},e}}}$ the Cameron-Martin subspace related to $(b_t^H, w_t)_{0\le t\le T}$. Then, $(\phi, \psi)\in \mathcal{H}$ is of finite $q$-variation with $(H+1/2)^{-1}<q<2$. 

We first give some necessary notations and variational representation formula for mixed FBM. 
Denote $
S_N=\{(\phi,\psi) \in \mathcal{H}: 
\frac{1}{2}\|(\phi,\psi)\|_{\mathcal{H}}^2 := \frac{1}{2} 
(\|\phi\|_{\mathcal{H}^{H,d}}^{2} +\|\psi \|_{\mathcal{H}^{\frac{1}{2},e}}^{2})
\leq N\}
$ for $N\in\mathbb{N}$, then $S_N$ is a compact Polish space under the weak topology of $\mathcal{H}$.

We denote the set of all 
$\mathbb{R}^{d+e}$-valued processes $(\phi_t,\psi_t)_{t \in [0,T]}$ on 
$(\Omega, {\mathcal{F}}, {\mathbb{P}})$ by ${\mathcal{A}}_b^N$ for $N\in \mathbb{N}$. Let ${\mathcal{A}}_b =\cup_{N\in \mathbb{N}} {\mathcal{A}}_b^N$.
Due to that every $(\phi,\psi)\in {\mathcal{A}}_b^N$ are random variables taking values in the compact ball ${S}_N$, the probability measure family $\{ \mathbb{P}\circ (\phi,\psi)^{-1} : 
(\phi,\psi) \in {\mathcal{A}}_b^N\}$ is automatically tight. 
Thanks to Girsanov's formula, for each $(\phi,\psi) \in {\mathcal{A}}_b$,
the law of $(b^H +\phi,w+\psi)$ is 
mutually absolutely continuous to that of $(b^H,w)$.
Next, we recall the variational representation formula for the mixed FBM, whose proof refers to \cite[Proposition 2.3]{2023Inahama}.
\begin{prop}\label{prop2-1}
Let $\alpha\in (0,H)$. For a bounded Borel measurable function $\Phi:\Omega\to \mathbb{R}$,
\begin{eqnarray}\label{2-1}
-\log \mathbb{E}[\exp (-\Phi(b^H, w))]=\inf _{(\phi, \psi) \in \mathcal{A}_b} \mathbb{E}\big[\Phi(b^H+\phi, w+\psi)+\frac{1}{2}\|(\phi, \psi)\|_{\mathcal{H}}^2\big].
\end{eqnarray}
\end{prop}

\subsection{Anisotropic rough path}\label{sec-3-2}
Firstly, we introduce the anisotropic rough path (RP).
\begin{defn}(Anisotropic rough path)\label{URP}
Set $1/4<\alpha<1/3$ and $2\alpha+\gamma>1$. A continuous map $\left(\Xi^1, \Xi^2, B^3\right): \triangle_T \rightarrow \mathcal{V} \oplus \mathcal{V}^{\otimes 2} \oplus \mathcal{V}_1^{\otimes 3}$ is called anisotropic RP of roughness $(\alpha,\gamma)$ if it satisfies the following conditions:
\begin{itemize}
\item (H\"older condition) $\max \{\|B^1\|_\alpha,\|B^2\|_{2 \alpha},\|B^3\|_{3 \alpha},\|W^1\|_\gamma,\|W^2\|_{2 \gamma},\|I[B, W]\|_{\alpha+\gamma},\|I[W, B]\|_{\alpha+\gamma}\}<\infty$;
\item (Chen's relation) 
$ \Xi_{s, t}^1=\Xi_{s, u}^1+\Xi_{u, t}^1, \quad \Xi_{s, t}^2=\Xi_{s, u}^2+\Xi_{u, t}^2+\Xi_{s, u}^1 \otimes \Xi_{u, t}^1$, $ B_{s, t}^3=B_{s, u}^3+B_{u, t}^3+B_{s, u}^1 \otimes B_{u, t}^2+B_{s, u}^2 \otimes B_{u, t}^1 .$
\end{itemize}
\end{defn}
Here, $\mathcal{V}^{\otimes2}=\oplus_{i,j=1,2}\mathcal{V}_i\otimes\mathcal{V}_j$ and $\mathcal{V}^{\otimes3}=\oplus_{i,j,k=1,2}\mathcal{V}_{i}\otimes\mathcal{V}_{j}\otimes\mathcal{V}_{k}$. Set $\pi_{i}: \mathcal{V}\to\mathcal{V}_i$ the projection onto the $i$ component, similarly, $\pi_{ij}: {\mathcal{V}^{\otimes2}}\to\mathcal{V}_{i}\otimes\mathcal{V}_{j}(i,j=1,2)$ and $\pi_{ijk}:\mathcal{V}^{\otimes3}\to\mathcal{V}_{i}\otimes\mathcal{V}_{j}\otimes\mathcal{V}_{k}(i,j,k=1,2)$.
We denote $ \hat{\Omega}_{\alpha,\gamma}(\mathcal{V})$ the set of all anisotroic RPs of roughness $(\alpha, \gamma)$. 

For a path $z\in \mathcal{C}^1_0(\mathcal{V})$, denote its natural lift by $\hat{S}(z)\in \hat{\Omega}_{\alpha,\gamma}(\mathcal{V})$, which could be defined in Riemann-Stieltjes integral. Set $G\hat{\Omega}_{\alpha,\gamma}(\mathcal{V})$ to be the closure of the $\hat{S}(z)$, which is the geometric anisotropic RP space. Note that $G\hat{\Omega}_{\alpha,\gamma}(\mathcal{V})$ is a complete separable metric space.
There is a canonical continuous injection from $G\hat{\Omega}_{\alpha,\gamma}(\mathcal{V})$ to $G\Omega_\alpha(\mathcal{V})$. Next, the GRP could be constructed according to the anisotropic RP. Write $\Xi=\mathbf{Ext}({\Xi^{1}},\Xi^{2},B^{3})$, so $\Xi=(\Xi^1,\Xi^2,\Xi^3)\in G\Omega_\alpha(\mathcal{V})$. Denote $\pi_{ijk}{:}\mathcal{V}^{\otimes3}\to\mathcal{V}_{i}\otimes\mathcal{V}_{j}\otimes\mathcal{V}_{k}$. Here, $\Xi^3$ is constructed as below: 
\begin{eqnarray}\label{GRP3}
\Xi_{s,t}^{3,[ijk]}:=\lim_{|\mathcal{P}|\searrow0}\sum_{l=1}^{N}\pi_{ijk}\langle\Xi_{s,t_{l-1}}^{1}\otimes\Xi_{t_{l-1},t_{l}}^{2}+\Xi_{s,t_{l-1}}^{2}\otimes\Xi_{t_{l-1},t_{l}}^{1}\rangle
\end{eqnarray}
for $i,j,k\neq(1,1,1),(2,2,2)$ and $\mathcal{P}:=\{s=t_0<t_1<\cdots<t_N=t\}$ is the partition of $[s,t]$. Then, $\Sigma^{3,ijk}:=B^{H,3}$ for $i,j,k=(1,1,1)$ and $\Sigma^{3,ijk}:=W^{3}$ for $i,j,k=(2,2,2)$. It is easy to verify that if $(\Xi^1,\Xi^2,B^3)=\hat{S}(z)$ for $z\in\mathcal{C}^1_0(\mathcal{V})$, then $(\Xi^1,\Xi^2,\Xi^3)={S}_3(z)\in G\Omega_\alpha(\mathcal{V})$, whose proofs refer to \cite[Propsition 4.1, Proposition 4.2]{2025Inahama}.

\subsection{The translation of GRP lifted by Mixed FBM}\label{sec-3-4}
In this subsection, we will show that the mixed FBM could be lifted to GRPs and anisotropic GRPs. Firstly, some prior estimates are given before. 

Then, according to \cite[Theorem 15.45, Proposition 15.5]{2010Friz} or \cite[Lemma 4.10]{2025Inahama}, we have the following result.
\begin{prop}\label{MFBM-RP}
Let $H\in (\frac{1}{4},\frac{1}{3})$ and $\alpha\in (0,H)$. The mixed FBM $(b^H,w)\in \mathbb{R}^{d+e}$ could be lifted to anisotropic RP $(\Sigma^1,\Sigma^2,B^{H,3})$, where the first level path and second level path are defined as below:
\begin{eqnarray}\label{MFBM-RP1}
\Sigma^{1}=(b^H_{s,t}, w_{s,t}), \quad \Sigma^{2}=\left(B^{H,2}+I[b^H,w]+I[w,b^H]+W^2\right)_{s, t},
\end{eqnarray}
with
\begin{eqnarray}\label{MFBM-RP2}
I[b^H,w]_{s,t}=\int_s^t b^H_{s,r}, dw^{I}_{r}, \quad I[w,b^H]_{s,t}=b^H_{s,t}\otimes w_{s,t}-\int_s^t(d^\mathrm{I}w_u)\otimes B_{s,u}^{H,1}.
\end{eqnarray}
$B^{H,3}$ is the third level path of the GRP lifted by the FBM $b^H\in \mathbb{R}^{d}$. Then, $(\Xi^1,\Xi^2,B^3)\in G\hat{\Omega}_{\alpha,\gamma}(\mathbb{R}^{d+e})$ a.s. Moreover, $(\Xi^1,\Xi^2,\Xi^3)\in G{\Omega}_{\alpha}(\mathbb{R}^{d+e})$ a.s.
\end{prop}

Next, we will state that C-M elements of FBM can be lifted to GRP.

\begin{prop}\label{FBMCM-RP}
Let $H\in (1/4,1/3)$ and $\alpha\in (0,H)$. The elements $u\in \mathcal{H}^{H,d}$ can be lifted to RP $U=(U^1,U^2, U^3)\in G\Omega_{\alpha}(\mathcal{H}^{H,d})$ with $U=S_3(u)$ which is well-defined in the variation setting. Moreover, $U=(U^1,U^2, U^3)$ is a locally Lipschitz continuous mapping from $\mathcal{H}^{H,d}$ to $G\Omega_\alpha(\mathbb{R}^d)$.
\end{prop}
\para{Proof}. It is known that each element $g\in \mathcal{H}^{H,d}$ is finite $(H+1/2)^{-1}<q<2$-variation, then $S_3(u)$ is well-defined in the variation setting, that is for $k=1,2,3$
\begin{eqnarray}\label{FBMCM-RP1}
|S^i(u)_{s,t}|\le C_q\|u\|_{p \text {-var }} \le C'_q\|u\|^k_{\mathcal{H}^{H,d}} |t-s|^{i(\frac{1}{q}-\frac{1}{2})}
\end{eqnarray}
for all $(s,t)\in \Delta_T$. Here, $C_q>0$ is a constant depending on $q$.

Besides, it is direct to verify that $S^k(u)$ is an $\alpha$-H\"older continuous weakly GRPs for every $1/4<\alpha<H$ by \eqref{FBMCM-RP1} . With aid of the result \cite[Theorem 8.22]{2010Friz}, $S^k(u)\in G\Omega_\alpha(\mathbb{R}^d)$ for every $1/4<\alpha<H$. Furthermore, it is natural to see that $U=(U^1,U^2, U^3)$ is a locally Lipschitz continuous mapping from $\mathcal{H}^{H,d}$ to $G\Omega_\alpha(\mathbb{R}^d)$. The proof is completed. \qed

Similarly, it deduces that C-M elements of BM could be lifted to GRPs.
\begin{rem}\label{BMCM-RP}
The elements $v\in \mathcal{H}^{\frac{1}{2},e}$ can be lifted to RP $V=(V^1,V^2)\in G\Omega_{\alpha}(\mathbb{R}^e)$ with
\begin{eqnarray*}\label{3-10}
V^1_{s,t}=v_{s,t},\quad V^2_{s,t}=\int_{s}^{t}v_{s,r}dv_r
\end{eqnarray*}
where $V^2$ is well-defined in the Young sense since $v$ is differentiable.
\end{rem}

\begin{prop}\label{FBM+CM-RP}
Let $H\in (1/4,1/3)$ and $\alpha\in (0,H)$. Let $b^H+u$ be the translation of FBM $B^H\in \mathbb{R}^{d}$ in the direction $u \in \mathcal{H}^{H,d}$. Then, $b^H+u$ can be lifted to GRP $\mathcal{T}^U(B^H)=(\mathcal{T}^{U,1}(B^H),\mathcal{T}^{U,2}(B^H), \mathcal{T}^{U,3}(B^H))\in G\Omega_{\alpha}(\mathbb{R}^{d})$.
\end{prop}
\para{Proof}. It is easy to verify that in our situation, the complementary Young regularity condition in \cite[Condition 15.56]{2010Friz} is satisfied. Then with the help of \cite[Lemma 15.58]{2010Friz}, for $\mathbb{P}$-a.e. $\omega$, we have for $u\in \mathcal{H}^{H,d}$, $B^{H}(\omega+u)=\mathcal{T}^{u}(B^H)$ where $B^H=\lim_{m\to\infty}S_3(b^H(m))$ is GRPs and $\mathcal{T}$ is the translation operator for GRPs. So it remains to verify that $\mathcal{T}^{U,i}(B^H)$ is well-defined and of $i\alpha$-H\"older continuous for $i=1,2,3$, which is a complement to the result \cite[Lemma 15.58]{2010Friz}.

The first and second level paths are defined as below:
\begin{eqnarray}\label{FBM-CMRP1}
\mathcal{T}_{s,t}^{U,1}(B^H)=(b^H+u)_{s,t}, \quad \mathcal{T}_{s, t}^{U,2}(B^H)=\left(B^{H,2}+I[b^H,u]+I[u,b^H]+U^2\right)_{s, t}.
\end{eqnarray}
Since the trajectories of $b^H$ are of ${p}$-variation almost surely for $\lfloor 1 / H \rfloor<{p}<\lfloor 1 / H \rfloor+1$ and $u\in \mathcal{H}^{H,d}$ is
of finite $(H+1/2)^{-1}<q<2$-variation. Since $\frac{1}{p}+\frac{1}{q}>1$, the integral $\int_{s}^{t}b^H_{s,r}du_r$ is well-defined in the Young sense. Similar to \cite[Remark 2.6]{2025Yang}, it could be verified that $I[b^H,u]$ is $2\alpha$-H\"older continuous. By taking similar estimates, the other remaining terms in the second level path are also well-defined in the Young sense and of $2\alpha$-H\"older continuous.

The analysis for the third level path will be a little complex. It is defined by the following sense,
\begin{eqnarray}\label{FBM-CMRP2}
&&\mathcal{T}_{s,t}^{U,3}(B^H)=B_{s,t}^{H,3}+U^3_{s,t}+\sum_{i=1}^{4}D^i_{s,t}, 
\end{eqnarray}
where
\begin{eqnarray}\label{FBM-CMRP3}
D^1_{s,t}&:=&\int_{s\le t_1\le t_2\le t_3\le t} du_{t_1}\otimes
du_{t_2}\otimes db^H_{t_3}+\int_{s\le t_1\le t_2\le t_3\le t} du_{t_1}\otimes db^H_{t_2} \otimes du_{t_3}\cr
&&+\int_{s\le t_1\le t_2\le t_3\le t} db^H_{t_1}\otimes du_{t_2}\otimes du_{t_3}, \cr
D^2_{s,t}&:=&\int_{s\le t_1\le t_2\le t_3\le t} db^H_{t_1}\otimes du_{t_2}\otimes db^H_{t_3}, \cr
D^3_{s,t}&:=& \int_{s\le t_1\le t_2\le t_3\le t} db^H_{t_1}\otimes db^H_{t_2}\otimes du_{t_3}, \quad D^4_{s,t}:=\int_{s\le t_1\le t_2\le t_3\le t} du_{t_1}\otimes db^H_{t_2}\otimes db^H_{t_3}, 
\end{eqnarray}
Then, $D^1$ is well-defined in the Young integral sense, and for every $b^H, u$ and $(s,t)\in \Delta_T$,
\begin{eqnarray}\label{FBM-CMRP4}
|D^1_{s,t}|\le C_4\|u\|^2_{q \text {-var }}\|b^H\|_{p \text {-var }} \le C_4\|u\|_{\mathcal{H}^{H,d}}^2\|b^H\|_{\alpha-\operatorname{hld}} (t-s)^{\alpha+2(\frac{1}{q}-\frac{1}{2})}.
\end{eqnarray}
Since that we could choose suitable $q$ such that $\frac{1}{q}-\frac{1}{2}$ is close to $H$, then $D^1$ is of $3\alpha$-H\"older continuous. The above estimation also works for the term $D^2$. For the term $D^3$, we have
$$D^3_{s,t}=\int_{s\le t_1\le t_2\le t_3\le t} db^H_{t_1}\otimes db^H_{t_2}\otimes du_{t_3}=B^{H,2}_{s,r}\otimes du_{r},$$
Along the similar approach as in \eqref{FBM-CMRP4}, we could see that $D^3$ are also well-defined in the Young sense and of $3\alpha$-H\"older continuous.

Then, we turn to estimate $D^4$. Firstly we will estimate $D^4$ for differentiable sample paths. Precisely, for $b^H(m),u(m)\in \mathcal{C}^1$, by leveraging the Fubini theorem, it deduces,
\begin{eqnarray}\label{FBM-CMRP5}
D^4_{s,t}=\int_{s\le t_1\le t_2\le t_3\le t} (u(m))'_{t_1}d_{t_1}\otimes (b^H(m))'_{t_2}d{t_2}\otimes (b^H(m))'_{t_3}d{t_3}=\int_{s}^t du_{t_1}\otimes B^{H,2}_{t_1,t}.
\end{eqnarray}
Then, by the limiting argument and the above representation for GRPs, it deduces that
\begin{eqnarray}\label{FBM-CMRP6}
|D^4_{s,t}|\le C_5(\|B^{H,1}\|_{\alpha-\mathrm{hld}}^{2}+\|B^{H,2}\|_{2\alpha-\mathrm{hld}})\|u\|_{\mathcal{H}^{H,d}}(t-s)^{2\alpha+(\frac{1}{q}-\frac{1}{2})}.
\end{eqnarray}
By choosing suitable $q$ such that $\frac{1}{q}-\frac{1}{2}$ is close to $H$, then $D^4$ is of $3\alpha$-H\"older continuous. 
The proof is completed. \qed

Thanks to the \cite[Proposition 13.37]{2010Friz}, the translation of standard BM in the direction of C-M space still could be lifted to GRPs.
\begin{rem}\label{BM+CM-RP}
Let $w+v$ be the translation of BM $w\in \mathbb{R}^{e}$ in the direction $v \in \mathcal{H}^{\frac{1}{2},e}$ and $\alpha\in (\frac{1}{3}, \frac{1}{2})$. Then, $w+v$ can be lifted to GRP $\mathcal{T}^V(W)=(\mathcal{T}^{V,1}(W),\mathcal{T}^{V,2}(W))\in G\Omega_{\alpha}(\mathbb{R}^{e})$, which is defined as following:
\begin{eqnarray}\label{BM-CMRP1}
\mathcal{T}_{s,t}^{V,1}(W)=(w+v)_{s,t}, \quad \mathcal{T}_{s, t}^{V,2}(W)=\left(W^2+I[w,v]+I[v,w]+V^2\right)_{s, t}.
\end{eqnarray}
\end{rem}

Finally, we intend to show that the translation of mixed FBM in the direction of C-M space can still be lifted to a anisotropic GRPs.
\begin{prop}\label{MFBM+CM-RP}
Let $H\in (1/4,1/3)$ and $\alpha\in (0,H)$. Let $(b^H+u,w+v)$ be the translation of $(b^H,w)^\mathrm{T}\in \mathbb{R}^{d+e}$ in the direction $h:=(u,v)\in \mathcal{H}$. The translation $(b^H+u,w+v)$ could be lifted to anisotropic GRP $(\mathcal{T} \Sigma^1,\mathcal{T} \Sigma^2,(B^H+U)^{3})\in G\hat{\Omega}_{\alpha,\gamma}( \mathbb{R}^{d+e})$. Directly, it also could be lifted to GRP $\mathcal{T}\Sigma:=(\mathcal{T} \Sigma^1,\mathcal{T} \Sigma^2,\mathcal{T} \Sigma^3):=\mathbf{Ext} (\mathcal{T} \Sigma^1,\mathcal{T} \Sigma^2,(B^H+U)^{3})\in G\Omega_{\alpha}(\mathbb{R}^{d+e})$.
\end{prop}
\para{Proof}. 
Similar to \propref{FBM+CM-RP}, if $(b^H+u,w+v)\in \mathbb{R}^{d+e}$ could be lifted to $(\mathcal{T} \Sigma^1,\mathcal{T} \Sigma^2,\mathcal{T}_{s,t}^{U,3}(B^H))\in G\hat{\Omega}_{\alpha,\gamma}( \mathbb{R}^{d+e})$ with $1/4<\alpha<1/3$ and $2\alpha+\gamma>1$. Moreover, due to the canonical continuous injection from $G\hat{\Omega}_{\alpha,\gamma}(\mathbb{R}^{d+e})$ to $G\Omega_\alpha(\mathbb{R}^{d+e})$, it suffices to show that $\mathcal{T} \Sigma:=(\mathcal{T} \Sigma^1,\mathcal{T} \Sigma^2,\mathcal{T} \Sigma^3)\in G{\Omega}_{\alpha}(\mathbb{R}^{d+e})$. 

The former two level path are defined as below:
\begin{eqnarray}\label{MFBM+CM-RP7}
&&\mathcal{T} \Sigma_{s,t}^{1}=(b^H+u,w+v)_{s,t}, \cr
&&\mathcal{T} \Sigma_{s,t}^{2}\cr
&=&\left(\begin{array}{ll}
B^{H,2}+I[b^H,u]+I[u,b^H]+ U^2& I[b_H, w]+I[b_H, v]+I[u, w]+I[u, v] \\
I[w, b_H]+I[w, u]+ I[v, b_H]+ I[v, u] & W^2+I[w,v]+I[v,w]+V^2
\end{array}\right)_{s, t}\cr
&=&(B^H,W^2)_{s t}+
\left(\begin{array}{ll}
I[b^H,u]+I[u,b^H]+ U^2& I[b_H, v]+I[u, w]+I[u, v] \\
I[w, u]+ I[v, b_H]+ I[v, u] & I[w,v]+I[v,w]+V^2
\end{array}\right)_{s, t}.
\end{eqnarray}
Thanks to \propref{FBMCM-RP}, and \remref{BMCM-RP} we have already shown that $B^H, W$, $U^2$ and $V^2$ are well-defined. 
By \propref{FBM+CM-RP}, we have $\mathcal{T}^U(B^H)=(\mathcal{T}^{U,1}(B^H),\mathcal{T}^{U,2}(B^H), \mathcal{T}^{U,3}(B^H))\in \Omega_{\alpha}(\mathbb{R}^{d})$. Define $\mathcal{T} \Sigma^{2}(m)$ as follow, 
\begin{eqnarray}\label{MFBM+CM-RP8}
\mathcal{T} \Sigma^{2}(m)&=&I[b^{H}(m)+u(m),b^{H}(m)+u(m)]+I[w(m)+v(m),w(m)+v(m)]\cr
&&+I[b^{H}(m)+u(m),w(m)+v(m)]+I[w(m)+v(m),b^{H}(m)+u(m)]\cr
&=:&I[b^{H}(m)+u(m),b^{H}(m)+u(m)]+I[w(m)+v(m),w(m)+v(m)]\cr
&&+E_{s,t}^1(m)+E_{s,t}^2(m).
\end{eqnarray}
The next step is to show that the following convergence holds,
\begin{eqnarray}\label{MFBM+CM-RP9}
\lim_{m\to\infty}\|E_{s,t}^1(m)-I[b^{H}+u,b^{H}+u]\|_{\alpha+\gamma-2\kappa-\operatorname{hld}}=0.\cr
\lim_{m\to\infty}\|E_{s,t}^2(m)-I[w+v,w+v]\|_{\alpha+\gamma-2\kappa-\operatorname{hld}}=0.
\end{eqnarray} 
Firstly, we will estimate $E_{s,t}^1(m)$, which could be rewritten as below,
\begin{eqnarray}\label{MFBM+CM-RP10}
E_{s,t}^1(m)&=&I[b^H(m),w(m)]+I[u(m),v(m)]+I[u(m),w(m)]+I[b^H(m),v(m)]\cr
&=:&\sum_{i=1}^4 E_{s,t}^{1,i}(m).
\end{eqnarray} 
For terms $E_{s,t}^{1,1}(m)$ and $E_{s,t}^{1,3}(m)$, the convergence could be obtained by the similar estimate in \cite[Section 4]{2015Inahama}. 

Then it remain to verify the convergence of $E_{s,t}^{1,2}(m)$ and $E_{s,t}^{1,4}(m)$ as $m\to \infty$. Since details of the proof is similar, so we show the proof for term $E_{s,t}^{1,2}(m)$ as an example.
Recall that $I[u,v]$ is well-defined in the variation setting since $2/q>1$ and $\alpha+\gamma$-H\"older continuous. For $m\in \mathbb{N}$, we denote $(u(m),v(m))$ be the $m$-th dyadic piecewise linear approximation of $(u,v)$ associated with the partition $\{kT2^{-m}|0\le k\le 2^m\}$ of time interval $[0,T]$. Hence, it suffices to show
\begin{eqnarray}\label{MFBMCM-RP2}
\lim_{m\to\infty}\|I[u(m),v(m)]-I[u,v]\|_{\alpha+\gamma-2\kappa-\operatorname{hld}}=0
\end{eqnarray}
for any small parameter $\kappa\in(0,\frac{\alpha+\gamma}{2})$.

To show the convergence \eqref{MFBMCM-RP2}, it requires to show that 
\begin{eqnarray}\label{MFBMCM-RP3}
\sup_{m\in\mathbb{N}}(\|I[u(m),v(m)]\|_{\alpha+\gamma-\operatorname{hld}}+\|I[u,v(m)]\|_{\alpha+\gamma-\operatorname{hld}})\le C,
\end{eqnarray}
and 
\begin{eqnarray}\label{MFBMCM-RP4}
\lim_{m\to\infty}\|I[u,v(m)]-I[u,v]\|_{\alpha+\gamma-2\kappa-\operatorname{hld}}=0
\end{eqnarray}
for any small parameter $\kappa\in(0,\frac{\alpha+\gamma}{2})$.

Then, it deduces
\begin{eqnarray*}\label{MFBMCM-RP22}
\|I[u(m),v(m)]-I[u,v]\|_{\alpha+\gamma-2\kappa-\operatorname{hld}}&\le& \sup_{m\in\mathbb{N}}(\|I[u(m),v(m)]\|_{\alpha+\gamma-\operatorname{hld}}+\|I[u,v(m)]\|_{\alpha+\gamma-\operatorname{hld}})\times \frac{1}{2^{\kappa m}}\cr
&&+\|I[u,v(m)]-I[u,v]\|_{\alpha+\gamma-2\kappa-\operatorname{hld}}.
\end{eqnarray*}
Our tasks are to show \eqref{MFBMCM-RP3} and \eqref{MFBMCM-RP4}. Firstly, we will show that \eqref{MFBMCM-RP3} holds. We set $[u]_k^m:=2^m\int_{t_{k-1}^m}^{t_k^m}u_rdr$ with $m\in\mathbb{N},1\leq k\leq2^{m}$. Then, when $s,t\in[t_{k-1}^m,t_k^m]$,
\begin{eqnarray}\label{MFBMCM-RP11}
|I[u,v(m)]_{s,t}|=\big|\int_{s}^{t}u_{s,r}\frac{\Delta_{k}^{m}v}{2^{-m}}dr\big|\le\|u\|_{\alpha-\operatorname{hld}}\frac{(t-s)^{\alpha+\frac{1}{2}}}{\alpha+1}\big|\frac{\Delta_{k}^{m}v}{2^{-m/2}}\big|\le C\|u\|_{\alpha-\operatorname{hld}}(t-s)^{\alpha+\frac{1}{2}}
\end{eqnarray}
for some constant $C>0$.

When $s\in[t_{k-1}^m,t_k^m]$ and $t\in[t_l^m,t_{l+1}^m]$ with $k\leq l$, by using Chen's relation, we have
\begin{eqnarray}\label{MFBMCM-RP12}
I[u,v(m)]_{s,t}&=&I[u,v(m)]_{s,t_{k}^{m}}+I[u,v(m)]_{t_{k}^{m},t_{l}^{m}}+I[u,v(m)]_{t_{l}^{m},t}\cr
&&+u_{s,t_{k}^{m}}^{1}v(m)_{t_{k}^{m},t_{l}^{m}}^{1}+u_{s,t_{k}^{m}}^{1}v(m)_{t_{l}^{m},t}^{1}+u_{t_{k}^{m},t_{l}^{m}}^{1}v(m)_{t_{l}^{m},t}^{1}\cr
&=:& \sum_{i=1}^6 A_{i}.
\end{eqnarray}
For the terms $A_{1}$ and $A_{3}$, the estimate is similar to that in \eqref{MFBMCM-RP11}. Since $v$ is differentiable, it is easy to see that $\sum_{i=4}^6 |A_4|<C\|u\|_{\alpha-\operatorname{hld}}(t-s)^{\alpha+\frac{1}{2}}$. Then, we will estimate $A_2$ as below, 
\begin{eqnarray}\label{MFBMCM-RP13}
A_{2}=\sum_{j=k+1}^{l}[u_{\cdot}-u_{t_{k}^{m}}]_{j}^{m}(\Delta_{j}^{m}v)=\int_{t_{k}^{m}}^{t_{l}^{m}}\sum_{j=k+1}^{l}[u_{\cdot}-u_{t_{k}^{m}}]_{j}^{m}\mathbf{1}_{[t_{j-1}^{m},t_{j}^{m}]}(u)dv_{r}.
\end{eqnarray}
By the H\"older inequality and the property that $v$ is differentiable, it deduces that 
\begin{eqnarray}\label{MFBMCM-RP14}
|A_2|&\le& C\bigg(\int_{t_k^m}^{t_l^m}\big|\sum_{j=k+1}^l[u_\cdot-u_{t_k^m}]_j^m\mathbf{1}_{[t_{j-1}^m,t_j^m]}(r)\big|^2dr\bigg)^{\frac{1}{2}}\big(\int_{t_k^m}^{t_l^m}(v'_r)^2dr\big)^{\frac{1}{2}}\cr
&\le& C\|u\|_{\alpha-\operatorname{hld}}\|v\|_{\mathcal{H}^{\frac{1}{2},e}}(t_l^m-t_k^m)^{\alpha+1/2},
\end{eqnarray}
Then, similar estimate could be obtained for the $\sup_{m\in\mathbb{N}}\|I[u(m),v(m)]\|_{\alpha+\gamma-\operatorname{hld}}$. Then, the estimate \eqref{MFBMCM-RP3} has already been verified.

Next, we will show that \eqref{MFBMCM-RP4} holds. Then, when $s,t\in[t_{k-1}^m,t_k^m]$, 
\begin{eqnarray}\label{MFBMCM-RP5}
|I[u,v(m)]_{s,t}-I[u,v]_{s,t}|&\le&|I[u,v(m)]_{s,t}|+|I[u,v]_{s,t}|\cr
&\le& C_{1}\|u\|_{\alpha-\operatorname{hld}}(t-s)^{\alpha+\frac{1}{2}}\cr
&\le& C_{1}(2^{-m})^{\kappa}\|u\|_{\alpha-\operatorname{hld}}(t-s)^{\alpha+\frac{1}{2}-\kappa}.
\end{eqnarray}
When $s\in[t_{k-1}^m,t_k^m]$ and $t\in[t_l^m,t_{l+1}^m]$ with $k\leq l$, then
\begin{eqnarray}\label{MFBMCM-RP6}
I[u,v(m)]_{s,t}-I[u,v]_{s,t}&=&\{I[u,v(m)]_{s,t_{k}^{m}}-I[u,v]_{s,t_{k}^{m}}\}+\{I[u,v(m)]_{t_k^m,t_l^m}-I[u,v]_{t_k^m,t_l^m}\}\cr
&&+\{I[u,v(m)]_{t_{l}^{m},t}-I[u,v]_{t_{l}^{m},t}\}+g_{s,t_k^m}^1\{v(m)_{t_k^m,t_l^m}^1-v_{t_k^m,t_l^m}^1\}\cr
&&+g_{s,t_{k}^{m}}^{1}\{v(m)_{t_{l}^{m},t}^{1}-v_{t_{l}^{m},t}^{1}\}+g_{t_{k}^{m},t_{l}^{m}}^{1}\{v(m)_{t_{l}^{m},t}^{1}-v_{t_{l}^{m},t}^{1}\}=:\sum_{i=1}^{6} B_{i}.
\end{eqnarray}
Since $B_1,B_3$ has been estimated similar to \eqref{MFBMCM-RP5}. According to the definition of $v(m)$, we have that $B_4=0$. Then, we estimate $B_6$ as below,
\begin{eqnarray}\label{MFBMCM-RP7}
|B_{6}|&\le&\|u\|_{\alpha-\operatorname{hld}}(t_{l}^{m}-t_{k}^{m})^{\alpha}\{|v(m)_{t_{l}^{m},t}^{1}|+|v_{t_{l}^{m},t}^{1}|\}\cr
&\le& C_2\|u\|_{\alpha-\operatorname{hld}}(t_l^m-t_k^m)^\alpha(t-t_l^m)^{\frac{1}{2}}\cr
&\le& C_2(2^{-m})^{\kappa}\|u\|_{\alpha-\operatorname{hld}}(t-s)^{\alpha+\frac{1}{2}-\kappa}.
\end{eqnarray}
$B_5$ also could be estimated by following same approach.

For $B_2$, directly, we have
\begin{eqnarray}\label{MFBMCM-RP8}
B_2=\int_{t_k^m}^{t_l^m}\big[\sum_{j=k+1}^l[u_\cdot-u_{t_k^m}]_j^m\mathbf{1}_{[t_{j-1}^m,t_j^m]}(r)-(u_r-u_{t_k^m})\big]dv_r.
\end{eqnarray}
Then, by leveraging the result that the integrand could be dominated by $\|u\|_\alpha(2^{-m})^\alpha$, we have
\begin{eqnarray}\label{MFBMCM-RP8}
|B_2|\leq C_3\|u\|_{\alpha-\operatorname{hld}}(2^{-m})^\alpha(t_l^m-t_k^m)^{\frac{1}{2}}\leq C_3(2^{-m})^\kappa\|u\|_{\alpha-\operatorname{hld}}(t-s)^{\alpha+\frac{1}{2}-\kappa}
\end{eqnarray}
where $0<\kappa<\alpha$ and $t-s\ge 2^{-m}$.

Thanks to the following estimate,
\begin{eqnarray}\label{MFBMCM-RP10}
\|I[u(m),v(m)-v]\|_{\alpha+\gamma-2\kappa-\operatorname{hld}} \le \varepsilon_{m}\|u\|_{\alpha}(t-s)^{\alpha+\frac{1}{2}-\kappa}
\end{eqnarray}
where $\varepsilon_{m}$ tends to $0$ when $m\to \infty$. Then, \eqref{MFBMCM-RP4} has been proved. Note that $\kappa$ is small enough, so $\alpha-\kappa$ and $\gamma-\kappa$ are arbitrarily close to $\alpha$ and $\gamma$.

Due to the result that $E_{s,t}^2(m)=(b^{H}(m)+u(m))_{s,t}\otimes(w(m)+v(m))_{s,t}-E_{s,t}^1(m)$, we could show that \eqref{MFBM+CM-RP9} holds. Then, $\{\hat S^3(b^H(m)+u(m), w(m)+v(m))\}_{m\in \mathbb{N}}$ converges to $(\mathcal{T} \Sigma^1,\mathcal{T} \Sigma^2,\mathcal{T}_{s,t}^{U,3}(B^H))\in G\hat{\Omega}_{\alpha-\kappa,\gamma-\kappa}( \mathbb{R}^{d+e})$, then $\mathcal{T}\Sigma:=(\mathcal{T} \Sigma^1,\mathcal{T} \Sigma^2,\mathcal{T} \Sigma^3)\in G\Omega_{\alpha-\kappa}(\mathbb{R}^{d+e})$. Note that $\kappa$ is small enough, so $\alpha-\kappa$ and $\gamma-\kappa$ are arbitrarily close to $\alpha$ and $\gamma$.
The proof is completed. \qed

The above result concerning GRP which is lifted from the translation of mixed FBM in the C-M direction coincides with \cite[Lemma 15.58]{2010Friz}. Nevertheless, in contrast to \cite[Lemma 15.58]{2010Friz}, the Gaussianity of the first element in the mixed process is not used in our proof. As a consequence, the result remains valid in the more general case where the first element in the mixed process need not be Gaussian.
Let $(\Omega, \mathcal{F}, \mathbb{P})$ be a probability space. Let $w=(w_t)_{0\le t\le T}$ be a standard $e$-dimensional Brownian motion, and denote by $w+v$ its translation in the direction of C-M elements $v\in \mathcal{H}^{\frac{1}{2},e}$. Let $B=\{(B_{s,t}^1, B_{s,t}^2, B_{s,t}^3)\}_{0\le s\le t\le T}$ be an $\mathcal{G}_\Omega^\beta(\mathbb{R}^d)$-valued random RP defined on $(\Omega, \mathcal{F}, \mathbb{P})$ for each $\beta\in(1/4, \alpha)$, and suppose that $B$ is independent of $w$. We assume that $\vertiii{B}_\alpha$ possesses moments of all orders. Let $\{\mathcal{F}_t\}_{0\le t\le T}$ be a usual filtration satisfying the following two properties: (i) $w$ is an $\{\mathcal{F}_t\}$-Brownian motion; and (ii) the mapping $t\mapsto (B_{0,t}^1, B_{0,t}^2, B_{0,t}^3)$ is $\{\mathcal{F}_t\}$-adapted.
\begin{rem}\label{extension}
Let $1/4 < \alpha \le 1/3$ and $1/3 < \gamma \le 1/2 $. Then, the mixed process $(B,w+v)$ could be lifted to anisotropic GRP $((B,W+V)^1,(B,W+V)^2, B^{3})\in G\hat{\Omega}_{\beta,\gamma}( \mathbb{R}^{d+e})$. Directly, it also could be lifted to GRP $(B,W+V):=((B,W+V)^1,(B,W+V)^2,(B,W+V)^3):=\mathbf{Ext} ((B,W+V)^1,(B,W+V)^2, B^{3})\in G\Omega_{\beta}(\mathbb{R}^{d+e})$.
\end{rem}

\section{Slow-fast RDE with GRP lifted by mixed FBM ($1/4<H<1/3$)}\label{sec-5}
In this section, we introduce our main system. 

Consider the slow-fast RDE in time interval $[0,T]$:
\begin{eqnarray}\label{5-1-1}
\left
\{
\begin{array}{ll}
dX^{\varepsilon, \delta}_t = f(X^{\varepsilon, \delta}_t, Y^{\varepsilon, \delta}_t)dt + \sqrt \varepsilon \sigma( X^{\varepsilon, \delta}_t)dB^H_{t},\\
dY^{\varepsilon, \delta}_t =\frac{1}{\delta} F( X^{\varepsilon, \delta}_t, Y^{\varepsilon, \delta}_t)dt + \frac{1}{{\sqrt \delta}}G( X^{\varepsilon, \delta}_t, Y^{\varepsilon, \delta}_t)dW_{t},\\
(X^{\varepsilon, \delta}_0, Y^{\varepsilon, \delta}_0)=(x, y)\in \mathbb{R}^{{m}}\times \mathbb{R}^{{n}}.
\end{array}
\right.
\end{eqnarray}
Here, the GRP $(B^H,W)$ is lifted from the mixed FBM $(b^H,w)$ with Hurst parameter $H\in (1/4,1/3)$. $\varepsilon, \delta\in(0,1]$ are small parameters which will be assumed tends to zero. The coefficients $f:\mathbb{R}^{{m}}\times \mathbb{R}^{{n}}\to\mathbb{R}^{{m}}$, $F:\mathbb{R}^{{m}}\times \mathbb{R}^{{n}}\to\mathbb{R}^{{n}}$, {$\sigma: \mathbb{R}^{{m}}\to \mathbb{R}^{{m}\times{d}}$ and $G: \mathbb{R}^{{m}}\times \mathbb{R}^{{n}}\to \mathbb{R}^{{n}\times{e}} $} are nonlinear functions. When we assume that $\sigma, G \in \mathcal{C}^4$ and $f,F$ are locally Lipschitz continuous, it deduces from \cite[Remark 3.4]{2025Inahama} that the RDE \eqref{5-1-2} has a unique local solution.

We write $V^{\varepsilon,\delta}:=(X^{\varepsilon,\delta},Y^{\varepsilon,\delta})$. Then, the slow-fast RDE \eqref{5-1-1} can be rewritten as following: for $t\in[0,T]$, 
\begin{eqnarray}\label{5-1-2}
V_t^{\varepsilon,\delta}=V_0+\int_0^t F_{\varepsilon,\delta}\big(V_s^{\varepsilon,\delta}\big) d s+\int_0^t G_{\varepsilon,\delta}\big(V_s^{\varepsilon,\delta}\big) d\Sigma^{\sqrt{\varepsilon}}_s,
\end{eqnarray}
with 
\begin{eqnarray}\label{5-1-3}
\big(V^{\varepsilon,\delta}\big)_t^{\dagger}=G_{\varepsilon,\delta}\big(V_t^{\varepsilon,\delta}\big), \quad \big(V^{\varepsilon,\delta}\big)_t^{\dagger\dagger}=(\nabla_VG_{\varepsilon,\delta}G_{\varepsilon,\delta})(V_t^{\varepsilon,\delta})
\end{eqnarray}
with the initial value $V_0=(X_0,Y_0)$ and
\begin{eqnarray*}\label{5-1-4}
F_{\varepsilon,\delta}(x, y)=\left(\begin{array}{c}
f(x, y) \\
\delta^{-1} F(x, y)
\end{array}\right), \quad G_{\varepsilon,\delta}(x, y)=\left(\begin{array}{cc}
\sigma(x) & O \\
O & (\varepsilon\delta)^{-1 / 2} G(x, y)
\end{array}\right) .
\end{eqnarray*}
Here, $\Sigma^{\sqrt{\varepsilon}}=(\sqrt{\varepsilon}\Sigma^1,\varepsilon\Sigma^2, \varepsilon^{3/2}\Sigma^3)\in \Omega_{\alpha}(\mathbb{R}^{d+e})$ is the dilation of $\Sigma=(\Sigma^1,\Sigma^2,\Sigma^3)\in \Omega_{\alpha}(\mathbb{R}^{d+e})$, which is defined in \propref{MFBM-RP}.
Then, $(V^{\varepsilon,\delta},(V^{\varepsilon,\delta})^{\dagger},(V^{\varepsilon,\delta})^{\dagger\dagger})\in \mathcal{Q}_{\Sigma^{\sqrt{\varepsilon}}}^\beta([a, b], \mathbb{R}^{m+n})$ with $1/4<\beta<\alpha<1/3$ is a controlled RP.

Then, we set 
\begin{eqnarray}\label{5-1-9}
\tilde F(x,y):= F(x,y)+\frac{1}{2}\{\sum_{i=1}^{n}\sum_{j=1}^{e}\frac{\partial G^{kj}}{\partial y_j}G^{ij}(x,y)\}_{1\le k\le n}.
\end{eqnarray}
We give the It\^o-Stratonovich correction term when $x$ is viewed as a parameter.
Define $\tau_N^\varepsilon=inf\{t\ge 0| |V_t^\varepsilon|\ge N\}$ for each $N\in \mathbb{N}$ and $\tau_\infty^{\varepsilon}=\lim_{N\to \infty}\tau_N^\varepsilon$. In fact, $Y^{\varepsilon,\delta}$ coincide with the It\^o SDE as below up to the explosion time $\tau_\infty^{\varepsilon}$ \cite[Proposition 5.5]{2025Inahama},
\begin{eqnarray}\label{5-1-10}
Y_t^{\varepsilon,\delta}=Y_0+\frac{1}{\delta} \int_0^{t} \tilde F\left(X_s^{\varepsilon,\delta}, Y_s^{\varepsilon,\delta}\right) d s+\frac{1}{\sqrt{\delta}} \int_0^{t} G\left(X_s^{\varepsilon,\delta}, Y_s^{\varepsilon,\delta}\right) d^{\mathrm{I}} w_s.
\end{eqnarray}

To ensure the existence and uniqueness of the global solution to the RDE \eqref{5-1-1}, we impose the following further conditions.
\begin{itemize}
\item[\textbf{A1}.] $\sigma\in \mathcal{C}_b^4$, $G\in \mathcal{C}^4$.
\item[\textbf{A2}.] Assume that $f$ is globally bounded. Further assume that there exists a constant $L> 0$ such that for any $ (x_1,y_1) $, $ (x_2,y_2)\in \mathbb{R}^{m} \times\mathbb{R}^{n}$, 
\begin{equation*}
|f(x_1, y_1)-f(x_2, y_2)|+ |\tilde F(x_1, y_1)-\tilde F(x_2, y_2)|+\left|G\left(x_1, y_1\right)-G\left(x_2, y_2\right)\right|\leq L(|x_1-x_2|+|y_1-y_2|), 
\end{equation*} 
hold.
\end{itemize} 
Then, the following results comes from \cite[Proposition 5.6]{2025Inahama}, the details of the proof are omitted. 
\begin{prop}\label{solution}
Let $H\in(\frac{1}{4}, \frac{1}{3})$ and $\beta\in(\frac{1}{4}, H)$ and assume \textbf{(A1)}--\textbf{(A2)}. Let $(V^{\varepsilon,\delta},(V^{\varepsilon,\delta})^{\dagger},(V^{\varepsilon,\delta})^{\dagger\dagger})\in \mathcal{Q}_{\Sigma^{\sqrt{\varepsilon}}}^\beta([a, b], \mathbb{R}^{m+n})$ be a unique solution to the RDE \eqref{5-1-2}. Then, wee have that $\left(X^{\varepsilon,\delta}, \sigma\left({X}^{\varepsilon,\delta}\right)\right)\in \mathcal{Q}_{B^{\varepsilon,H}}^\beta\left([0, T], \mathbb{R}^m\right)$ is a unique global solution of the RDE driven by $B^{\varepsilon,H}=(\sqrt{\varepsilon}B^{H,1},\varepsilon B^{H,2}, \varepsilon^{3/2} B^{H,3})$ as following:
$$
X_t^{\varepsilon,\delta}=X_0+\int_0^t f(X_s^{\varepsilon,\delta}, Y_s^{\varepsilon,\delta}) d s+\int_0^t \sigma(X_s^{\varepsilon,\delta}) dB^{\varepsilon,H}_s
$$
with 
\begin{eqnarray*}\label{5-1-11}
\big(X^{\varepsilon,\delta}\big)_t^{\dagger}=\sigma(X_t^{\varepsilon,\delta}),\quad \big(X^{\varepsilon,\delta}\big)_t^{\dagger\dagger}=(\nabla_x\sigma\sigma)(X_t^{\varepsilon,\delta}).
\end{eqnarray*}
for $t \in[0, T] $.
\end{prop}

Furthermore, 
for each $0<\delta, \varepsilon\leq 1$, the fast one $Y^{\varepsilon,\delta}$ satisfies the It\^o SDE \eqref{5-1-10}. Below, in order to give the averaged drift term $\bar f$ in the limit of the slow one $X^{\varepsilon,\delta}$, we consider the frozen It\^o SDE with fixed component $X$.
\begin{eqnarray}\label{fast}
d{Y}^{X,Y_0}_t = F({X}, {Y}^{X,Y_0}_t) dt +G( {X}, {Y}^{X,Y_0}_t)d\tilde w_{t}, \quad Y^{X,Y_0}_0=Y_0\in\mathbb{R}^n,
\end{eqnarray}
where $\{w_t\}_{t\in[0,T]}$ is a $d_2$-dimensional standard Brownian motion defined on another probability space $(\tilde \Omega, \tilde F, \{\tilde F\}_{t\ge0}, \mathbb{\tilde P} )$. Denote the expectation operator related to $\mathbb{\tilde P}$ by $\mathbb{\tilde E}$. Let $\{P^X_t\}_{t\in[0,T]}$ be the transition semigroup of $\{ Y^{X,Y_0}_t\}_{t\in[0,T]}$, i.e. for any bounded measurable function $\varphi:\mathbb{R}^m\to \mathbb{R}$:
$$
P_s^X \varphi(y):=\mathbb{\tilde E}[\varphi(Y_s^{X, Y_0})], \quad Y_0 \in \mathbb{R}^{n}, s \geq 0.
$$
To ensure the existence and uniqueness of invariant probability measures as well as the exponential ergodicity with respect to the above frozen equation, the following condition is assumed.
\begin{itemize}
\item[(\textbf{A3}).] Assume that there exist positive constants $C>0$ and $\beta_i >0 \,(i=1,2)$ such that for any $ (x,y_1),(x,y_2)\in \mathbb{R}^{m} \times\mathbb{R}^{n}$
\begin{equation*}
\begin{aligned}
2\langle y_1-y_2, \tilde F(x, y_1)-\tilde F(x, y_2)\rangle+|G(x, y_1)-G(x, y_2)|^2 
&\leq-\beta_1\left|y_1-y_2\right|^2, \\
2\langle y_1, \tilde F(x, y_1)\rangle+|G(x, y_1)|^2 & \leq-\beta_2|y_1|^2+C|x|^2+C
\end{aligned}
\end{equation*}
hold. 
\end{itemize} 
Following the similar arguments as in \cite[Lemma 3.6 and Lemma 3.7]{Liu2020}, we could derive the estimates for the frozen equation as below.
\begin{rem}\label{y_2}
Under Assumption {(\textbf{A3})}, for any given $X\in \mathbb{R}^m$, $Y_0\in \mathbb{R}^n$ and $t\in[0,T]$, it deduces 
$$
\mathbb{\tilde E}[|Y_t^{X, Y_0}|^2] \le e^{-\beta_2 t}|Y_0|^2+c(1+|X|^2).
$$
Moreover, for any $y_1,y_2\in \mathbb{R}^n$, it has 
$$
\mathbb{\tilde E}[|Y_t^{X, y_1}-Y_t^{X, y_2}|^2] \le e^{-\beta_2 t}|y_1-y_2|^2.
$$
\end{rem}
The Krylov-Bogoliubov argument and \remref{y_2} yield the existence of an invariant probability measure for $\{P^X_t\}_{t\in[0,T]}$ for every $X$. Morevoer, via the similar estimates as in \cite[Proposition 3.8]{Liu2020}, the exponential ergodicity for the frozen equation is given as below.
\begin{rem} \label{y}
{Suppose that {(\textbf{A1})}--{(\textbf{A3})} hold. For any given $X\in\mathbb{R}^m$ and initial value $Y_0\in \mathbb{R}^n$, the semigroup $\{P^X_t\}_{t\in[0,T]}$ has a unique invariant probability measure $\mu^X$. Furthermore, the following statements hold:}

{(a)There exists a constant $c>0$ such that 
$$\int_{\mathbb{R}^{n}}|y|^2 \mu^X(d y) \le c\left(1+|X|^2\right).$$
Here, $C>0$ is independent of $X$.}

{(b)There exists $c'> 0$ such that for any Lipschitz function $\varphi:\mathbb{R}^n \to \mathbb{R}$:
$$
\big|P_s^X \varphi(y)-\int_{\mathbb{R}^{n}} \varphi(z) \mu^X(d z)\big| \le c'(1+|X|+|Y_0|) e^{-\beta_1 s}|\varphi|_{\text {Lip }}, \quad s \geq 0,
$$
where $|\varphi|_{\text {Lip }}$ is the Lipschitz coefficient of $\varphi$ and $\beta_1 > 0$ is in Assumption (\textbf{A3}). }
\end{rem}

Then we give the averaged equation:
\begin{eqnarray}\label{5-1-12}
d \bar X_t =\bar f(X_t) dt, \quad \bar X_0=X_0,
\end{eqnarray}
with $\bar{f}(x)=\int_{\mathbb{R}^{n}}f(x, { y})\mu^{x}(d{ y})$. The above averaged system \eqref{5-1-12} is a deterministic ODE. 
Then, it suffices to show that $\bar f$ is Lipschitz continuous. Firstly, by Assumption {(\textbf{A2})} and \remref{y}, we have that for all for any $ (x_1, x_2)\in \mathbb{R}^{m} $ and initial value $Y_0\in \mathbb{R}^{n} $, 
\begin{eqnarray}\label{average}
\left|\bar{f}\left(x_1\right)-\bar{f}\left(x_2\right)\right| &\le &|\int_{\mathbb{R}^n} f(x_1, y) \mu^{x_1}(d y)-\mathbb{E}[f(x_1, Y_t^{x_1,Y_0})]| +\big|\int_{\mathbb{R}^n} f(x_2, y) \mu^{x_1}(d y)-\mathbb{E}[f(x_2, Y_t^{x_2,Y_0})]\big| \cr
&&+\big|\mathbb{E}[f(x_1, Y_t^{x_1,Y_0})]-\mathbb{E}[f(x_2, Y_t^{x_2,Y_0})]\big| \cr
&\le & C e^{-\beta_1 s}(1+|x_1|+|x_2|+|Y_0|)+L|x_1-x_2|.
\end{eqnarray}
Let $s\to \infty$, it deduces that $\bar f$ is Lipschitz continuous.

Then, the slow component $X^{\varepsilon,\delta}$ will converges to the solution to the averaged equation \eqref{5-1-12} as $\varepsilon$, $\delta$ tend to zero which is precisely shown as below, whose proof comes from \cite[Theorem 2.1]{2025Inahama} with some small extension, so details of the proof are omitted.
\begin{prop}\label{AP}
Let $H\in(\frac{1}{4}, \frac{1}{3})$, $\beta\in(\frac{1}{4}, H)$ and $\varepsilon,\delta\in (0,1)$. Under Assumption \textbf{(A1)}--\textbf{(A3)}, we have that for $\nu\ge 1$,
\begin{eqnarray}\label{5-1-13}
\mathbb{E}[\|X^{\varepsilon,\delta}-\bar X\|_{\beta-\operatorname{hld}}^\nu ]=0.
\end{eqnarray}
as $\delta$ tends to zero.
\end{prop}

\section{LDP and MDP for slow-fast RDE}\label{sec-5}
In this section, we give the LDP and MDP for the slow one in the RDE \eqref{5-1-1}. Firstly, we state our main results.

Before stating the LDP result, we first introduce a measurable map related to the slow one in the RDE \eqref{5-1-1}, that is 
\[
\mathcal{G}^{(\varepsilon, \delta)}: \mathcal{C}_0\left([0,T], \mathbb{R}^d\right) \rightarrow \mathcal{C}^{\beta-\operatorname{hld}}\left([0,T], \mathbb{R}^m\right)
\]
such that
$X^{\varepsilon,\delta}:=\mathcal{G}^{\varepsilon,\delta}(\sqrt \varepsilon b^H, \sqrt \varepsilon w)$. 
Furthermore, to study an LDP for the slow component in slow-fast RDE \eqref{5-1-1}, we assume the following conditions.
\begin{itemize}
\item[(\textbf{A4}).] Further assume that there exists a constant $L> 0$ such that for any $x\in \mathbb{R}^m$,
\begin{equation*}
\sup_{y\in\mathbb{R}^{n}}\left|G\left(x, y\right)\right| \leq L\left(1+\left|x\right|\right)
\end{equation*}
holds.
\end{itemize} 
Next, we define the skeleton equation in the rough sense as follows
\begin{eqnarray}\label{5-1-6}
d\tilde{X}_t = \bar{f}(\tilde{X}_t)dt + \sigma( \tilde{X}_t)dU_t
\end{eqnarray} 
where $\tilde{X}_0=X_0$, $U=(U^1,U^2, U^3)\in G\Omega_{\alpha}(\mathbb{R}^d)$ and $\bar{f}(x)=\int_{\mathbb{R}^{n}}f(x, { y})\mu^{x}(d{ y})$ for $x\in\mathbb{R}^m$. Since $\bar f$ is Lipschtiz continuous, it is not too difficult to see that there exists a unique global solution $(\tilde{X},\tilde{X}^{\dagger},\tilde{X}^{\dagger\dagger})\in \mathcal{Q}_{U}^\beta([0, T], \mathbb{R}^m)$ to the RDE \eqref{5-1-6}. Moreover, we have for $0<\beta<\alpha<H$ that 
$$\|\tilde{X}\|_{\beta-\operatorname{hld}}\le c,$$
with the constant $c>0$ independent of $U$. 
Therefore, we also define a map
\[
\mathcal{G}^{0}: S_{N} \rightarrow \mathcal{C}^{\beta-\operatorname{hld}}\left([0,T],\mathbb{R}^m\right)
\]
such that its solution $\tilde{X}=\mathcal{G}^{0}(u, v)$. 

\begin{rem}\label{rem3}
The above RDE \eqref{5-1-6} here also coincides with the Young ODE as below:
\begin{eqnarray}\label{5-1-7}
d\tilde{X}_t = \bar{f}(\tilde{X}_t)dt + \sigma( \tilde{X}_t)du_t,\quad \tilde{X}_t=X_0.
\end{eqnarray} 
Since $(H+1/2)^{-1}<q<2$, it has $\|(u,v)\|_{q-\operatorname{var}}<\infty$, it suffices to verify that there exists a unique pathwise solution $\check {X} \in \mathcal{C}^{p-\operatorname{var}}\left([0,T],\mathbb{R}^m\right)$ to \eqref{5-1-7} for $(u, v) \in S_{N}$. Furthermore, 
$\|\tilde{X}\|_{p-\operatorname{var}}\le c$
where the postitive constant $c$ is independent of $(u, v)$. 
\end{rem}

Now, we give the statement of our LDP result. 
\begin{thm}\label{thm}
{Let $H\in(1/4,1/3)$ and fix $1/4<\beta<H$. Assume {(\textbf{A1})}--{(\textbf{A4})} and $\delta=o(\varepsilon)$. Let $\varepsilon \to 0$, the slow component $X^{\varepsilon,\delta}$ of system (\ref{5-1-1}) satisfies a LDP on $\mathcal{C}^{\beta-\operatorname{hld}}([0,T],\mathbb{R}^{m})$ with a good rate function $I: \mathcal{C}^{\beta-\operatorname{hld}}([0,T],\mathbb{R}^m)\rightarrow [0, \infty)$}
\begin{eqnarray*}\label{rate}
I(\xi) &=& {
\inf\Big\{\frac{1}{2}\|u\|^2_{\mathcal{H}^{H,d}}~:~{ u\in \mathcal{H}^{H,d} 
\quad\text{such that} \quad\xi =\mathcal{G}^{0}(u, 0)}\Big\} 
}
\cr
&=& \inf\Big\{\frac{1}{2}\|(u,v)\|^2_{\mathcal{H}} ~:~ {( u,v)\in \mathcal{H}\quad\text{such that} \quad\xi =\mathcal{G}^{0}(u, v)}\Big\},
\qquad 
\xi\in \mathcal{C}^{\beta-\operatorname{hld}}\left([0,T], \mathbb{R}^m\right). 
\end{eqnarray*}
\end{thm} 
\begin{rem} 
The variational formula cannot be applied directly in $\mathcal{C}^{\beta-\operatorname{hld}}([0,T], \mathbb{R}^{k})$ since it is not separable. So we give the ``little H\"older space", 
$H^\beta = H^\beta ([0,T], \mathbb{R}^{k})$ is the space that for all {$\phi\in \mathcal{C}^{\beta-\operatorname{hld}}([0,T], \mathbb{R}^{k})$}, equipped with the norm
\begin{eqnarray*}\label{111}
\lim_{\delta\to0+}\sup_{\substack{|t-s|\le \delta\\0 \leq s<t \le T}} \frac{|\phi_t-\phi_s|}{(t-s)^\beta}=0. 
\end{eqnarray*}
{The space $H^\alpha$ is a separable Banach space}. Moreover, 
$
H^\beta=\overline{\bigcup_{\kappa>0}\mathcal{C}^{(\beta+\kappa)-\operatorname{hld}}}
$
with the closure being taken in the norm $\|\cdot\|_{\alpha-\operatorname{hld}}$ and $H^\beta \hookrightarrow \mathcal{C}^{\beta-\operatorname{hld}}([0,T], \mathbb{R}^{k})$ \cite{Ciesielski}.
So the variational formula {can} be used well, and it only asks to prove the weak convergence method under the $\beta$-H\"older norm. Then the conventional contraction principle works to show the LDP on the H\"older continuous space \cite[Theorem 4.2.1]{Dembo2009Large}.
\end{rem}
Then, we intend to demonstrate the MDP result. Further assumptions are necessary. Let $a:=GG^{T}$. 
\begin{itemize}
\item [(\textbf{A5}).] Assume that for every $x\in \mathbb{R}^m$, $f(x,\cdot), F(x,\cdot), a(x,\cdot)\in C^2_b$.
\item [(\textbf{A6}).] Assume that the coefficient $a$ is nondegenerate in $y$ uniformly with respect to $x$, that is, there exists $\hat L>1$ such that for any $x\in \mathbb{R}^m$, 
$$\hat L^{-1}|\xi|^2\le |a(x,y)\xi|^2\le \hat L|\xi|^2,\quad \xi \in \mathbb{R}^n.$$
\item [(\textbf{A7}).] Assume that 
$\lim_{|y| \to \infty}\sup_{ x}\langle y,F(x,y)\rangle =-\infty$.
\end{itemize} 
\begin{rem}\label{bar_f}
If Assumptions (\textbf{A2}), (\textbf{A5})--(\textbf{A7}) hold, thanks to the result \cite[Theorem 2.1]{2021Rockner}, we have that
for every $\hat f(x,y):=f(x,y)-\bar f(x)$, there exists a solution $u(x,\cdot)$ satisfying that $u(x,\cdot) \in C^3_b$ for every $x\in \mathbb{R}^m$ and $u(\cdot,y)\in C_b^1$ for every $y\in \mathbb{R}^n$, to the following Poisson equation in $ \mathbb{R}^n$:
$$\mathcal{L}(x,y)u(x,y)=\hat f(x,y),\quad y\in \mathbb{R}^n$$
where
$$\mathcal{L}(x,y): = \sum_{i,j=1}^{n}a^{ij}(x,y)\frac{\partial^2}{\partial y_i \partial y_j}+\sum_{i}^{n}\hat f^{i}(x,y)\frac{\partial}{\partial y_i}.$$
Then, according to \cite[Lemma 3.2]{2021Rockner}, it deduces that $\bar f\in C^2_b$. When $F(x,y)=F(y)$ and $G(x,y)=G(y)$, it ensures that $\bar f\in C^2_b$ by only assuming $\nabla_{x} f\left( x,y\right)$ is Lipschitz with respect to $x\in\mathbb{R}^m$ for each $y\in\mathbb{R}^n$.
\end{rem}
Consider the following RDE,
\begin{eqnarray}\label{5-1-10}
d\check{Z}_t = \nabla {f}(\bar {X}_t)\check Z_tdt + G( \bar X_t)dU_t 
\end{eqnarray}
where $\check{Z}_0=0$. Here, $\bar X$ is a deterministic differentiable path. By \propref{FBM+CM-RP}, the GRP $U:=(U^1,U^2,U^3)\in G\Omega_{\alpha}(\mathbb{R}^{d})$ is lifted from C-M elements $u\in \mathcal{H}^{H,d}$. Then, due to \remref{bar_f}, we have that under assumptions (\textbf{A2}), (\textbf{A5})--(\textbf{A7}), there exists a unique solution $(\check{Z},\check{Z}^{\dagger},\check{Z}^{\dagger\dagger})\in \mathcal{Q}_{U}^\beta([0, T], \mathbb{R}^m)$ to the RDE \eqref{5-1-10} globally. Moreover, it has for $0<\beta<\alpha<H$ that 
$\|\check{Z}\|_{\beta-\operatorname{hld}}\le c$ for constant $c>0$ which is independent of $U$. 
Therefore, we also define a map
\[
\mathcal{\check G}^{0}: S_{N} \rightarrow \mathcal{C}^{\beta-\operatorname{hld}}\left([0,T],\mathbb{R}^m\right)
\]
such that its solution $\check{Z}=\mathcal{\check G}^{0}(u, v)$. In fact, the map $\mathcal{\check G}^{0}(u, v)$ is independent of $v$. 
Besides, the above RDE \eqref{5-1-10} coincides with the Young ODE as following:
\begin{eqnarray}\label{5-1-11}
\left
\{
\begin{array}{ll}
d\check{Z}_t = \nabla \bar {f}(\bar {X}_t)\check Z_tdt + G(\bar {X}_t)du_t, \\
dX_t = f(X_t)dt,
\end{array}
\right.
\end{eqnarray}
with $(\bar X_0,\check{Z}_0)=(X_0,0)$.
Recall that $(H+1/2)^{-1}<q<2$, it has $\|u\|_{q-\operatorname{var}}<\infty$. By the Young integral theory, it suffices to see that there exists a unique solution $\check {Z} \in \mathcal{C}^{p-\operatorname{var}}\left([0,T],\mathbb{R}^d\right)$ to (\ref{5-1-11}) for $(u, v) \in S_{N}$. Moreover, we have
$\|\check{Z}\|_{p-\operatorname{var}}\le c$
where the constant $c>0$ is independent of $(u, v)$. 
Now, we give the statement of our MDP result. 
\begin{thm}\label{mdpthm}
{Let $H\in(1/4,1/3)$ and fix $1/4<\beta<H$. Assume {(\textbf{A1})}--{(\textbf{A7})} and $\delta=o(\varepsilon)$. Let $\varepsilon \to 0$, the slow component $X^{\varepsilon,\delta}$ of system (\ref{5-1-1}) satisfies a MDP on $\mathcal{C}^{\beta-\operatorname{hld}}([0,T],\mathbb{R}^{m})$ with a good rate function $I: \mathcal{C}^{\beta-\operatorname{hld}}([0,T],\mathbb{R}^m)\rightarrow [0, \infty)$}
\begin{eqnarray*}\label{rate}
I(\xi) &=& {
\inf\Big\{\frac{1}{2}\|u\|^2_{\mathcal{H}^{H,d}}~:~{ u\in \mathcal{H}^{H,d} 
\quad\text{such that} \quad\xi =\mathcal{G}^{0}(u, 0)}\Big\} 
}
\cr
&=& \inf\Big\{\frac{1}{2}\|(u,v)\|^2_{\mathcal{H}} ~:~ {( u,v)\in \mathcal{H}\quad\text{such that} \quad\xi =\mathcal{\check G}^{0}(u, v)}\Big\},
\qquad 
\xi\in \mathcal{C}^{\beta-\operatorname{hld}}\left([0,T], \mathbb{R}^m\right). 
\end{eqnarray*}
\end{thm}

We now comment on the reason for the assumption $\lim_{\varepsilon\to 0}\delta/\varepsilon=0$.
\begin{rem}\label{delta}
The relationship between the time-scale parameter $\delta$ and the noise intensity $\varepsilon$ converge to zero can be characterized by the following three regimes,
$$\lim_{\varepsilon\to0}\frac\delta\varepsilon=\left\{
\begin{array}
{cl}0,\quad\quad\quad~~ & \mathrm{Regime}~1; \\
\gamma\in(0,\infty), & \mathrm{Regime}~2; \\
\infty, \quad\quad\quad~& \mathrm{Regime}~3.
\end{array}\right.
$$
Here we assume that $\delta$ decays faster than $\varepsilon$, which is corresponding to Regime 1. 
In this work, we establish the LDP and MDP by applying the weak convergence method and the variational representation for mixed FBM with $H\in(\frac{1}{4}, \frac{1}{3})$. This approach requires a quantitative analysis of the dynamical behavior of the controlled slow component, which is perturbed by the controlled fast variable. However, no result shows that there exists a unique invariant measure for the controlled equation. To overcome this, we introduce a ``replaced" fast equation that admits a unique invariant measure.

In the other two regimes, the long-term behavior of the controlled fast equation depends non-trivially on the control term and the $y$-marginal of the limiting occupation measure is no longer decoupled from the control term $v$, so the above averaging methods break down. When $H> 1/2$, the LDP for the three regimes was established by using the viable pair method \cite{2025Gailus}, which is under the Young integral setting. However, for $H<1/2$, the situation poses significant challenges, as the LDP in the other two regimes would need to be formulated within the RP framework. It is unclear whether the viable pair method remains applicable in this RP setting, and this is what we will to explore in future work.
\end{rem}

\subsection{Proof of the LDP}
Firstly, let $(u^{\varepsilon, \delta}, v^{\varepsilon, \delta})\in \mathcal{A}_{b}$. In order to apply the variational representation \eqref{2-1}, we give the following controlled slow-fast RDE associated to the original slow-fast component $({X}^{\varepsilon, \delta}, {Y}^{\varepsilon, \delta})$.
\begin{eqnarray}\label{2}
\left
\{
\begin{array}{ll}
d\tilde {X}^{\varepsilon, \delta}_t =& f(\tilde {X}^{\varepsilon, \delta}_t, \tilde {Y}^{\varepsilon, \delta}_t)dt + \sigma(\tilde {X}^{\varepsilon, \delta}_t)d[\mathcal{T}_t^{u^{\varepsilon, \delta}}(B^{\varepsilon,H}) ]\\
d\tilde {Y}^{\varepsilon, \delta}_t =& \frac{1}{\delta}F( \tilde {X}^{\varepsilon, \delta}_t, \tilde {Y}^{\varepsilon, \delta}_t)dt +\frac{1}{\sqrt{\delta\varepsilon}}G(\tilde {X}^{\varepsilon, \delta}_t, \tilde {Y}^{\varepsilon, \delta}_t)dv^{\varepsilon, \delta}_t + \frac{1}{\sqrt{\delta}}G(\tilde {X}^{\varepsilon, \delta}_t, \tilde {Y}^{\varepsilon, \delta}_t)dw_{t}.
\end{array}
\right.
\end{eqnarray}
Here, $\mathcal{T}^{u^{\varepsilon, \delta}}(B^{\varepsilon,H})$ is an GRP which has already been defined in \propref{FBM+CM-RP}.
Here, $(u^{\varepsilon, \delta}, v^{\varepsilon, \delta})\in \mathcal{A}_{b}$ is called a pair of control. 

We divide $[0,T]$ into subintervals of equal length $\Delta$. For $t\in[0,T]$, we set $t(\Delta)=\left\lfloor\frac{t}{\Delta}\right\rfloor \Delta$, which denotes the nearest breakpoint preceding $t$. Next, we construct the auxiliary process as below:
\begin{eqnarray}\label{4-3}
d\hat {Y}^{\varepsilon, \delta}_t &=& \frac{1}{\delta} F( \tilde {X}^{\varepsilon, \delta}_{t(\Delta)}, \hat {Y}^{\varepsilon, \delta}_t)dt + \frac{1}{{\sqrt \delta }}G(\tilde {X}^{\varepsilon, \delta}_{t(\Delta)}, \hat {Y}^{\varepsilon, \delta}_t)dw_{t}
\end{eqnarray}
with $\hat {Y}^{\varepsilon, \delta}_0=Y_0$.

To prove \thmref{thm}, a-priori estimates should be given. 
\begin{lem}\label{lem1}
Assume (\textbf{A1})--(\textbf{A4}) and let $\nu\ge 1$ and $N\in\mathbb{N}$. Then, for all $\varepsilon,\delta\in(0,1]$, we have
\begin{equation}\label{4-6}
\mathbb{E}\big[\|\tilde X^{\varepsilon,\delta}\|^\nu_{\beta-\operatorname{hld}}\big] \le C_{\nu,N}.
\end{equation}
Here, $C_{\nu,N}$ is a constant.
\end{lem}
\para{Proof}. Now, $(\tilde {X}^{\varepsilon,\delta},(\tilde {X}^{\varepsilon,\delta})^{\dagger},(\tilde {X}^{\varepsilon,\delta})^{\dagger\dagger})\in \mathcal{Q}_{\mathcal{T}\Sigma^{\sqrt{\varepsilon}}, [0,T]}^\beta$ 
satisfies the following RDE driven by $\mathcal{T}\Sigma^\varepsilon$:
\begin{eqnarray}\label{4-2-4}
\tilde {X}^{\varepsilon,\delta}_t &=&X_0+ \int_{0}^{t}f(\tilde {X}^{\varepsilon,\delta}_s, \tilde {Y}^{\varepsilon,\delta}_s)ds + \int_{0}^{t}\sigma(\tilde {X}^{\varepsilon,\delta}_t)d[\mathcal{T}_t^{u^{\varepsilon, \delta}}(B^{\sqrt{\varepsilon},H}) ],\cr 
(\tilde {X}^{\varepsilon,\delta}_t)^{\dagger}&=&\sigma(\tilde {X}^{\varepsilon,\delta}_t),\cr
(\tilde X^{\varepsilon,\delta})_t^{\dagger\dagger}&=&(\nabla\sigma\cdot \sigma)( \tilde X^{\varepsilon,\delta}_t).
\end{eqnarray}
For every $(\tilde {X}^{\varepsilon,\delta},(\tilde {X}^{\varepsilon,\delta})^{\dagger},(\tilde {X}^{\varepsilon,\delta})^{\dagger\dagger})\in \mathcal{Q}_{\mathcal{T}_t^{u^{\varepsilon, \delta}}(B^{\sqrt{\varepsilon},H}) , [0,T]}^\beta$, we observe that the right hand side of \eqref{4-2-4} also belongs to $ \mathcal{Q}_{\mathcal{T}_t^{u^{\varepsilon, \delta}}(B^{\sqrt{\varepsilon},H}) , [0,T]}^\beta$.
{Set $\tilde {X}^{\varepsilon,\delta}_{s,t}=\tilde {X}^{\varepsilon,\delta}_t-\tilde {X}^{\varepsilon,\delta}_s$. Let $\tau \in[0,T]$ and set $$\mathcal{M}_{[0, \tau]}^1,\mathcal{M}_{[0, \tau]}^2:\mathcal{Q}_{\mathcal{T}_t^{u^{\varepsilon, \delta}}(B^{\sqrt{\varepsilon},H}) }^\beta([0, \tau], \mathbb{R}^m)\mapsto \mathcal{Q}_{\mathcal{T}_t^{u^{\varepsilon, \delta}}(B^{\sqrt{\varepsilon},H}) }^\beta([0, \tau], \mathbb{R}^m)$$ by
\begin{eqnarray}\label{2-42}
\mathcal{M}_{[0, \tau]}^1(\tilde X^{\varepsilon,\delta}, (\tilde X^{\varepsilon,\delta})^{\dagger},(\tilde X^{\varepsilon,\delta})^{\dagger\dagger})&=&\left(\int_0^{\cdot} \sigma(\tilde X^{\varepsilon,\delta}_s)d\mathcal{T}_s^{u^{\varepsilon, \delta}}(B^{\sqrt{\varepsilon},H}), \nabla\sigma( \tilde X^{{\varepsilon,\delta}}),(\nabla\sigma\cdot \sigma)( \tilde X^{\varepsilon,\delta})\right),\cr
\mathcal{M}_{[0, \tau]}^2(\tilde X^{\varepsilon,\delta}, (\tilde X^{\varepsilon,\delta})^{\dagger},(\tilde X^{\varepsilon,\delta})^{\dagger\dagger})&=&\left(\int_0^{\cdot} f(\tilde X^{\varepsilon,\delta}_s,\tilde Y^{\varepsilon,\delta}_s) ds, 0,0\right)
\end{eqnarray}
and $(\tilde {X}^{{\varepsilon,\delta}},(\tilde {X}^{{\varepsilon,\delta}})^{\dagger},(\tilde {X}^{\varepsilon})^{\dagger\dagger}):=\mathcal{M}_{[0, \tau]}^{\xi}:=(\xi, 0,0)+\mathcal{M}_{[0, \tau]}^1+\mathcal{M}_{[0, \tau]}^2$. Moreover, we stress the fact that the fixed point of $\mathcal{M}_{[0, \tau]}^{\xi}$ is the solution to the \eqref{2-39} on the time interval $[0,\tau]$ for $0<\tau\le T$. We also set
\begin{eqnarray*}\label{4-2-5}
B^{X_0}_{0,\tau}=\{(\tilde {X}^{{\varepsilon,\delta}},(\tilde {X}^{{\varepsilon,\delta}})^{\dagger},(\check {X}^{{\varepsilon,\delta}})^{\dagger\dagger})\in \mathcal{Q}^\beta_{\mathcal{T}^{u^{\varepsilon, \delta}}(B^{\sqrt{\varepsilon},H}),[0, \tau]}|\|(\tilde {X}^{{\varepsilon,\delta}},(\check {X}^{{\varepsilon,\delta}})^{\dagger},(\tilde {X}^{{\varepsilon,\delta}})^{\dagger\dagger})\|_{\mathcal{Q}^\beta_{\mathcal{T}^{u^{\varepsilon, \delta}}(B^{\sqrt{\varepsilon},H}),[0, \tau]}}\le 1 \}.
\end{eqnarray*}
In fact, the above ball is of radius 1 centered at $t \mapsto (X_0+\sigma(X_0)\mathcal{T}^{u^{\varepsilon, \delta},1}(B^{\sqrt{\varepsilon},H})_{0,\tau}+ \nabla\sigma( X_0)\mathcal{T}^{u^{\varepsilon, \delta},2}(B^{\sqrt{\varepsilon},H})_{0,\tau}, \sigma(X_0)+ \nabla\sigma\cdot \sigma( X_0)\mathcal{T}^{u^{\varepsilon, \delta},1}(B^{\sqrt{\varepsilon},H})_{0,\tau},\nabla\sigma\cdot \sigma_{1}( X_0)).$
From now, we only focus on the time interval $[0,\tau]$. By some direct computation and Assumption (\textbf{H1}), it deduces that for all $(\tilde {X}^{{\varepsilon,\delta}},(\tilde {X}^{{\varepsilon,\delta}})^{\dagger},(\tilde {X}^{\varepsilon})^{\dagger\dagger})\in B^{X_0}_{0,\tau}$,
\begin{eqnarray*}\label{4-2-6}
\|(\tilde {X}^{{\varepsilon,\delta}})^{\dagger}\|_{\infty,[0,\tau]} \le |\sigma(X_0)|+\sup_{0 \le s\le \tau}|(\tilde {X}^{{\varepsilon,\delta}})^{\dagger}_s-(\tilde {X}^{{\varepsilon,\delta}})^{\dagger}_0|
\le K+ \|(\tilde {X}^{{\varepsilon,\delta}})^{\dagger}\|_{\beta-\operatorname{hld},[0,\tau]}\tau^{\beta} 
\le K+1.
\end{eqnarray*}
where the constant $K:=\|\sigma\|_{C_b^3} \vee\|f\|_{\infty} \vee L$ where $L$ is defined in (\textbf{H2}).}
Similarly, we have
\begin{eqnarray*}\label{4-2-7}
\|(\tilde {X}^{{\varepsilon,\delta}})^{\dagger\dagger}\|_{\infty,[0,\tau]} &\le& | (\nabla\sigma\cdot\sigma)( X_0)|+\sup_{0 \le s\le \tau}|(\tilde {X}^{{\varepsilon,\delta}})^{\dagger\dagger}_s-(\tilde {X}^{{\varepsilon,\delta}})^{\dagger\dagger}_0|\cr
&\le& K+ \|(\tilde {X}^{{\varepsilon,\delta}})^{\dagger\dagger}\|_{\beta-\operatorname{hld},[0,\tau]}\tau^{\beta} \cr
&\le& K+1.
\end{eqnarray*}
Then, with aid of \propref{MFBM+CM-RP}, it has
\begin{eqnarray}\label{4-2-8}
|(\tilde {X}^{{\varepsilon,\delta}})^\dagger_{s,t}|&\le& |(\tilde {X}^{{\varepsilon,\delta}})^{\dagger\dagger}_{s}\mathcal{T}^{u^{\varepsilon, \delta},1}(B^{\sqrt{\varepsilon},H})_{s,t}|+|(\tilde{X}^{\varepsilon})^{\sharp\sharp}_{s,t}|\cr
&\le& (K+1)\|\mathcal{T}^{u^{\varepsilon, \delta},1}(B^{\sqrt{\varepsilon},H})\|_{\alpha-\operatorname{hld}} (t-s)^{\alpha}+\|(\tilde{X}^{\varepsilon, \delta})^{\sharp\sharp}\|_{3\beta-\operatorname{hld},[0,\tau]} (t-s)^{3\beta}\cr
&\le& (K+1)(\|\mathcal{T}^{u^{\varepsilon, \delta},1}(B^{\sqrt{\varepsilon},H})\|_{\alpha-\operatorname{hld}}+1) (t-s)^{\alpha}.
\end{eqnarray}
Then, with the above estimation, it deduces that 
\begin{eqnarray}\label{4-2-9}
|\check {X}^{{\varepsilon,\delta}}_{s,t}|&\le& |(\tilde {X}^{{\varepsilon,\delta}})^{\dagger}_{s}\mathcal{T}^{u^{\varepsilon, \delta},1}(B^{\sqrt{\varepsilon},H})_{s,t}|+|(\tilde {X}^{{\varepsilon,\delta}})^{\dagger\dagger}_{s}\mathcal{T}^{u^{\varepsilon, \delta},2}(B^{\sqrt{\varepsilon},H})_{s,t}|+|(\check{X}^{{\varepsilon,\delta}})^\sharp_{s,t}|\cr
&\le& (K+1)(\|\mathcal{T}^{u^{\varepsilon, \delta},1}(B^{\sqrt{\varepsilon},H})\|_{\alpha-\operatorname{hld}}+1) (t-s)^{\alpha}
\cr
&&+(K+1)\|\mathcal{T}^{u^{\varepsilon, \delta},2}(B^{\sqrt{\varepsilon},H})\|_{\alpha-\operatorname{hld}} (t-s)^{2\alpha}+\|(\tilde{X}^{\varepsilon, \delta})^\sharp\|_{3\beta-\operatorname{hld},[0,\tau]} (t-s)^{3\beta}\cr
&\le&(K+1)(\|\mathcal{T}^{u^{\varepsilon, \delta},1}(B^{\sqrt{\varepsilon},H})\|_{\alpha-\operatorname{hld}}+1) (t-s)^{\alpha}\cr
&&+ (K+1)(\|\mathcal{T}^{u^{\varepsilon, \delta},2}(B^{\sqrt{\varepsilon},H})\|_{2\alpha-\operatorname{hld}}+1) (t-s)^{2\alpha}.
\end{eqnarray}
By taking $\tau$ is small enough, it is enough to show that $\mathcal{M}$ maps $B^{X_0}$ to itself and is a contraction on $B^{X_0}$. For sake of the simplicity, we choose suitable $\tau <\lambda:= \{C_\beta(K+1)^{\hat \nu}(\vertiii{\mathcal{T}^{u^{\varepsilon, \delta}}(B^{\sqrt{\varepsilon},H})}_{\alpha-\operatorname{hld}}+1)^{\hat \nu}\}^{-1/(\alpha-\beta)}<1$ with some constants ${\hat \nu}>2$ and $C_\beta>4$. Then, then $\beta$-H\"older norm of $\tilde {X}^{\varepsilon}$ on subinterval $[0,\tau]$ can be dominated by $1$. Then we estimate $\|\tilde {X}^{\varepsilon}\|_{\beta-\operatorname{hld}}=\|\tilde {X}^{\varepsilon}\|_{\beta-\operatorname{hld},[0,T]}$. By the H\"older inequality and the fact that there are $\lfloor\frac{T}{\lambda}\rfloor+1$ subintervals on $[0,T]$, we have
\begin{eqnarray}\label{4-2-10}
\|\tilde {X}^{{\varepsilon,\delta}}\|_{\beta-\operatorname{hld}}\le (\lfloor\frac{T}{\lambda}\rfloor+1)^{1-\beta}\le c_{\alpha,\beta}\{(K+1){(\vertiii{\mathcal{T}^{u^{\varepsilon, \delta}}(B^{\sqrt{\varepsilon},H})}_{\alpha-\operatorname{hld}}+1)}\}^{\iota}
\end{eqnarray}
for constants $c_{\alpha,\beta}$ and $\iota>0$ only depending on $\alpha$ and $\beta$. Then, for all $\nu\ge 1$, by taking expectation of $\nu$-moments of \eqref{4-2-9} and the property that $\mathbb{E}[\vertiii{\mathcal{T}^{u^{\varepsilon, \delta}}(B^{\sqrt{\varepsilon},H})}_{\alpha-\operatorname{hld}}^\nu]<\infty$ for every $1/4<\alpha<H$ and all $\nu\ge 1$, the estimate \eqref{4-6} is arrived. The proof is completed.\qed

Following the similar arguments as in \cite[Lemma 4.3]{2025Yang} and \cite[Lemma 4.5]{2025Yang}, we could derive the below estimates.
\begin{lem}\label{lem3}
Assume (\textbf{A1})--(\textbf{A4}) and let $N\in\mathbb{N}$, we have 
$$\int_{0}^{T} \mathbb{E}\big[|\tilde Y_t^{\varepsilon,\delta}|^2\big] dt\le C_{\alpha, \beta, N},$$
furthermore, 
$$\int_{0}^{T}\mathbb{E}\big[|\tilde Y_{t}^{\varepsilon,\delta}-\hat Y_{t}^{\varepsilon,\delta}|^2\big]dt \le C_{\alpha, \beta, N}(\frac{\sqrt{\delta}}{{\sqrt{\varepsilon}}}+\Delta^{2\beta}).$$
Here, $C_{\alpha, \beta, N}>0$ is a constant.
\end{lem}

Due to the equivalence between the LDP and the Laplace principle, it turns to show the Laplace principle. The convenient and sufficient condition to prove the Laplace principle is given below, whose proof could refer to \cite[Theorem 3.1]{2023Inahama}.
\begin{prop}\label{weak_convergence}
If the following two claims hold:

(\textbf{Condition (i)}) Let $(u^{(j)}, v^{(j)}),(u, v)\in S_N$ such that $(u^{(j)}, v^{(j)})\rightarrow(u, v)$ as $j\rightarrow\infty$ with the weak topology in $\mathcal{H}$. 
In this step, we will prove that 
\begin{eqnarray} \label{3.26}
\mathcal{G}^{0}( u^{(j)} ,v^{(j)} )\rightarrow\mathcal{G}^{0}(u,v) \label{step2}
\end{eqnarray}
in $\mathcal{C}^{\beta-\operatorname{hld}}([0,T],\mathbb{R}^{m})$ as $j \to \infty$.

(\textbf{Condition (ii)}) Assume $(u^{\varepsilon,\delta}, v^{\varepsilon,\delta})\in \mathcal{A}^{N}_b$ such that $(u^{\varepsilon,\delta}, v^{\varepsilon,\delta})$ weakly converges to $(u, v)$ as $\varepsilon \to 0$. 
In this step, we will prove that $\tilde {X}^{\varepsilon,\delta}$ weakly converges to $\tilde {X}$ in $\mathcal{C}^{\beta-\operatorname{hld}}([0,T],\mathbb{R}^{m})$ as $\varepsilon\rightarrow 0$, that is,
\begin{eqnarray} \label{step3}
\mathcal{G}^{(\varepsilon,\delta)}(\sqrt \varepsilon b^H+u^{\varepsilon,\delta}, \sqrt \varepsilon w+v^{\varepsilon,\delta})\xrightarrow{\textrm{weakly}}\mathcal{G}^{0}(u , v) \quad \text{as }\varepsilon\rightarrow 0.
\end{eqnarray} 
Then for every bounded and continuous function $\Phi: \mathcal{C}^{\beta-\operatorname{hld}}([0,T],\mathbb{R}^m)\to \mathbb{R}$, we have that the Laplace lower bound 
\begin{eqnarray}\label{5-36}
\liminf _{\varepsilon \rightarrow 0} -\varepsilon \log \mathbb{E}\big[e^{-\frac{\Phi(X^{\varepsilon,\delta})}{\varepsilon} }\big] \geq\inf _{\psi:=\mathcal{G}^0(u,v) \in \mathcal{C}^{{\beta-\operatorname{hld}}}([0,T],\mathbb{R}^m)}\left[\Phi(\psi)+I(\psi)\right]
\end{eqnarray}
and the Laplace upper bound
\begin{eqnarray}\label{5-37}
\limsup _{\varepsilon \rightarrow 0} -\varepsilon\log \mathbb{E}\big[e^{-\frac{\Phi(X^{\varepsilon,\delta})}{\varepsilon} }\big] \le\inf _{\psi:=\mathcal{G}^0(u,v) \in \mathcal{C}^{{\beta-\operatorname{hld}}}([0,T],\mathbb{R}^m)}\left[\Phi(\psi)+I(\psi)\right]
\end{eqnarray}
hold and the goodness of rate function $I$. 
\end{prop}

\para{Proof of \thmref{thm}}.
We will prove that \textbf{Condition (i)} and \textbf{Condition (ii)} in \propref{weak_convergence} hold. 

\textbf{Verification of \textbf{Condition (i)}}. The skeleton equation satisfies the RDE as follows
\begin{eqnarray}\label{4-2-12}
d\tilde{X}^{(j)}_t = \bar{f}(\tilde{X}^{(j)}_t)dt +\sigma( \check{X}^{(j)}_t)dU^{(j)}_t
\end{eqnarray} 
where $\tilde{X}^{(j)}_t=X_0$, $U^{(j)}\in G\Omega_{\alpha}(\mathbb{R}^{d_1})$ is GRP lifted by $u\in\mathcal{H}^{H,d_1}$. Under assumptions {(\textbf{H1})}--{(\textbf{H2})}, there exists a unique global solution $(\tilde{X}^{(j)},(\tilde{X}^{(j)})^{\dagger},(\tilde{X}^{(j)})^{\dagger\dagger})\in \mathcal{Q}_{U}^\beta([0, T], \mathbb{R}^m)$ to the \eqref{4-2-12}. Moreover, we have that
$$\|\tilde{X}^{(j)}\|_{\beta-\operatorname{hld}}\le c$$
holds for $0<\beta<H$. Here, the constant $c>0$ is independent of $U$. 

Due to the compact embedding that $\mathcal{C}^{\beta-\operatorname{hld}}([0,T],\mathbb{R}^{m}) \subset \subset \mathcal{C}^{\beta'-\operatorname{hld}}([0,T],\mathbb{R}^{m})$ for any small parameter $\beta'<\beta$, we have that the family $\{\tilde{X}^{(j)}\}_{ j\ge 1}$ is pre-compact in $ \mathcal{C}^{\beta'-\operatorname{hld}}([0,T],\mathbb{R}^{m})$. Let $\check X$ be any limit point. Then, there exists a subsequence of $\{\check{X}^{(j)}\}_{ j\ge 1}$ (denoted by the same symbol) converges to $\tilde X$ in $ \mathcal{C}^{\beta'-\operatorname{hld}}([0,T],\mathbb{R}^{m})$. In the following, 
we will prove that the limit point $\tilde X$ satisfies the RDE as follows,
\begin{eqnarray}\label{4-2-13}
d\tilde{X}_t = \bar f(\tilde{X}_t)dt + \sigma( \tilde{X}_t)dU_t.
\end{eqnarray}

According to \remref{rem3}, we emphasize that $\{\check{X}^{(j)}\}_{ j\ge 1}$ solves the following ODE:
\begin{eqnarray}\label{4-2-14}
d\tilde{X}^{(j)}_t = \bar{f}(\tilde{X}^{(j)}_t)dt + \sigma( \tilde{X}^{(j)}_t)du^{(j)}_t.
\end{eqnarray} 
where $\|(u^{(j)},v^{(j)})\|_{q-\operatorname{var}}<\infty$ with $(H+1/2)^{-1}<q<2$ for all $j\ge 1$.
Due to the Young integral theory, it is not too difficult to verify that for all $(u, v) \in S_{N}$, there exists a unique solution $\{\tilde{X}^{(j)}\}_{ j\ge 1} \in \mathcal{C}^{p-\operatorname{var}}\left([0,T],\mathbb{R}^m\right)$ to \eqref{4-2-14} in the Young sense. In fact, $\{\tilde{X}^{(j)}\}_{ j\ge 1}$ is independent of $\{v^{(j)}\}_{ j\ge 1}$. Moreover, we have
$$\|\check{X}^{(j)}\|_{q-\operatorname{var}}\le c,$$
where the constant $c>0$ is independent of $(u^{(j)}, v^{(j)})$. 
Note that the Young integral $u^{(j)}\mapsto\int_{0}^{\cdot}\sigma(\tilde{X}_s^{(j)})du^{(j)}_s$ is a linear continuous map from $\mathcal{H}^{H,d_1}$ to $\mathcal{C}^{q-\operatorname{var}}\left([0,T],\mathbb{R}^m\right)$. 
Furthermore, we have 
$\|\tilde{X}^{(j)}\|_{p-\operatorname{var}}\le \|\tilde{X}^{(j)}\|_{q-\operatorname{var}}\le c$
for $p>q$. 

Let us show that the limit point $\tilde{X}$ satisfies the skeleton equation. By the direct computation, we derive
\begin{eqnarray}\label{4-2-15}
\big|\tilde{X}^{(j)}_t-\tilde{X}_t\big|&\le& \bigg|\int_{0}^{t} \big[\bar{f}(\tilde{X}^{(j)}_s)-\bar {f}(\tilde{X}_s)\big]ds\bigg|+\bigg|\int_{0}^{t} \big[\sigma( \tilde{X}^{(j)}_s)-\sigma( \tilde{X}_s)\big]d u^{(j)}_s\bigg|\cr
&&+\bigg|\int_{0}^{t} \sigma( \tilde{X}_s)\big[du^{(j)}_s-du_s\big]\bigg|\cr
&=:& K_1+K_2+K_3.
\end{eqnarray}
For the first term $K_1$, by using the assumption {(\textbf{H2})}, we have
\begin{eqnarray}\label{4-2-16}
K_1\le L\int_{0}^{t} |\tilde{X}^{(j)}_s-\tilde{X}_s|ds\le C_7\sup_{0\le s\le t}|\tilde{X}^{(j)}_s-\tilde{X}_s|.
\end{eqnarray}
Since $\{\tilde X^{(j)}\}_{ j\ge 1}$ converges to $\tilde X$ in the uniform norm, it is an immediate consequence that
$K_1\to 0$ as $j \to \infty$.

After that, we estimate $J_2$ as following:
\begin{eqnarray}\label{4-2-17}
K_2\le C_8\|u^{(j)}\|_{q-\operatorname{var}}\|\sigma( \tilde{X}^{(j)})-\sigma( \tilde{X}^{(j)})\|_{p-\operatorname{var}}.
\end{eqnarray}
Here, for any $\beta<H$, we could choose suitable $\beta', p,q$. For example, when $\beta=H-\frac{\kappa}{2}$ for any $\varepsilon<H$, we could set $\beta'=H-\kappa$, $\frac{1}{p}=H-2\kappa$ and $\frac{1}{q}=H+\frac{1}{2}-\kappa$ such that $\frac{1}{p}+\frac{1}{q}=1+2\left(H-\frac{1}{4}\right)-3 \kappa$. Then, since we state that $\{\tilde{X}^{(j)}\}_{ j\ge 1}$ converges to $\tilde X$ in $ \mathcal{C}^{\beta'-\operatorname{hld}}([0,T],\mathbb{R}^{m})$, so we could obtain that $K_2\to 0$ as $j \to \infty$.

Next, we estimate $K_3$. To this end, set
\[
B(u^{(j)},\tilde{X}):=\int_0^t \sigma(\tilde{X}_s)\,du^{(j)}_s,
\]
which is a continuous bilinear map from $\mathcal{H}^{H,d_1}\times \mathcal{C}^{p-\mathrm{var}}([0,T],\mathbb{R}^m)$ to $\mathbb{R}$. By the Riesz representation theorem, there exists a unique element in $\mathcal{H}^{H,d_1}$, denoted by $B(\cdot,\tilde{X})$, such that
\[
B(u^{(j)},\tilde{X})=\langle B(\cdot,\tilde{X}), u^{(j)}\rangle_{\mathcal{H}^{H,d_1}}
\quad\text{for all } u^{(j)}\in \mathcal{H}^{H,d_1}.
\]
Note that $B(\cdot,\tilde{X})\in (\mathcal{H}^{H,d_1})^*\cong \mathcal{H}^{H,d_1}$. Consequently, we have
\begin{eqnarray}\label{4-2-18}
K_3= |B(u^{(j)},\tilde{X})-B(u,\tilde{X})|
= |\langle B(\cdot,\tilde{X}),u^{(j)}\rangle_{\mathcal{H}^{H,d_1}}-\langle B(\cdot,\tilde{X}),u\rangle_{\mathcal{H}^{H,d_1}}|.
\end{eqnarray}
Since $(u^{(j)}, v^{(j)})\rightarrow(u, v)$ as $j\rightarrow\infty$ with the weak topology in $\mathcal{H}$, we prove that $K_3$ converges to 0 as $j\rightarrow\infty$. 
By combining the above analysis \eqref{4-2-14}-\eqref{4-2-18}, it is clear that the limit point $\tilde X$ satisfies the ODE \eqref{4-2-14}. 

Subsequently, we will show that the \textbf{Condition (ii)} holds.

\textbf{Verification of \textbf{Condition (ii)}}. Firstly, we construct an auxiliary process:
\begin{eqnarray}\label{5-8}
d\hat {X}^{\varepsilon,\delta}_t = \bar f(\hat {X}^{\varepsilon,\delta}_t)dt +\sigma(\hat {X}^{\varepsilon,\delta}_t)d[\mathcal{T}_t^{u^{\varepsilon, \delta}}(B^{\varepsilon,H}) ]
\end{eqnarray}
with initial value $\hat {X}^{\varepsilon,\delta}_0=X_0$. Similar to esitmates in \lemref{lem1}, it is direct to see that $\mathbb{E}[\|\hat {X}^{\varepsilon,\delta}\|^2_{\beta-\operatorname{hld}}]\le C_{\alpha, \beta, N}$
for some constant $C_{\alpha, \beta, N}>0$.

Then, according to some direct computation, it deduces that
\begin{eqnarray}\label{5-9}
&&\tilde {X}^{\varepsilon,\delta}_t-\hat {X}^{\varepsilon,\delta}_t\cr
&=&\int_{0}^{t}[f(\tilde {X}^{\varepsilon, \delta}_s, \tilde {Y}^{\varepsilon, \delta}_s)-f(\tilde {X}^{\varepsilon, \delta}_{s(\Delta)}, \tilde {Y}^{\varepsilon, \delta}_s)]dt +\int_{0}^{t}[f(\tilde {X}^{\varepsilon, \delta}_{s(\Delta)}, \tilde {Y}^{\varepsilon, \delta}_s)-f(\tilde {X}^{\varepsilon, \delta}_{s(\Delta)}, \hat {Y}^{\varepsilon, \delta}_s)]ds \cr
&&+\int_{0}^{t}[f(\tilde {X}^{\varepsilon, \delta}_{s(\Delta)}, \hat {Y}^{\varepsilon, \delta}_s)-\bar f(\tilde {X}^{\varepsilon, \delta}_{s(\Delta)})]ds +\int_{0}^{t}[\bar f(\tilde {X}^{\varepsilon, \delta}_{s(\Delta)})-\bar f(\tilde {X}^{\varepsilon, \delta}_s)]ds \cr
&& +\int_{0}^{t}[\bar f(\tilde {X}^{\varepsilon,\delta}_{s})-\bar f(\hat {X}^{\varepsilon,\delta}_{s})]ds+\int_{0}^{t}[\sigma(\tilde {X}^{\varepsilon, \delta}_s)-\sigma(\hat {X}^{\varepsilon,\delta}_s)]d[\mathcal{T}_s^{u^{\varepsilon, \delta}}(B^{\varepsilon,H}) ]\cr
&=:&M_t+\int_{0}^{t}[\bar f(\tilde {X}^{\varepsilon,\delta}_{s})-\bar f(\hat {X}^{\varepsilon,\delta}_{s})]ds+\int_{0}^{t}[\sigma(\tilde {X}^{\varepsilon, \delta}_s)-\sigma(\hat {X}^{\varepsilon,\delta}_s)]d[\mathcal{T}_s^{u^{\varepsilon, \delta}}(B^{\varepsilon,H}) ].
\end{eqnarray}
Then, by \propref{psi2}, that is
\begin{eqnarray}\label{5-31}
\|\tilde {X}^{\varepsilon, \delta}-\hat {X}^{\varepsilon,\delta}\|_{\beta-\operatorname{hld}} \le c \exp \big[c\left(K^{\prime}+1\right)^\nu\big(\vertiii{\mathcal{T}^{u^{\varepsilon, \delta}}(B^{\varepsilon,H})}_{\alpha-\operatorname{hld}}+1\big)^\nu\big]\|M\|_{3\beta-\operatorname{hld}}.
\end{eqnarray}
Here, $K^{\prime}=\max \{\|\sigma\|_{C_b^4},\|f\|_{\infty}, L\}$. The subsequent proof is to give estimates for $M$. 
To this end, we divide $M$ into four terms as below:
\begin{eqnarray}\label{5-3-10}
M_1&:=&\int_{0}^{t}[f(\tilde {X}^{\varepsilon, \delta}_s, \tilde {Y}^{\varepsilon, \delta}_s)-f(\tilde {X}^{\varepsilon, \delta}_{s(\Delta)}, \tilde {Y}^{\varepsilon, \delta}_s)]dt, \cr
M_2&:=& \int_{0}^{t}[f(\tilde {X}^{\varepsilon, \delta}_{s(\Delta)}, \tilde {Y}^{\varepsilon, \delta}_s)-f(\tilde {X}^{\varepsilon, \delta}_{s(\Delta)}, \hat {Y}^{\varepsilon, \delta}_s)]ds ,\cr
M_3&:=&\int_{0}^{t}[f(\tilde {X}^{\varepsilon, \delta}_{s(\Delta)}, \hat {Y}^{\varepsilon, \delta}_s)-\bar f(\tilde {X}^{\varepsilon, \delta}_{s(\Delta)})]ds , \cr
M_4&:=&\int_{0}^{t}[\bar f(\tilde {X}^{\varepsilon, \delta}_{s(\Delta)})-\bar f(\tilde {X}^{\varepsilon, \delta}_s)]ds.
\end{eqnarray}
Firstly, we estimate $M_1$ with the H\"older inequality, (\textbf{A2}) and \lemref{lem1}, 
\begin{eqnarray}\label{5-3-11}
\mathbb{E}[\sup_{0 \le t \le T}|M_1|^2] =\mathbb{E}\bigg[\sup_{0 \le t \le T}\big|\int_{0}^{t}[f(\tilde {X}^{\varepsilon, \delta}_s, \tilde {Y}^{\varepsilon, \delta}_s)-f(\tilde {X}^{\varepsilon, \delta}_{s(\Delta)}, \tilde {Y}^{\varepsilon, \delta}_s)]ds\big|^2\bigg]
\le LT^2 \mathbb{E}[\|\tilde {X}^{\varepsilon, \delta}\|^2_{\beta-\operatorname{hld}}]\Delta^{2\beta}.
\end{eqnarray}
For the second term $M_2$, with aid of the H\"older inequality and \lemref{lem3}, we get 
\begin{eqnarray}\label{5-3-12}
\mathbb{E}[\sup_{0 \le t \le T}|M_2|^2] &=&\mathbb{E}\bigg[\sup_{0 \le t \le T}\big|\int_{0}^{t}[f(\tilde {X}^{\varepsilon, \delta}_{s(\Delta)}, \tilde {Y}^{\varepsilon, \delta}_s)-f(\tilde {X}^{\varepsilon, \delta}_{s(\Delta)}, \hat {Y}^{\varepsilon, \delta}_s)]ds\big|^2\bigg] \cr
&\le&TL \int_{0}^{T}\mathbb{E}\big[\big|\tilde Y_{s}^{\varepsilon,\delta}-\hat Y_{s}^{\varepsilon,\delta}\big|^2\big]ds \le C_{15}T^2(\frac{\sqrt{\delta}}{{\sqrt{\varepsilon}}}+\Delta^{2\beta}).
\end{eqnarray}
For the term $M_3$, we set $\tilde M_{s,t}=\int_{s}^{t}[f(\tilde {X}^{\varepsilon, \delta}_{r(\Delta)}, \hat {Y}^{\varepsilon, \delta}_r)-\bar f(\tilde {X}^{\varepsilon, \delta}_{r(\Delta)})]dr$. Set $1/2<\eta<1$.
When $0<t-s<2\Delta$, it is straightforward to see that
\begin{eqnarray}\label{5-2-18}
|\tilde M_{s,t}|&\le&L {(2\Delta)^{1-\eta}}{(t-s)^\eta}
\end{eqnarray}
When $t-s>2\Delta$, by applying the Schwarz inequality, it deduces that 
\begin{eqnarray}\label{5-2-19}
\frac{\left|\tilde M_{s, t}\right|^2}{(t-s)^{2\eta}} & \le&\frac{\big|\tilde M_{s,(\lfloor s / \Delta\rfloor+1) \Delta}+\sum_{k=\lfloor s / \Delta\rfloor+1}^{\lfloor t / \Delta\rfloor-1} \tilde M_{k \Delta,(k+1) \Delta}+\tilde M_{\lfloor t / \Delta\rfloor \Delta, t}\big|^2}{(t-s)^{2\eta}} \cr
& \le& C \Delta^{2-2\eta}+\frac{2C(t-s)^{1-2\eta}}{\Delta} \sum_{k=0}^{\lfloor T / \Delta\rfloor-1}\left|\tilde M_{k \Delta,(k+1) \Delta}\right|^2.
\end{eqnarray}
Then, by \eqref{5-2-18} and \eqref{5-2-19}, it deduces that
\begin{eqnarray}\label{5-2-20}
\mathbb{E}[\|M_3\|_{\beta-\operatorname{hld}}^2]&=&\mathbb{E}\bigg[\bigg\|\int_{0}^{\cdot}[f(\tilde {X}^{\varepsilon, \delta}_{s(\Delta)}, \hat {Y}^{\varepsilon, \delta}_s)-\bar f(\tilde {X}^{\varepsilon, \delta}_{s(\Delta)})]ds\bigg\|_{\beta-\operatorname{hld}}^2\bigg]\cr
&\le &\frac{CT}{\Delta^{(1+2\eta)}} \max _{0 \le k \le\left\lfloor\frac{T}{\Delta}\right\rfloor-1} \mathbb{E}\big[\big|\int_{k \Delta}^{(k+1) \Delta}\big(f( \tilde X_{k \Delta}^{\varepsilon, \delta},\hat {Y}^{\varepsilon, \delta}_s)-\bar{f}( \tilde X_{k \Delta}^{\varepsilon, \delta})\big) d s\big|^2\big]\cr
&&+C\Delta^{2(1-\eta)}.
\end{eqnarray}
With some direct but cumbersome computation, we arrive at
\begin{eqnarray}\label{5-2-13}
&&\max _{0 \le k \le\left\lfloor\frac{T}{\Delta}\right\rfloor-1} \mathbb{E}\big[\big|\int_{k \Delta}^{(k+1) \Delta}\big(f( \tilde X_{k \Delta}^{\varepsilon, \delta},\hat {Y}^{\varepsilon, \delta}_s)-\bar{f}( \tilde X_{k \Delta}^{\varepsilon, \delta})\big) d s\big|^2\big]\cr
&\le & C \delta^2 \max _{0 \le k \le\left\lfloor\frac{T}{\Delta}\right\rfloor-1} \int_0^{\frac{\Delta}{\delta}} \int_r^{\frac{\Delta}{\delta}} \mathbb{E}\big[\big\langle f( \tilde X_{k \Delta}^{\varepsilon,\delta}, \hat{Y}_{s \varepsilon+k \Delta}^{\varepsilon,\delta})-\bar{f}( \tilde X_{k \Delta}^{\varepsilon,\delta}),\big.\big.\cr
&&\qquad\qquad\qquad\qquad\qquad\big.\big.f( \tilde X_{k \Delta}^{\varepsilon,\delta}, \hat{Y}_{r \varepsilon+k \Delta}^{\varepsilon,\delta})-\bar{f}( \tilde X_{k \Delta}^{\varepsilon,\delta})\big\rangle\big] d s d r\cr
& \le& C \delta^2 \max _{0 \le k \le\left\lfloor\frac{T}{\Delta}\right\rfloor-1} \int_0^{\frac{\Delta}{\delta}} \int_r^{\frac{\Delta}{\delta}} e^{-\frac{\beta_1}{2}(s-r)} d s d r \cr
& \le &C \delta^2\left(\frac{2}{\beta_1} \frac{\Delta}{\delta}-\frac{4}{\beta_1^2}+e^{\frac{-\beta_1}{2} \frac{\Delta}{\delta}}\right)\cr
& \le &C \delta\Delta.
\end{eqnarray}
Here, we exploit the exponential ergodicity of $\hat Y^{\varepsilon,\delta}$, that is 
\begin{eqnarray}\label{5-2-56}
&&\mathbb{E}\big[\big\langle f( \check X_{k \Delta}^{\varepsilon, \delta}, \hat{Y}_{s \varepsilon+k \Delta}^{\varepsilon, \delta})-\bar{f}( \check X_{k \Delta}^{\varepsilon, \delta}),f( \tilde X_{k \Delta}^{\varepsilon, \delta}, \hat{Y}_{r \varepsilon+k \Delta}^{\varepsilon, \delta})-\bar{f}( \check X_{k \Delta}^{\varepsilon, \delta})\big\rangle\big]\cr&\le& C(1+\mathbb{E}[| \check X_{k \Delta}^{\varepsilon, \delta}|^2]+\mathbb{E}[| \hat Y_{k \Delta}^{\varepsilon, \delta}|^2]) e^{-\frac{\beta_1}{2}(s-r)}\cr
&\le& C e^{-\frac{\beta_1}{2}(s-r)},
\end{eqnarray}
where $\beta_1$ is in (\textbf{A4}). For the first inequality, we refer to \cite[Appendix B]{2023Pei} for instance. The final inequality comes from \lemref{lem1} and \lemref{lem3}. 
So according to the estimates \eqref{5-2-20}--\eqref{5-2-56}, we have
\begin{eqnarray}\label{5-3-13}
\mathbb{E}[\|M_3\|_{\eta-\operatorname{hld}}^2]\le C_{16}\Delta^{2(1-\eta)}+\frac{C_{16}T\delta}{\Delta^{2\eta}}.
\end{eqnarray}
For the forth term $M_4$, recall that $\bar{f}$ is Lipschitz continuous and globally bounded, we obtain
\begin{eqnarray}\label{5-3-14}
\mathbb{E}[\sup_{0 \le t \le T}|M_4|^2] =\mathbb{E}\bigg[\sup_{0 \le t \le T}\bigg|\int_{0}^{t}[\bar f(\tilde {X}^{\varepsilon, \delta}_{s(\Delta)})-\bar f(\tilde {X}^{\varepsilon, \delta}_s)]ds\bigg|^2\bigg]
\le LT^2 \mathbb{E}[\|\tilde {X}^{\varepsilon, \delta}\|^2_{\beta-\operatorname{hld}}]\Delta^{2\beta}.
\end{eqnarray}
We obtain that $M\in \mathcal{C}^{1-\operatorname{hld}}([0,T],\mathbb{R}^m)$ and
\begin{eqnarray}\label{5-32}
\mathbb{E}\left[\|M\|_{3 \beta}^2\right] \le C_{17}\big(\Delta^{2 \beta}+\Delta^{2(1-3 \beta)}+\Delta^{-6 \beta} \delta+\frac{\sqrt{\delta}}{\sqrt{\varepsilon}}\big) .
\end{eqnarray}
By choosing some suitable $\Delta>0$ such that $\mathbb{E}[\|M\|_{3\beta-\operatorname{hld}}^2] \to 0$ as $\varepsilon \to 0$ such as set $\Delta:=\delta^{1 /(4 \beta)} \log \delta^{-1}$.
Hence, it deduces that $\|\tilde {X}^{\varepsilon, \delta}-\hat {X}^{\varepsilon,\delta}\|^2_{\beta-\operatorname{hld}}$ converges to $0$ in probability as $\varepsilon \to 0$.

Then, it remains to show the convergence of the auxiliary process {$\hat X^{\varepsilon,\delta}$ at \eqref{5-9}. 
By the condition that $(u^{\varepsilon,\delta}, v^{\varepsilon,\delta})\in \mathcal{A}^{N}_b$, we have that $\mathcal{T}^{u^{\varepsilon, \delta}}(B^{\varepsilon,H}) : \mathcal{H}^{H,d_1}\mapsto \Omega_\alpha(\mathbb{R}^{d_1})$ is a Lipschitz continuous mapping. Next, by \propref{prop2.5}, it has that $\tilde {X}^{\varepsilon,\delta}$ is a continuous solution map with respect to RP $ \Sigma^{\varepsilon}$. Thanks to the condition that $(u^{\varepsilon,\delta}, v^{\varepsilon,\delta})$ weakly converges to $(u, v)$ as $\varepsilon \to 0$ and continuous mapping theorem, it directly deduces that $\tilde X^{\varepsilon,\delta}$ converges in distribution to $\tilde X$ as $\varepsilon \to 0$.
With aid of the Portmanteau theorem \cite[Theorem 13.16]{2020Klenke}, it deduces that for any bounded Lipschitz functions $F:\mathcal{C}^{\beta-\operatorname{hld}}\left([0,T], \mathbb{R}^m\right) \to\mathbb{R}$, that 
\begin{eqnarray*}\label{5-34}
|\mathbb{E}[F(\tilde X^{\varepsilon,\delta})]-\mathbb{E}[F(\tilde X)]|&\le& |\mathbb{E}[F(\tilde X^{\varepsilon,\delta})]-\mathbb{E}[F(\hat X^{\varepsilon,\delta})]|+|\mathbb{E}[F(\hat X^{\varepsilon,\delta})]-\mathbb{E}[F(\tilde X)]|\cr
&\le& \|F\|_\textrm{Lip}\mathbb{E}[\|\tilde X^{\varepsilon,\delta}-\hat X^{\varepsilon,\delta}\|_{\beta\textrm{-hld}}^2]^{\frac{1}{2}}+|\mathbb{E}[F(\hat X^{\varepsilon,\delta})]-\mathbb{E}[F(\tilde X)]|\to 0
\end{eqnarray*} 
as $\varepsilon \to 0$. Here, $\|F\|_\textrm{Lip}$ stands the Lipschitz constant of $F$. 
The proof is completed.\qed

\subsection{Proof of the MDP}
Denote the deviation between $X^{\varepsilon,\delta}$ and $\bar X$ by $Z^{\varepsilon,\delta}$, so it is defined as below:
\begin{eqnarray}\label{deviation-5}
Z^{\varepsilon,\delta}:=\frac{X^{\varepsilon,\delta}-\bar X}{\sqrt{\varepsilon}h(\varepsilon)}.
\end{eqnarray}
Here, $h:(0,1]\to (0,\infty)$ is continuous function satisfying that $\lim_{\varepsilon \to 0} h(\varepsilon)= \infty$ and $\lim_{\varepsilon \to 0}\sqrt{\varepsilon}h(\varepsilon)=0$ for all $\varepsilon\in (0,1]$. We also assume that $\lim_{\varepsilon \to 0}\sqrt{\varepsilon}h(\varepsilon)=0$. Directly, we have
\begin{eqnarray}\label{5-3-1}
d{Z}^{\varepsilon,\delta}_t &=& \frac{1}{\sqrt{\varepsilon}h(\varepsilon)}({f}({X}^{\varepsilon,\delta}_t)-{f}(\bar {X}_t))dt +\frac{1}{h(\varepsilon)} G( X^{\varepsilon,\delta}_t)db^H_t,\quad {Z}^\varepsilon_0=0.
\end{eqnarray}

It is equivalent between showing the MDP for $X^{\varepsilon,\delta}$ and the LDP for $Z^{\varepsilon,\delta}$. To this end, we intend to apply the variational representation \propref{prop2-1}, then the controlled slow-fast RDE associated to the original slow-fast component $({X}^{\varepsilon, \delta}, {Y}^{\varepsilon, \delta})$ is given firstly. 
\begin{eqnarray}\label{5-3-7}
\left
\{
\begin{array}{ll}
d\check {X}^{\varepsilon, \delta}_t =& f(\check {X}^{\varepsilon, \delta}_t, \check {Y}^{\varepsilon, \delta}_t)dt + \sigma(\check {X}^{\varepsilon, \delta}_t)d\mathcal{T}_t^{\sqrt{\varepsilon}h(\varepsilon)u^{\varepsilon,\delta}}(B^{\sqrt{\varepsilon},H})\\
d\check {Y}^{\varepsilon, \delta}_t =& \frac{1}{\delta}F( \check {X}^{\varepsilon, \delta}_t, \check {Y}^{\varepsilon, \delta}_t)dt +\frac{h(\varepsilon)}{\sqrt{\delta}}G(\check {X}^{\varepsilon, \delta}_t, \check {Y}^{\varepsilon, \delta}_t)dv^{\varepsilon, \delta}_t + \frac{1}{\sqrt{\delta}}G(\check {X}^{\varepsilon, \delta}_t, \check {Y}^{\varepsilon, \delta}_t)dw_{t}.
\end{array}
\right.
\end{eqnarray}
Here, $\mathcal{T}^{\sqrt{\varepsilon}h(\varepsilon)u^{\varepsilon,\delta}}(B^{\sqrt{\varepsilon},H})=(\mathcal{T}^{\sqrt{\varepsilon}h(\varepsilon)u^{\varepsilon,\delta},1}(B^{\sqrt{\varepsilon},H}),\mathcal{T}^{\sqrt{\varepsilon}h(\varepsilon)u^{\varepsilon,\delta},2}(B^{\sqrt{\varepsilon},H}), \mathcal{T}^{\sqrt{\varepsilon}h(\varepsilon)u^{\varepsilon,\delta},3}(B^{\sqrt{\varepsilon},H}))\in \Omega_{\alpha}(\mathbb{R}^{d})$ is GRP lifted from the translation of FBM $\sqrt{\varepsilon}b^H$ in the direction of C-M elements $\sqrt{\varepsilon}h(\varepsilon)u^{\varepsilon,\delta}$, which is according to \propref{FBM+CM-RP}. Now, $(u^{\varepsilon, \delta}, v^{\varepsilon, \delta})\in \mathcal{A}_{b}$ is called a pair of control. 

Then, the controlled deviation component is defined as below:
\begin{eqnarray}\label{5-3-8}
\check{Z}^{\varepsilon,\delta}_t :=\frac{\check {X}^{\varepsilon, \delta}_t-\bar X_t}{\sqrt{\varepsilon}h(\varepsilon)},\quad\check{Z}^{\varepsilon,\delta}_0=0.
\end{eqnarray}

Before proving \thmref{mdpthm}, some prior estimates are given as below.
\begin{lem}\label{lem5-1}
Assume (\textbf{A1})--(\textbf{A4}) and let $\nu\ge 1$ and $N\in\mathbb{N}$. Then, for all $\varepsilon,\delta\in(0,1]$, we have
\begin{equation}\label{5-6}
\mathbb{E}\big[\|\check X^{\varepsilon,\delta}\|^\nu_{\beta-\operatorname{hld}}\big] \le C_{\nu,N}, \quad \int_{0}^{T} \mathbb{E}\big[|\check Y_t^{\varepsilon,\delta}|^2\big] dt\le C_{\alpha, \beta, N}.
\end{equation}
Here, $C_{\nu,N}$ and $C_{\alpha, \beta, N}>0$ are constants.
\end{lem}
\para{Proof}. Same proof as that in \lemref{lem1} still works. \qed

\para{Proof of \thmref{mdpthm}}.
It is equivalent between showing the MDP for $X^{\varepsilon,\delta}$ and the LDP for $Z^{\varepsilon,\delta}$. So we prove that \textbf{Condition (i)} and \textbf{Condition (ii)} in \propref{weak_convergence} hold. Similar to the proof in \thmref{thm}, \textbf{Condition (i)} is easy to verify. 

Everything in this step is deterministic.
Let $(u^{(j)}, v^{(j)})$ and $(u, v)$ belong to $S_N$ such that $(u^{(j)}, v^{(j)})\rightarrow(u, v)$ as $j\rightarrow\infty$ in the weak topology of $\mathcal{H}$. 
In this step, we will show 
\begin{eqnarray} \label{3.26}
\mathcal{G}^{0}( u^{(j)} ,v^{(j)} )\rightarrow\mathcal{G}^{0}(u,v) \label{step2}
\end{eqnarray}
in $C^{\alpha}([0,T],\mathbb{R}^m)$ as $j \to \infty$.

Due to the {continuous embedding} 
$\mathcal{H}^{H,d_1}\subset C^{1-\alpha}([0,T], \mathbb{R}^{d_1})$, 
we obtain that $\{u^{(j)}\}_{j\ge 1} \subset C^{1-\alpha}([0,T], \mathbb{R}^{d_1})$. Let $\{\check {Z}^{(j)}_t\}_{j\ge 1}$ be a family of solutions to the below skeleton equation,
\begin{eqnarray*} \label{3.25}
d\check {Z}^{(j)}_t = D\bar{f}( \bar{X}_t)\check {Z}^{(j)}dt + \sigma( \bar{X}_t)du^{(j)}_t,\quad \tilde {z}^{(j)}_0=0.
\end{eqnarray*} 
Due to Assumption (\textbf{A1'}), (\textbf{A4'}), the fact that $\bar X$ is a deterministic differentiable path, and the result \cite[Proposition 3.6]{2020Budhiraja}, there exists a unique solution $\tilde {Z}^{(j)} \in W_0^{\alpha, \infty}([0,T], \mathbb{R}^{m})$ to the skeleton equation \eqref{3.25}. Furthermore, there exists a positive constant $C:=C_N$ such that
\begin{equation*} \label{3.28}
\sup_{ j\geq 1}\|\tilde Z^{(j)}\|_{1-\alpha} \leq C. 
\end{equation*} 

Next, we will verify that there exists a subsequence of $\{\check Z^{(j)}\}_{ j\ge 1}$ converges to the limit point $\check Z$ in $C^{\alpha}([0,T], \mathbb{R}^{m})$ which satisfies the following ODEs}:
\begin{eqnarray} \label{3.44}
d\check {Z}_t = D\bar{f}( \bar{X}_t)\check Z_tdt + \sigma( \bar{X}_t)du_t, \quad \check {Z}_0=0.
\end{eqnarray} 
Note that the Young integral $u^{(j)}\mapsto\int_{0}^{\cdot}\sigma(\check{X}_s)du^{(j)}_s$ is a linear continuous map from $\mathcal{H}^{H,d_1}$ to $\mathcal{C}^{q-\operatorname{var}}\left([0,T],\mathbb{R}^m\right)$. 
Furthermore, we have 
$\|\check{Z}^{(j)}\|_{p-\operatorname{var}}\le \|\tilde{X}\|_{q-\operatorname{var}}\le c$
for $p>q$. 

Let us show that the limit point $\check{Z}$ satisfies the skeleton equation. By the direct computation, we derive
\begin{eqnarray}\label{4-2-15}
\big|\check{Z}^{(j)}_t-\check{Z}_t\big|&\le& \bigg|\int_{0}^{t} D\bar{f}(\bar{X}_s)(\check{Z}^{(j)}_s-\check Z_s)ds\bigg|+\bigg|\int_{0}^{t} \sigma( \bar{X}_s)\big[du^{(j)}_s-du_s\big]\bigg|\cr
&=:& H_1+H_2.
\end{eqnarray}
We first estimate the term $H_1$. With aid of the \remref{bar_f}, under assumptions (\textbf{A5})--(\textbf{A7}), it deduces that $\bar f\in C_b^2$. Therefore, we have
\begin{eqnarray}\label{4-2-16}
H_1\le L\int_{0}^{t} |\check{Z}^{(j)}_s-\check{Z}_s|ds\le C_7\sup_{0\le s\le t}|\check{X}^{(j)}_s-\check{X}_s|.
\end{eqnarray}
Since $\{\check Z^{(j)}\}_{ j\ge 1}$ converges to $\check Z$ in the uniform topology, it is an immediate consequence that
$K_1\to 0$ as $j \to \infty$.

Next, it proceeds to estimates $H_2$. To do this, we set $B(u^{(j)},\bar{X}):=\int_{0}^{t} \sigma( \bar{X}_s)du^{(j)}_s$, which is a bilinear continuous map from $\mathcal{H}^{H,d_1}\times \mathcal{C}^{1}\left([0,T],\mathbb{R}^m\right)$ to $\mathbb{R}$ due to the result that $\bar X$ is a deterministic differentiable path. According to the Riesz representation theorem, there exists a unique element in $ \mathcal{H}^{H,d_1}$ (denoted by $B(\cdot,\bar{X})$) such that $B(u^{(j)},\bar{X})=\langle B(\cdot,\bar{X}),u^{(j)}\rangle_{\mathcal{H}^{H,d_1}}$ for all $u^{(j)}\in \mathcal{H}^{H,d_1}$. Note that $B(\cdot,\bar{X})\in (\mathcal{H}^{H,d_1})^{*}\cong \mathcal{H}^{H,d_1}$. Then, we have 
\begin{eqnarray}\label{4-2-18}
H_2= |B(u^{(j)},\bar{X})-B(u,\bar{X})|
= |\langle B(\cdot,\bar{X}),u^{(j)}\rangle_{\mathcal{H}^{H,d_1}}-\langle B(\cdot,\bar{X}),u\rangle_{\mathcal{H}^{H,d_1}}|.
\end{eqnarray}
Since $(u^{(j)}, v^{(j)})\rightarrow(u, v)$ as $j\rightarrow\infty$ with the weak topology in $\mathcal{H}$, we prove that $H_2$ converges to 0 as $j\rightarrow\infty$. 
By combining the above analysis \eqref{4-2-14}-\eqref{4-2-18}, it is clear that the limit point $\check Z$ satisfies the ODE \eqref{4-2-14}.

Subsequently, we will show \textbf{Condition (ii)} holds.

\textbf(Verification of \textbf{Condition (ii)}.)
Assume that $(u^{\varepsilon,\delta}, v^{\varepsilon,\delta})\in \tilde{\mathcal{A}}^{N}_b$ converges in distribution to $(u, v)$ as $\varepsilon \to 0$. We rewrite the controlled deviation component as below:
\[
\check{Z}^{\varepsilon,\delta}=\mathcal{G}^{\varepsilon,\delta}\bigg(\frac{b^H}{h(\varepsilon)}+u^{\varepsilon} , \frac{W}{h(\varepsilon)} +v^{\varepsilon}\bigg).
\]
We aim to show that $\check{Z}^{\varepsilon,\delta}$ converges in distribution to $\check {Z}$ in $ C^{\beta}([0,T],\mathbb{R}^m)$ as $\varepsilon\rightarrow 0$, that is,
\begin{eqnarray}\label{5-3-2}
\mathcal{G}^{\varepsilon,\delta}\bigg(\frac{b^H}{h(\varepsilon)}+u^{\varepsilon} , \frac{w}{h(\varepsilon)} +v^{\varepsilon}\bigg)\rightarrow\mathcal{G}^{0}(u , v) \quad\text{(in distribution)}.
\end{eqnarray} 
Firstly, we define the auxiliary controlled slow process 
\begin{eqnarray*}\label{5-3-3}
d\bar {X}^{\varepsilon, \delta}_t = \bar f(\bar {X}^{\varepsilon, \delta}_t)dt + \sigma(\bar {X}^{\varepsilon, \delta}_t)d\mathcal{T}_t^{\sqrt{\varepsilon}h(\varepsilon)u^{\varepsilon,\delta}}(B^{\sqrt{\varepsilon},H})
\end{eqnarray*}
Following the same approach as in the proof of Lemma \ref{lem1}, for $\nu\ge 1$ and $N\in\mathbb{N}$, we have, for all $\varepsilon,\delta\in(0,1]$,
\begin{equation}\label{6-6}
\mathbb{E}\big[\|\bar X^{\varepsilon,\delta}\|^\nu_{\beta-\operatorname{hld}}\big] \le C_{\nu,N}.
\end{equation}
Here, $C_{\nu,N}$ and $C_{\alpha, \beta, N}>0$ are constants.

Denote $\check Z^{\varepsilon,\delta}=\hat Z^{\varepsilon,\delta}+\bar Z^{\varepsilon,\delta}$, where 
\begin{eqnarray}\label{5-3-4}
\hat {Z}^{\varepsilon,\delta}_t =\frac{\check{X}^{\varepsilon,\delta}_t-\bar {X}^{\varepsilon,\delta}_t}{\sqrt{\varepsilon}h(\varepsilon)},\qquad \bar Z^{\varepsilon,\delta}_t =\frac{\bar{X}^{\varepsilon,\delta}_t-\bar {X}_t}{\sqrt{\varepsilon}h(\varepsilon)}.
\end{eqnarray}
Here, $\hat {Z}^{\varepsilon}_0=0$ and $\bar {Z}^{\varepsilon}_0=0$. To show \eqref{5-3-2}, it suffices to verify two statements \textbf{(a)} and \textbf{(b)} as below.
\begin{itemize}
\item[\textbf{(a)}.] 
For any $\Delta_1>0$,
\begin{eqnarray}\label{5-3-5}
\lim_{\varepsilon \to 0}\mathbb{P}\big(\|\hat Z^{\varepsilon,\delta}\|_{\beta\textrm{-hld}}>\Delta_1\big)=0.
\end{eqnarray}
\item[\textbf{(b)}.]
As $\varepsilon \to 0$,
\begin{eqnarray}\label{5-3-6}
\bar Z^{\varepsilon,\delta} \rightarrow\check Z \quad (\text{in} \,\, C^{\beta}([0,T],\mathbb{R}^m) \,\text{in distribution}).
\end{eqnarray}
\end{itemize}

We firstly show the statement \textbf{(a)} holds.

(\textbf{Verification of Statement \textbf{(a)}}). 
Before verifying the Statement \textbf{(a)}, we first construct the auxiliary process as below:
\begin{eqnarray}\label{5-3-9}
d\hat {Y}^{\varepsilon, \delta}_t &=& \frac{1}{\delta} F( \tilde {X}^{\varepsilon, \delta}_{t(\Delta)}, \hat {Y}^{\varepsilon, \delta}_t)dt + \frac{1}{{\sqrt \delta }}G(\tilde {X}^{\varepsilon, \delta}_{t(\Delta)}, \hat {Y}^{\varepsilon, \delta}_t)dw_{t}
\end{eqnarray}
with $\hat {Y}^{\varepsilon, \delta}_0=Y_0$.

Under (\textbf{A1})--(\textbf{A4}) and let $N\in\mathbb{N}$, then, we have 
\begin{eqnarray}\label{6-3-9}
\int_{0}^{T}\mathbb{E}[{| \tilde Y_{t}^{\varepsilon,\delta}-\hat Y_{t}^{\varepsilon,\delta} |^2}]dt \le C_{\alpha, \beta, N}(\delta h^2(\varepsilon)+\Delta^{2\beta}).
\end{eqnarray}
Here, $C_{\alpha, \beta, N}>0$ is a constant which depends only on $N,\alpha,\beta$. The proof for \eqref{6-3-9} is similar to \lemref{lem3}, with minor modifications to adapt the argument to our setting. So precise details of the proof is given in Appendix A.

Then, according to some direct computation, it deduces that
\begin{eqnarray}\label{5-3-10}
&&\check {X}^{\varepsilon,\delta}_t-\bar {X}^{\varepsilon,\delta}_t\cr
&=&\int_{0}^{t}[f(\check {X}^{\varepsilon, \delta}_s, \check {Y}^{\varepsilon, \delta}_s)-f(\check {X}^{\varepsilon, \delta}_{s(\Delta)}, \check {Y}^{\varepsilon, \delta}_s)]dt +\int_{0}^{t}[f(\check {X}^{\varepsilon, \delta}_{s(\Delta)}, \check {Y}^{\varepsilon, \delta}_s)-f(\check {X}^{\varepsilon, \delta}_{s(\Delta)}, \hat {Y}^{\varepsilon, \delta}_s)]ds \cr
&&+\int_{0}^{t}[f(\check {X}^{\varepsilon, \delta}_{s(\Delta)}, \hat {Y}^{\varepsilon, \delta}_s)-\bar f(\check {X}^{\varepsilon, \delta}_{s(\Delta)})]ds +\int_{0}^{t}[\bar f(\check {X}^{\varepsilon, \delta}_{s(\Delta)})-\bar f(\check {X}^{\varepsilon, \delta}_s)]ds \cr
&& +\int_{0}^{t}[\bar f(\check {X}^{\varepsilon,\delta}_{s})-\bar f(\bar {X}^{\varepsilon,\delta}_{s})]ds+\int_{0}^{t}[\sigma(\check {X}^{\varepsilon, \delta}_s)-\sigma(\bar {X}^{\varepsilon,\delta}_s)]d[\mathcal{T}_s^{\sqrt{\varepsilon}h(\varepsilon)u^{\varepsilon,\delta}}(B^{\sqrt{\varepsilon},H}) ]\cr
&=:&Q_t+\int_{0}^{t}[\bar f(\check {X}^{\varepsilon,\delta}_{s})-\bar f(\bar {X}^{\varepsilon,\delta}_{s})]ds+\int_{0}^{t}[\sigma(\check {X}^{\varepsilon, \delta}_s)-\sigma(\bar {X}^{\varepsilon,\delta}_s)]d[\mathcal{T}_s^{\sqrt{\varepsilon}h(\varepsilon)u^{\varepsilon,\delta}}(B^{\sqrt{\varepsilon},H})].
\end{eqnarray}
Then, by the result in \propref{psi2}, that is
\begin{eqnarray}\label{5-3-21}
\|\hat {Z}^{\varepsilon,\delta}\|_{\beta-\operatorname{hld}} &=&\|\frac{\check {X}^{\varepsilon, \delta}-\bar {X}^{\varepsilon,\delta}}{\sqrt{\varepsilon}h(\varepsilon)}\|_{\beta-\operatorname{hld}} \cr
&\le& c \exp \big[c\left(K^{\prime}+1\right)^\nu\big(\vertiii{\mathcal{T}^{\sqrt{\varepsilon}h(\varepsilon)u^{\varepsilon,\delta}}(B^{\sqrt{\varepsilon},H})}_{\alpha-\operatorname{hld}}+1\big)^\nu\big]\|\frac{Q}{\sqrt{\varepsilon}h(\varepsilon)}\|_{3\beta-\operatorname{hld}}.
\end{eqnarray}
Here, $K^{\prime}=\max \{\|\sigma\|_{C_b^4},\|f\|_{\infty}, L\}$. The subsequent proof is to give the estimation for $Q$. To this end, we divide it into four terms as below:
\begin{eqnarray}\label{5-3-10}
Q_1&=&\int_{0}^{t}[f(\check {X}^{\varepsilon, \delta}_s, \check {Y}^{\varepsilon, \delta}_s)-f(\check {X}^{\varepsilon, \delta}_{s(\Delta)}, \check {Y}^{\varepsilon, \delta}_s)]dt, \cr
Q_2&=& \int_{0}^{t}[f(\check {X}^{\varepsilon, \delta}_{s(\Delta)}, \check {Y}^{\varepsilon, \delta}_s)-f(\check {X}^{\varepsilon, \delta}_{s(\Delta)}, \hat {Y}^{\varepsilon, \delta}_s)]ds ,\cr
Q_3&=&\int_{0}^{t}[f(\check {X}^{\varepsilon, \delta}_{s(\Delta)}, \hat {Y}^{\varepsilon, \delta}_s)-\bar f(\check {X}^{\varepsilon, \delta}_{s(\Delta)})]ds , \cr
Q_4&=&\int_{0}^{t}[\bar f(\check {X}^{\varepsilon, \delta}_{s(\Delta)})-\bar f(\check {X}^{\varepsilon, \delta}_s)]ds.
\end{eqnarray}
The estimates of terms $Q_1$, $Q_3$, $Q_4$ are similar to that in \eqref{5-3-11}--\eqref{5-3-14}, so we only show the results but omitting details of the proof for the breviety.
\begin{eqnarray}\label{5-3-11}
\mathbb{E}[\sup_{0 \le t \le T}|Q_1|^2] =\mathbb{E}\bigg[\sup_{0 \le t \le T}\big|\int_{0}^{t}[f(\check {X}^{\varepsilon, \delta}_s, \hat {Y}^{\varepsilon, \delta}_s)-f(\check {X}^{\varepsilon, \delta}_{s(\Delta)}, \hat {Y}^{\varepsilon, \delta}_s)]ds\big|^2\bigg]
\le LT^2 \mathbb{E}[\|\check {X}^{\varepsilon, \delta}\|^2_{\beta-\operatorname{hld}}]\Delta^{2\beta},
\end{eqnarray}
\begin{eqnarray}\label{6-3-13}
\mathbb{E}[\|Q_3\|_{\eta-\operatorname{hld}}^2]\le C_{16}\Delta^{2(1-\eta)}+\frac{C_{16}T\delta}{\Delta^{2\eta}},
\end{eqnarray}
and
\begin{eqnarray}\label{6-3-14}
\mathbb{E}[\sup_{0 \le t \le T}|Q_4|^2] =\mathbb{E}\bigg[\sup_{0 \le t \le T}\bigg|\int_{0}^{t}[\bar f(\check {X}^{\varepsilon, \delta}_{s(\Delta)})-\bar f(\check {X}^{\varepsilon, \delta}_s)]ds\bigg|^2\bigg]
\le LT^2 \mathbb{E}[\|\check {X}^{\varepsilon, \delta}\|^2_{\beta-\operatorname{hld}}]\Delta^{2\beta}.
\end{eqnarray}
Then it remains to estimate $Q_2$ which is different from LDP. With aid of the H\"older inequality and \lemref{lem3}, we get 
\begin{eqnarray}\label{6-3-12}
\mathbb{E}[\sup_{0 \le t \le T}|Q_2|^2] &=&\mathbb{E}\bigg[\sup_{0 \le t \le T}\big|\int_{0}^{t}[f(\check {X}^{\varepsilon, \delta}_{s(\Delta)}, \check {Y}^{\varepsilon, \delta}_s)-f(\check {X}^{\varepsilon, \delta}_{s(\Delta)}, \hat {Y}^{\varepsilon, \delta}_s)]ds\big|^2\bigg] \cr
&\le&TL \int_{0}^{T}\mathbb{E}\big[\big|\check Y_{s}^{\varepsilon,\delta}-\hat Y_{s}^{\varepsilon,\delta}\big|^2\big]ds \le C_{18}T^2(\delta h^2(\varepsilon)+\Delta^{3\beta}).
\end{eqnarray}
Then, we have $Q\in \mathcal{C}^{1-\operatorname{hld}}([0,T],\mathbb{R}^m)$ and
\begin{eqnarray}\label{5-3-12}
\mathbb{E}\left[\|\frac{Q}{\sqrt{\varepsilon}h(\varepsilon)}\|_{3 \beta}^2\right] \le C_{19}\big(\frac{\Delta^{2 \beta}}{{\varepsilon}h^2(\varepsilon)}+\frac{\Delta^{2(1-3 \beta)}}{{\varepsilon}h^2(\varepsilon)}+\frac{\delta}{\Delta^{6\beta}{\varepsilon}h^2(\varepsilon)}+\frac{{\delta}}{\varepsilon}\big) .
\end{eqnarray}
By choosing some suitable $\Delta>0$ such that $\mathbb{E}[\|\frac{Q}{\sqrt{\varepsilon}h(\varepsilon)}\|_{3\beta-\operatorname{hld}}^2] \to 0$ as $\varepsilon \to 0$. Then we could take suitable $\Delta$ such that as $\varepsilon$ goes to zero, 
\begin{eqnarray}\label{5-3-12}
\frac{\Delta^{2 \beta}}{{\varepsilon}h^2(\varepsilon)}\to 0, \quad \frac{\Delta^{2(1-3 \beta)}}{{\varepsilon}h^2(\varepsilon)} \to 0, \quad \frac{\delta}{\Delta^{6\beta}{\varepsilon}h^2(\varepsilon)} \to 0.
\end{eqnarray}
Here, we communicate about the choice of $\Delta$.
For example, when $h(\varepsilon)=\varepsilon^{-\frac{\theta}{2}}$ for $\theta\in(3\beta,1)$, we set $\Delta_1:=\varepsilon^{\gamma+1}h^{2}(\varepsilon)|ln \varepsilon|$ with $\gamma=(\theta-3\beta)/(3\beta)$. Then we take $\Delta:=(\Delta_1)^{a(\beta)}$ and $a(\beta):=\frac{1}{2\beta}\vee \frac{1}{2(1-3\beta)} \vee \frac{1}{6\beta}= \frac{1}{2(1-3\beta)}$. It is easy to see that as $\varepsilon$ tends to $0$, $\Delta\to 0$, furthermore, we have that \eqref{5-3-12} holds. The statement \textbf{(a)} is verified. 

Next it aims to verify the statement \textbf{(b)}.

(\textbf{Verification of Statement \textbf{(b)}}). 
$\bar Z^{\varepsilon,\delta}$ satisfying the following RDE,
\begin{eqnarray*}\label{6-3-3}
d\bar {Z}^{\varepsilon, \delta}_t =\frac{1}{\sqrt{\varepsilon}h(\varepsilon)} [\bar f(\bar {X}^{\varepsilon, \delta}_t)-\bar f(\bar {X}_t)\big]dt + \sigma(\bar {X}^{\varepsilon, \delta}_t)d\mathcal{T}_t^{u^{\varepsilon,\delta}}(B^{\frac{1}{h(\varepsilon)},H})
\end{eqnarray*}
Firstly, it requires to show the estimate as below
\begin{equation}\label{4-3-9}
\sup_{0<\varepsilon\le 1}\mathbb{E}\big[\|\bar Z^{\varepsilon,\delta}\|^\nu_{\beta-\operatorname{hld}}\big] \le C_{\nu,N}
\end{equation}
where $C_{\nu,N}>0$ is a constant depending only on $\nu$ and $N$. Firstly, we will show for all $\varepsilon\in(0,1]$, we have for all $\nu\ge 1$, 
\begin{equation}\label{4-3-10}
\mathbb{E}\big[\|\bar Z^{\varepsilon,\delta}\|^\nu_{\infty}\big] \le C_{\nu,N}.
\end{equation}
Here, $C_{\nu,N}>0$ is a constant depending only on $\nu$ and $N$. To this end, 
\begin{eqnarray}\label{4-3-11}
\bar{Z}^{\varepsilon,\delta}_t &=& \int_0^t\frac{1}{\sqrt{\varepsilon}h(\varepsilon)}(\bar {f}(\bar{X}^{\varepsilon,\delta}_s)-\bar{f}({X}_s))ds + \int_0^t\sigma(\bar {X}^{\varepsilon, \delta}_s)d\mathcal{T}_s^{u^{\varepsilon,\delta}}(B^{\frac{1}{h(\varepsilon)},H})ds\cr
&=&\int_0^t\left( \int_0^1 \nabla \bar f (X_s + \theta \sqrt{\varepsilon}h(\varepsilon) \check Z^{\varepsilon}_s ) \langle \bar Z^{{\varepsilon,\delta}}_s\rangle d\theta \right) ds + \int_0^t\sigma( \bar X^{\varepsilon,\delta}_s) d\mathcal{T}_s^{u^{\varepsilon,\delta}}(B^{\frac{1}{h(\varepsilon)},H})\cr
&=:&L_1+L_2.
\end{eqnarray}
According to Assumption (\textbf{H2'}) and H\"older's inequality, we obtain 
\begin{eqnarray}\label{4-3-12}
L_1& \le& C_9 \int_0^t|\bar{Z}^{\varepsilon,\delta}_s| d s,
\end{eqnarray}
Next, we aim to estimate the second term $L_2$. Since $(\bar {X}^{{\varepsilon,\delta}},(\bar {X}^{{\varepsilon,\delta}})^{\dagger},(\bar {X}^{{\varepsilon,\delta}})^{\dagger\dagger})\in \mathcal{Q}_{\mathcal{T}^{\sqrt{\varepsilon}h(\varepsilon)u^{\varepsilon,\delta}}(B^{\sqrt{\varepsilon},H}), [0,T]}^\beta$, furthermore, let $\nu\ge 1$ and $N\in\mathbb{N}$, and for all $\varepsilon\in(0,1]$, we have
$\mathbb{E}\big[\|\bar X^{{\varepsilon,\delta}}\|^\nu_{\beta-\operatorname{hld}}\big] \le C_{\nu,N}$.
By \remref{integral-CRP}, it is directly to see that for a given RP $\mathcal{T}^{\sqrt{\varepsilon}h(\varepsilon)u^{\varepsilon,\delta}}(B^{\sqrt{\varepsilon},H})\in \Omega_{\alpha}(\mathbb{R}^{d+e})$ and controlled RP $(\bar {X}^{{\varepsilon,\delta}},(\bar {X}^{{\varepsilon,\delta}})^{\dagger},(\bar {X}^{{\varepsilon,\delta}})^{\dagger\dagger})\in \mathcal{Q}_{\mathcal{T}_t^{\sqrt{\varepsilon}h(\varepsilon)u^{\varepsilon,\delta}}(B^{\sqrt{\varepsilon},H}), [0,T]}^\beta$, then we have $(\sigma( \bar X^{\varepsilon,\delta}), \sigma( \bar X^{\varepsilon,\delta})^\dagger,\sigma( \bar X^{\varepsilon,\delta})^{\dagger\dagger})\in \mathcal{Q}_{\mathcal{T}^{\sqrt{\varepsilon}h(\varepsilon)u^{\varepsilon,\delta}}(B^{\sqrt{\varepsilon},H}), [0,T]}^\beta$. With some direct computation, it is not too difficult to see that $\|L_2\|_{\beta-\operatorname{hld}}\le C_N$.
So with H\"older inequality, it is straightforward to see that \eqref{4-3-9} holds.

Then it remians to verify the statement \textbf{(b)}, we constructed a new process as below before
\begin{eqnarray}\label{6-1}
d\bar {Z}^{{\varepsilon}}_t =\nabla \bar f(\bar {X}_t)\bar {Z}^{{\varepsilon}}_tdt + \sigma(\bar {X}_t)d\mathcal{T}_t^{u^{\varepsilon,\delta}}(B^{\frac{1}{h(\varepsilon)},H})
\end{eqnarray}
So after some straightforward computation, we have 
\begin{eqnarray}\label{4-3-14}
\bar Z^{\varepsilon,\delta}_t-\bar {Z}^{\varepsilon}_t&=&\int_0^t\left( \int_0^1 \nabla \bar f (\bar X_s + \theta \sqrt{\varepsilon}h(\varepsilon) \bar Z^{{\varepsilon,\delta}}_s ) \langle \bar Z^{{\varepsilon,\delta}}_s\rangle d\theta -\nabla \bar{f}(\bar {X}_s) \langle\bar {Z}^{{\varepsilon,\delta}}_s\rangle\right) ds\cr
&&+\int_0^t\big[\sigma(\bar X^{\varepsilon,\delta}_s)-\sigma(\bar X_s)d\big]d\mathcal{T}_s^{u^{\varepsilon,\delta}}(B^{\frac{1}{h(\varepsilon)},H})\cr
&=:& J_1+J_2.
\end{eqnarray}
With Assumption (\textbf{H2'}) and estimation \eqref{4-3-9}, we have 
\begin{eqnarray}\label{4-3-15}
J_1&\le&\int_0^t\left( \int_0^1 \nabla \bar f ( \bar X_s + \theta \sqrt{\varepsilon}h(\varepsilon) \bar Z^{\varepsilon,\delta}_s ) \langle \bar Z^{\varepsilon,\delta}_s\rangle d\theta -\nabla \bar{f}(\bar{X}_s) \langle\bar {Z}^{{\varepsilon,\delta}}_s\rangle\right) ds\cr
&&+\int_0^t\left(\nabla \bar{f}(\bar{X}_s) \langle\bar {Z}^{{\varepsilon,\delta}}_s\rangle-\nabla \bar{f}(\bar{X}_s)\bar Z^{\varepsilon}_s\right) ds\cr
&\le&\int_0^t \|\nabla \bar f\|_{Lip} \sqrt{\varepsilon}h(\varepsilon) |\bar Z^{\varepsilon,\delta}_s |^2 ds+\int_0^t \nabla \bar{f}(\bar{X}_s)(\bar Z^{\varepsilon,\delta}_s-\bar Z^{\varepsilon}_s)ds\cr
&\le& C_{10}\sqrt{\varepsilon}h(\varepsilon)t+\int_0^t \nabla \bar{f}(\bar{X}_s)(\bar Z^{\varepsilon,\delta}_s-\bar Z^{\varepsilon}_s)ds.
\end{eqnarray}
We now turn to the estimate of $J_2$. To simplify notation, set $\tilde\sigma:=\sigma(\bar X^{\varepsilon,\delta})-\sigma(\bar X)$. Then, from Remark \ref{integral-CRP}, it follows that
\begin{eqnarray}\label{4-3-16}
&&\big|\int_s^t \tilde \sigma_r d \mathcal{T}_r^{u^{\varepsilon,\delta}}(B^{\frac{1}{h(\varepsilon)},H})-(\tilde \sigma_{s}\mathcal{T}_{s,t}^{u^{\varepsilon,\delta},1}(B^{\frac{1}{h(\varepsilon)},H})+\tilde \sigma_{s}^{\dagger}\mathcal{T}_{s,t}^{u^{\varepsilon,\delta},2}(B^{\frac{1}{h(\varepsilon)},H})+\tilde \sigma_{s}^{\dagger\dagger}\mathcal{T}_{s,t}^{u^{\varepsilon,\delta},3}(B^{\frac{1}{h(\varepsilon)},H}))\big|\cr
&\le& 2^{4\alpha}\zeta(4\alpha) (t-s)^{4\alpha}\big(\|\tilde \sigma^{\sharp}\|_{{2\alpha-\operatorname{hld}},[a, b]}\|\mathcal{T}^{u^{\varepsilon,\delta},1}(B^{\frac{1}{h(\varepsilon)},H})\|_{\alpha \mathrm{-hld}}\big. \cr
&&\big.+\|\tilde \sigma^{\sharp\sharp}\|_{{3\alpha-\operatorname{hld}},[a, b]} \|\mathcal{T}^{u^{\varepsilon,\delta},2}(B^{\frac{1}{h(\varepsilon)},H})\|_{2\alpha \mathrm{-hld}}+\|\tilde \sigma^{\dagger\dagger}\|_{{\alpha-\operatorname{hld}},[a, b]}\|\mathcal{T}^{u^{\varepsilon,\delta},3}(B^{\frac{1}{h(\varepsilon)},H})\|_{3\alpha \mathrm{-hld}}\big) 
\end{eqnarray}
where $\zeta$ is the Riemann zeta function. By using the estimate \eqref{4-3-9}, it deduces that $\|\bar {X}^{\varepsilon,\delta}-{ \bar{X}}\|_{\beta-\operatorname{hld}}\le C \sqrt{\varepsilon}h(\varepsilon)$.
By the \propref{function-stability}, we obtain
\begin{eqnarray}\label{4-3-17}
\|\tilde \sigma, \tilde \sigma^{\dagger}, \tilde \sigma^{\dagger\dagger}\|_{\mathcal{T}^{u^{\varepsilon,\delta}(B^{\frac{1}{h(\varepsilon)},H})}, 0, 3 \alpha} &\le& C_{11}\big(\rho_\alpha(\mathcal{T}^{u^{\varepsilon,\delta}}(B^{\frac{1}{h(\varepsilon)},H}), U^{\varepsilon,\delta})+\|\bar X^{\varepsilon,\delta}, (\bar X^{\varepsilon,\delta})^{\dagger}, (\bar X^{\varepsilon,\delta})^{\dagger\dagger} ; X, 0, 0\|_{\mathcal{T}^{u^{\varepsilon,\delta}(B^{\frac{1}{h(\varepsilon)},H})}, 0, 3 \alpha}\big) \cr
&\le& C_{11}\sqrt{\varepsilon}h(\varepsilon),
\end{eqnarray}
and 
\begin{eqnarray}\label{4-3-18}
\|\tilde \sigma\|_{\beta \mathrm{-hld}} &\le& C_{12}\big(\rho_\alpha(\mathcal{T}^{u^{\varepsilon,\delta}}(B^{\frac{1}{h(\varepsilon)},H}), \theta)+\|\bar X^{\varepsilon,\delta}, (\bar X^{\varepsilon,\delta})^{\dagger}, (\bar X^{\varepsilon,\delta})^{\dagger\dagger} ; X, 0, 0\|_{\mathcal{T}^{u^{\varepsilon,\delta}}(B^{\frac{1}{h(\varepsilon)},H}), 0, 3 \alpha}\big)\cr
& \le& C_{\alpha, \sigma_1, \sigma_2}\sqrt{\varepsilon}h(\varepsilon),
\end{eqnarray}
where $\Theta\in G\Omega_\alpha(\mathcal{H})$ are GRPs lifted by $(u,v)\in \mathcal{H}$ and 
$C_{\alpha, \sigma_1, \sigma_2}=C(\alpha, \sigma_1, \sigma_2)>0$ are constants. Then, by substituting \eqref{4-3-17}--\eqref{4-3-18} into Equation \eqref{4-3-16}, we obtain
$$\|\int_s^t \tilde \sigma_u d \mathcal{T}_u^{u^{\varepsilon,\delta},1}(B^{\frac{1}{h(\varepsilon)},H})\|_{\beta \mathrm{-hld}} \le C_{13}\sqrt{\varepsilon}h(\varepsilon).$$
Then with the Gronwall inequality, we arrive at
\begin{eqnarray}\label{4-3-20}
\|\bar Z^{\varepsilon,\delta}-\bar Z^{\varepsilon}\|_{\infty}&\le&C_{14}e^{\|\nabla f\|_{Lip}T}(\sqrt{\varepsilon}h(\varepsilon)t+\|\int_s^t \tilde \sigma_u d \mathcal{T}_u^{u^{\varepsilon,\delta},3}(B^{\frac{1}{h(\varepsilon)},H})\|_{\beta \mathrm{-hld}}t^\beta).
\end{eqnarray}
By the condition $(u^{\varepsilon}, v^{\varepsilon})\in \mathcal{A}^{N}_b$, the map
\[
\mathcal{T}^{u^{\varepsilon,\delta}}(B^{\frac{1}{h{\varepsilon}},H}): \mathcal{H}\mapsto G\Omega_\alpha(\mathbb{R}^{d+e})
\]
is Lipschitz continuous. Moreover, by Proposition \ref{prop2.5}, $\bar{Z}^{\varepsilon,\delta}$ is a continuous solution map with respect to the GRP $\mathcal{T}^{u^{\varepsilon,\delta}}(B^{\frac{1}{h(\varepsilon)},H})$. Since $(u^{\varepsilon}, v^{\varepsilon})$ converges weakly to $(u, v)$ as $\varepsilon\to 0$, the continuous mapping theorem yields the weak convergence of $\bar Z^\varepsilon$ to $\check Z$ as $\varepsilon\to 0$. Thus Condition (ii) is verified. The proof is then completed by Proposition \ref{weak_convergence}.\qed

The proof is completed. \qed




\begin{rem}\label{MDPE}
Our result could also uncover the MDP for the slow-fast RDE driven by RP $(B^H,W)$ which is lifted from the mixed FBM $(b^H,w)$ with $H\in (\frac{1}{3},\frac{1}{2})$. Precisely, consider the slow-fast RDE as follows, 
\begin{eqnarray}\label{6-1}
\left
\{
\begin{array}{ll}
dX^{\varepsilon, \delta}_t = f(X^{\varepsilon, \delta}_t, Y^{\varepsilon, \delta}_t)dt + \sqrt \varepsilon \sigma( X^{\varepsilon, \delta}_t)dB^H_{t},\\
dY^{\varepsilon, \delta}_t =\frac{1}{\delta} F( X^{\varepsilon, \delta}_t, Y^{\varepsilon, \delta}_t)dt + \frac{1}{{\sqrt \delta}}G( X^{\varepsilon, \delta}_t, Y^{\varepsilon, \delta}_t)dW_{t},\\
(X^{\varepsilon, \delta}_0, Y^{\varepsilon, \delta}_0)=(x, y)\in \mathbb{R}^{{m}}\times \mathbb{R}^{{n}}.
\end{array}
\right.
\end{eqnarray}
Now, because RP $(B^H, W)\in \Omega_{\alpha}(\mathbb{R}^{d+e})$ with $1/3<\alpha<H$ only includes the first and second level paths, the regularity of $\sigma$ and $G$ could be ``weaker". Precisely, we impose the following further condition.
\begin{itemize}
\item[\textbf{A1'}.] $\sigma\in \mathcal{C}_b^3$, $G\in \mathcal{C}^3$.
\end{itemize}
Now the solution to the above system could be ensured whose proof could be found in \cite[Proposition 3.3]{2022Inahama}.

Assume {(\textbf{A1'})}, {(\textbf{A1'})}--{(\textbf{A7})} and $\delta=o(\varepsilon)$. Let $\varepsilon \to 0$, the slow component $X^{\varepsilon,\delta}$ of system (\ref{6-1}) satisfies a MDP on $\mathcal{C}^{\beta-\operatorname{hld}}([0,T],\mathbb{R}^{m})$ with speed $1/h^2(\varepsilon)$ and a good rate function $I: \mathcal{C}^{\beta-\operatorname{hld}}([0,T],\mathbb{R}^m)\rightarrow [0, \infty)$
\begin{eqnarray*}\label{rate-1}
I(\xi) &=& {
\inf\Big\{\frac{1}{2}\|u\|^2_{\mathcal{H}^{H,d}}~:~{ u\in \mathcal{H}^{H,d} 
\quad\text{such that} \quad\xi =\mathcal{G}^{0}(u, 0)}\Big\} 
}
\cr
&=& \inf\Big\{\frac{1}{2}\|(u,v)\|^2_{\mathcal{H}} ~:~ {( u,v)\in \mathcal{H}\quad\text{such that} \quad\xi \text{satisfies \eqref{6-2}}}\Big\},
\qquad 
\xi\in \mathcal{C}^{\beta-\operatorname{hld}}\left([0,T], \mathbb{R}^m\right), 
\end{eqnarray*}
where the skeleton equation satisfies the Young ODE as following:
\begin{eqnarray}\label{6-2}
d\check{Z}_t = \nabla \bar {f}(\bar {X}_t)\check Z_tdt + G(\bar {X}_t)du_t, \quad \check{Z}_0=0.
\end{eqnarray}
where $\bar X$ is the solution to the deterministic equation,
$$d \bar X_t =\bar f(X_t) dt, \quad \bar X_0=X_0.$$
The proof is similar to that of \thmref{mdpthm}.
\end{rem}

\section*{Appendix}
\subsection*{A. Proof of \eqref{6-3-9}}
By It\^o's formula, we have
\begin{eqnarray}\label{6-1}
\mathbb{E}[{| \check Y_{t}^{\varepsilon,\delta}-\hat Y_{t}^{\varepsilon,\delta}|^2}] &=& \frac{2}{\delta }\mathbb{E}\bigg[\int_0^t {\langle\check Y_{s}^{\varepsilon,\delta}-\hat Y_{s}^{\varepsilon,\delta},F( \check {X}_{s}^{\varepsilon,\delta},\check Y_s^{\varepsilon,\delta})-F( \check {X}_{s(\Delta)}^{\varepsilon,\delta},\hat Y_s^{\varepsilon,\delta})\rangle ds}\bigg]\cr
&&+ \frac{1}{\delta }\mathbb{E}\bigg[\int_0^t |{G }( {\tilde {X}_{s}^{\varepsilon,\delta},\tilde {Y}_{s}^{\varepsilon,\delta}} )-{G }( {\tilde {X}_{s(\Delta)}^{\varepsilon,\delta},\hat {Y}_{s}^{\varepsilon,\delta}} )|^2ds\bigg] \cr
&&+ \frac{2h(\varepsilon)}{{\sqrt { \delta } }}\mathbb{E}\bigg[\int_0^t {\langle \tilde Y_{s}^{\varepsilon,\delta}-\hat Y_{s}^{\varepsilon,\delta},\sigma_{2} ( {\tilde {X}_{s}^{\varepsilon,\delta},\tilde {Y}_{s}^{\varepsilon,\delta}} ){\frac{dv_s^{\varepsilon,\delta}}{ds} }\rangle ds} \bigg].
\end{eqnarray}
By differentiating with respect to $t$ for \eqref{6-1}, we find that
\begin{eqnarray}\label{6-2}
&&\frac{d}{dt}\mathbb{E}[{| \tilde Y_{t}^{\varepsilon,\delta}-\hat Y_{t}^{\varepsilon,\delta}|^2}]\cr &=& \frac{2}{\delta }\mathbb{E}\big[ {\langle\tilde Y_{t}^{\varepsilon,\delta}-\hat Y_{t}^{\varepsilon,\delta},F( \tilde {X}_{t}^{\varepsilon,\delta},\tilde Y_t^{\varepsilon,\delta})-F( \tilde {X}_{t(\Delta)}^{\varepsilon,\delta},\hat Y_t^{\varepsilon,\delta})\rangle }\big]+ \frac{1}{\delta }\mathbb{E}\big[ |{G }( {\tilde {X}_{t}^{\varepsilon,\delta},\tilde {Y}_{t}^{\varepsilon,\delta}} )-{G }( {\tilde {X}_{t(\Delta)}^{\varepsilon,\delta},\hat {Y}_{t}^{\varepsilon,\delta}} )|^2\big] \cr
&&+\frac{2h(\varepsilon)}{{\sqrt { \delta } }}\mathbb{E}\big[ {\langle \tilde Y_{t}^{\varepsilon,\delta}-\hat Y_{t}^{\varepsilon,\delta},G ( {\tilde {X}_{t(\Delta)}^{\varepsilon,\delta},\tilde {Y}_{t}^{\varepsilon,\delta}} ){\frac{dv_s^{\varepsilon,\delta}}{dt} }\rangle } \big]\cr
&=& \frac{1}{\delta }\mathbb{E}\big[ {2\langle\tilde Y_{t}^{\varepsilon,\delta}-\hat Y_{t}^{\varepsilon,\delta},F(\tilde {X}_{t}^{\varepsilon,\delta},\tilde Y_t^{\varepsilon,\delta})-F( \tilde {X}_{t}^{\varepsilon,\delta},\hat Y_t^{\varepsilon,\delta})\rangle } +|{G}( {\tilde {X}_{t}^{\varepsilon,\delta},\tilde {Y}_{t}^{\varepsilon,\delta}} )-{G }( {\tilde {X}_{t}^{\varepsilon,\delta},\hat {Y}_{t}^{\varepsilon,\delta}} )|^2\big]\cr
&& + \frac{2}{\delta }\mathbb{E}\big[ {\langle\tilde Y_{t}^{\varepsilon,\delta}-\hat Y_{t}^{\varepsilon,\delta},F( \tilde {X}_{t}^{\varepsilon,\delta},\hat Y_t^{\varepsilon,\delta})-F( \hat {X}_{t(\Delta)}^{\varepsilon,\delta},\hat Y_t^{\varepsilon,\delta})\rangle }\big]\cr
&&+ \frac{2}{\delta }\mathbb{E}\big[ \langle{G}( {\tilde {X}_{t}^{\varepsilon,\delta},\tilde {Y}_{t}^{\varepsilon,\delta}} )-{G }( {\tilde {X}_{t}^{\varepsilon,\delta},\hat {Y}_{t}^{\varepsilon,\delta}} ),{G}( {\tilde {X}_{t}^{\varepsilon,\delta},\hat {Y}_{t}^{\varepsilon,\delta}} )-{G}( {\tilde {X}_{t(\Delta)}^{\varepsilon,\delta},\hat {Y}_{t}^{\varepsilon,\delta}} )\rangle\big] \cr
&&+ \frac{1}{\delta }\mathbb{E}\big[ |{G}( {\tilde {X}_{t}^{\varepsilon,\delta},\hat {Y}_{t}^{\varepsilon,\delta}} )-{G}( {\tilde {X}_{t(\Delta)}^{\varepsilon,\delta},\hat {Y}_{t}^{\varepsilon,\delta}} )|^2\big] + \frac{2h(\varepsilon)}{{\sqrt { \delta } }}\mathbb{E}\big[ {\langle \tilde Y_{t}^{\varepsilon,\delta}-\hat Y_{t}^{\varepsilon,\delta},G ( {\tilde {X}_{t}^{\varepsilon,\delta},\tilde {Y}_{t}^{\varepsilon,\delta}} ){\frac{dv_t^{\varepsilon,\delta}}{dt} }\rangle } \big]\cr
&=:&I_1+I_2+I_3+I_4+I_5.
\end{eqnarray}
For the term $K_1$, by leveraging the assumption (\textbf{A4}), it deduces that 
\begin{eqnarray}\label{6-3}
K_1 \le -\frac{\beta_1}{\delta}\mathbb{E}\big[|\tilde Y_{t}^{\varepsilon,\delta}-\hat Y_{t}^{\varepsilon,\delta}|^2\big].
\end{eqnarray}
Then, we compute the second term $K_2$ by applying (\textbf{A2}) and \lemref{lem1} as follows,
\begin{eqnarray}\label{6-4}
K_2 &\le& \frac{c_1}{\delta}\mathbb{E}\big[|\tilde Y_{t}^{\varepsilon,\delta}-\hat Y_{t}^{\varepsilon,\delta}|\cdot|\tilde X_{t}^{\varepsilon,\delta}-\tilde X_{t(\Delta)}^{\varepsilon,\delta}|\big]\cr
&\le& \frac{\beta_1}{4\delta}\mathbb{E}\big[|\tilde Y_{t}^{\varepsilon,\delta}-\hat Y_{t}^{\varepsilon,\delta}|^2\big]+\frac{c_2}{\delta}\Delta^{2\beta}\mathbb{E}[\| \tilde {X}^{\varepsilon,\delta} \|_{\beta-\operatorname{hld}}^2].
\end{eqnarray}
where $c_1,c_2>0$ is independent of $\varepsilon,\delta$.

For the terms $\sum_{i=3}^4 K_i$, by using (\textbf{A3}), \lemref{lem1} and the definition of H\"older norm, we have
\begin{eqnarray}\label{6-5}
I_3+i_4 &\le& \frac{c}{\delta}\mathbb{E}\big[|\tilde Y_{t}^{\varepsilon,\delta}-\hat Y_{t}^{\varepsilon,\delta}|\cdot|\tilde X_{t}^{\varepsilon,\delta}-\tilde X_{t(\Delta)}^{\varepsilon,\delta}|+|\tilde X_{t}^{\varepsilon,\delta}-\tilde X_{t(\Delta)}^{\varepsilon,\delta}|^2\big]\cr
&\le& \frac{\beta_1}{4\delta}\mathbb{E}\big[|\tilde Y_{t}^{\varepsilon,\delta}-\hat Y_{t}^{\varepsilon,\delta}|^2\big]+\frac{c_3}{\delta}\mathbb{E}\big[|\tilde X_{t}^{\varepsilon,\delta}-\tilde X_{t(\Delta)}^{\varepsilon,\delta}|^2\big]\cr
&\le& \frac{\beta_1}{4\delta}\mathbb{E}\big[|\tilde Y_{t}^{\varepsilon,\delta}-\hat Y_{t}^{\varepsilon,\delta}|^2\big]+\frac{c_3}{\delta}\Delta^{2\beta}\mathbb{E}[\| \tilde {X}^{\varepsilon,\delta} \|_{\beta-\operatorname{hld}}^2],
\end{eqnarray}
where $c_3>0$ is independent of $\varepsilon,\delta$. 

By applying (\textbf{A3}), it deduces that 
\begin{eqnarray}\label{6-6}
I_5 &\le& \frac{2ch(\varepsilon)}{{\sqrt { \delta } }}\mathbb{E}\big[\big|\tilde Y_{t}^{\varepsilon,\delta}-\hat Y_{t}^{\varepsilon,\delta}\big|\times \big|1+\tilde X_{t}^{\varepsilon,\delta}\big|\big|{\frac{dv_t^{\varepsilon,\delta}}{dt} }\big|\big]\cr
&\le& \frac{\beta_1}{4{\delta}}\mathbb{E}\big[\big|\tilde Y_{t}^{\varepsilon,\delta}-\hat Y_{t}^{\varepsilon,\delta}\big|^2\big]+{c_4h^2(\varepsilon)}\mathbb{E}\big[\big|1+\tilde X_{t}^{\varepsilon,\delta}\big|^2\big|{\frac{dv_t^{\varepsilon,\delta}}{dt} }\big|^2\big],
\end{eqnarray}
where $c_4>0$ is independent of $\varepsilon,\delta$.
Then, by combining \eqref{6-2}--\eqref{6-6}, we obtain
\begin{eqnarray}\label{6-7}
\frac{d}{dt}\mathbb{E}[{| \tilde Y_{t}^{\varepsilon,\delta}-\hat Y_{t}^{\varepsilon,\delta}|^2}]&\le&-\frac{\beta_1}{4\delta}\mathbb{E}\big[\big|\tilde Y_{t}^{\varepsilon,\delta}-\hat Y_{t}^{\varepsilon,\delta}\big|^2\big]+{c_4}{h^2(\varepsilon)}\mathbb{E}\big[\big|1+\tilde X_{t}^{\varepsilon,\delta}\big|^2\big|{\frac{dv_t^{\varepsilon,\delta}}{dt} }\big|^2\big]+\frac{c_2+c_3}{\delta}\Delta^{2\beta}\cr
&\le&-\frac{\beta_1}{4\delta}\mathbb{E}\big[\big|\tilde Y_{t}^{\varepsilon,\delta}-\hat Y_{t}^{\varepsilon,\delta}\big|^2\big]+{c_5}{h^2(\varepsilon)}\mathbb{E}\big[\big|{\frac{dv_t^{\varepsilon,\delta}}{dt} }\big|^2\big]+{c_5}{h^2(\varepsilon)}\mathbb{E}\big[\big|\tilde X_{t}^{\varepsilon,\delta}\big|^2\big|{\frac{dv_t^{\varepsilon,\delta}}{dt} }\big|^2\big]\cr
&&+\frac{c_2+c_3}{\delta}\Delta^{2\beta}
\end{eqnarray}
Moreover, consider the ODE as following:
\begin{eqnarray*}\label{6-11}
\begin{aligned}
\frac{d {A}_{t} }{dt} 
&= {\frac{{ - \beta_1 }}{4\delta } } {{A}_{t} } +{c_5}{h^2(\varepsilon)}\mathbb{E}\big[\big|{\frac{dv_t^{\varepsilon,\delta}}{dt} }\big|^2\big]+{c_5}{h^2(\varepsilon)}\mathbb{E}\big[\big|\tilde X_{t}^{\varepsilon,\delta}\big|^2\big|{\frac{dv_t^{\varepsilon,\delta}}{dt} }\big|^2\big]+\frac{c_2+c_3}{\delta}\Delta^{2\beta}
\end{aligned}
\end{eqnarray*}
with initial value $ A_0=0$. So its solution has an explicit expression as below:
\begin{eqnarray*}\label{3-234}
\begin{aligned}
{A_t}
&= {c_5}{h^2(\varepsilon)}\int_{0}^{t} e^{-\frac{\beta_1}{4\delta} (t-s)}{\mathbb{E}[{| \tilde {X}_{s}^{\varepsilon,\delta} |^2{| {\frac{dv_s^{\varepsilon, \delta}}{ds} } |^2}}]} ds+ {c_5}{h^2(\varepsilon)}\int_{0}^{t} e^{-\frac{\beta_1}{4\delta} (t-s)} { \mathbb{E}[|{\frac{dv_s^{\varepsilon, \delta}}{ds} } |^2]}ds\\
&\quad+ \frac{c_2+c_3}{\delta}\Delta^{2\beta}\int_{0}^{t} e^{-\frac{\beta_1}{4\delta} (t-s)}ds.
\end{aligned}
\end{eqnarray*}
A straightforward calculation yields
\begin{eqnarray*}\label{6-12}
\begin{aligned}
\frac{d(\mathbb{E}[{| \tilde Y_{t}^{\varepsilon,\delta}-\hat Y_{t}^{\varepsilon,\delta} |^2}]-{A_t})}{dt} 
&\le{\frac{{ - \beta_1 }}{4\delta } } \big[\mathbb{E}[{| \tilde Y_{t}^{\varepsilon,\delta}-\hat Y_{t}^{\varepsilon,\delta} |^2}]-{A_t}]\big]. 
\end{aligned}
\end{eqnarray*}
So, by comparison theorem, we have for all $t$ that
\begin{eqnarray*}\label{6-13}
\mathbb{E}[{| \tilde Y_{t}^{\varepsilon,\delta}-\hat Y_{t}^{\varepsilon,\delta} |^2}] \le A_t
&=& {c_5}{h^2(\varepsilon)}\int_{0}^{t} e^{-\frac{\beta_1}{4\delta} (t-s)}{\mathbb{E}[{| \tilde {X}_{s}^{\varepsilon,\delta} |^2{| {\frac{dv_s^{\varepsilon, \delta}}{ds} } |^2}}]} ds+ {c_5}{h^2(\varepsilon)}\int_{0}^{t} e^{-\frac{\beta_1}{4\delta} (t-s)} { \mathbb{E}[|{\frac{dv_s^{\varepsilon, \delta}}{ds} } |^2]}ds\cr
&&+ \frac{c_2+c_3}{\delta}\Delta^{2\beta}\int_{0}^{t} e^{-\frac{\beta_1}{4\delta} (t-s)}ds.
\end{eqnarray*}
Then, by leveraging the Fubini theorem, it deduces that 
\begin{eqnarray*}\label{6-14}
&&\int_{0}^{T}\mathbb{E}[{| \tilde Y_{t}^{\varepsilon,\delta}-\hat Y_{t}^{\varepsilon,\delta} |^2}]dt \cr
&\le& {c_5}{h^2(\varepsilon)}\int_0^T\int_{0}^{t} e^{-\frac{\beta_1}{4\delta} (t-s)}{\mathbb{E}[{| \tilde {X}_{s}^{\varepsilon,\delta} |^2{| {\frac{dv_s^{\varepsilon, \delta}}{ds} } |^2}}]} ds+ {c_5}{h^2(\varepsilon)}\int_0^T\int_{0}^{t} e^{-\frac{\beta_1}{4\delta} (t-s)} { \mathbb{E}[|{\frac{dv_s^{\varepsilon, \delta}}{ds} } |^2]}ds\cr
&&+ \frac{c_2+c_3}{\delta}\Delta^{2\beta}\int_0^T\int_{0}^{t} e^{-\frac{\beta_1}{4\delta} (t-s)}ds\cr
&\le& {c_5}{h^2(\varepsilon)}\mathbb{E}\big[\sup_{0 \leq t \leq T}| \tilde {X}_{t}^{\varepsilon,\delta} |^2\int_{0}^{T}\int_{s}^{T} e^{-\frac{\beta_1}{4\delta} (t-s)} dt| {\frac{dv_s^{\varepsilon, \delta}}{ds} } |^2ds\big] \\
&&+ {c_5}{h^2(\varepsilon)}\int_{0}^{T}\int_{s}^{T} e^{-\frac{\beta_1}{4\delta} (t-s)}dt { \mathbb{E}[|{\frac{dv_s^{(\varepsilon, \delta)}}{ds} } |^2]}ds + \frac{c_2+c_3}{\delta}\Delta^{2\beta}\int_{0}^{T}\int_{0}^{t} e^{-\frac{\beta_1}{4\delta} (t-s)}ds\\
&\le& \frac{4\delta {c_5}{h^2(\varepsilon)}}{\beta_1 }\mathbb{E}\big[\sup_{0 \leq t \leq T}| \tilde {x}_{t}^{\varepsilon,\delta} |^2\int_{0}^{T} e^{-\frac{\beta_1}{4\delta} (T-s)} | {\frac{dv_s^{\varepsilon, \delta}}{ds} } |^2ds\big] \\
&&+ \frac{4\delta {c_5}{h^2(\varepsilon)}}{\beta_1 }\int_{0}^{T} e^{-\frac{\beta_1}{4\delta} (T-s)} { \mathbb{E}[|{\frac{dv_s^{\varepsilon, \delta}}{ds} } |^2]}ds +c_6\Delta^{2\beta}.
\end{eqnarray*}
By the condition that $(u^{\varepsilon, \delta}, v^{\varepsilon, \delta})\in \mathcal{A}_{b}^N$ and \lemref{lem1}, we obtain that
\begin{eqnarray*}\label{6-15}
\int_{0}^{T}\mathbb{E}[{| \tilde Y_{t}^{\varepsilon,\delta}-\hat Y_{t}^{\varepsilon,\delta} |^2}]dt &\le& C_{\alpha, \beta, N}(\delta h^2(\varepsilon)+\Delta^{2\beta}).
\end{eqnarray*}
The proof is completed.
\qed

\section*{Acknowledgments}
This work was partly supported by the Key International (Regional) Cooperative Research Projects of the NSF of China (Grant 12120101002) and JSPS Grant-in-Aid for JSPS Fellows (Grant No.25KF0152). 
\section*{Declaration of competing interest}
The authors declare that they have no known competing financial interests or personal relationships that could have appeared to influence the work reported in this paper.

\end{document}